\documentclass{gtart}
\gtart

\usepackage{amscd}
\usepackage{amssymb,amsthm,amsmath}
\usepackage{latexsym}
\usepackage{enumerate}
\usepackage[all]{xy}

\newtheorem{thm}{Theorem}
\newtheorem{lem}{Lemma}
\newtheorem{prop}[lem]{Proposition}
\newtheorem{cor}[lem]{Corollary}
\newtheorem{conj}[lem]{Conjecture}

\newtheorem{clai}[lem]{Claim}
\newtheorem{expl}[lem]{Explanation}
\theoremstyle{definition}
\newtheorem{defi}[lem]{Definition}
\newtheorem{rema}[lem]{Remark}
\newtheorem{nota}[lem]{Notation}
\newtheorem{exam}[lem]{Example}

\newtheorem{add}[lem]{Addendum}

\def\Q {{\mathbb Q}}

\def\Z{{\mathbb Z}}
\def\Kaz{{\mathcal K}}
\def\BR{{\mathbb R}}

\def\Cob{{\text{\rm Cob}}}
\def\Emb{{\text{\rm Emb}}}

\def\ker{{\text{\rm ker}\,}}
\def\coker{{\text{\rm coker}\,}}
\def\Hom{{\text{\rm Hom}}}
\def\id{{\text{\rm id}}}
\def\fin{{\text{\rm fin}}}
\def\sk{{\text{\rm sk}}}
\def\wt{\widetilde}
\def\wh{\widehat}
\def\ol{\overline }

\newcommand{\defo}{\overset{\mathrm{def}}{=}}

\begin{document}

\title{Cobordism of singular maps}
\author{Andr\'as Sz\H ucs}
\address{E\"otv\"os Lor\'and University\\
P\'azm\'any P\'eter s\'et\'any 1/C, 3-206\\
1117 Budapest, HUNGARY}
\email{szucs@cs.elte.hu}

\begin{abstract}
Throughout
this paper we consider smooth maps of positive codimensions, having only
stable singularities (see \cite{Ar2},  $\S 1.4$ in Chapter 1).
We prove a conjecture due to M. Kazarian, connecting two classifying spaces in
singularity theory for this type of singular maps.. These spaces are:

1) Kazarian's space (generalising Vassiliev's algebraic complex and) showing
  which cohomology classes are represented by singularity strata.

2) The space $X_{\tau}$ giving homotopy representation of cobordisms of
singular maps with a given list of allowed singularities \cite{R--Sz},
\cite{Sz1}, \cite{Sz2}.

Our results are:\
\begin{enumerate}[a)]
\item 
We obtain that the ranks of cobordism groups of singular maps
with a given list of allowed stable singularities, and also their $p$-torsion parts
for big primes $p$
coincide with those of the homology groups of the corresponding Kazarian space.
 (A prime $p$ is ``big'' if it is greater than half of the
dimension of the source manifold.)
For all types of Morin maps (i.e.
when the list of allowed singularities contains only corank $1$ maps)
we compute these ranks explicitly.

\item 

We give a very transparent  homotopical description of the classifying space
$X_\tau$ as a fibration.
Using this fibration
we solve the problem of elimination of singularities by
cobordisms. (This is a modification of a question posed by Arnold, \cite{Ar1},
page 212.)

\end{enumerate}

\end{abstract}


\makeshorttitle

In this paper all smooth maps will have {\it {positive codimensions}} and {\it
{stable}} singularities.
The aim of the present paper is to establish a close
relationship between Kazarian's homological characteristic spectral sequence
(generalising Vassiliev's complex) and the classifying space for
cobordisms of maps having singularities only from a given fixed list(see \cite{Sz1},
\cite{R--Sz}).
Such a relationship was conjectured by M. Kazarian, and it can
be expressed as follows (in our formulation):
There is a spectrum with a filtration such that the arising
homological spectral sequence gives Kazarian's characteristic
spectral sequence while the homotopy groups of the spectrum give
the corresponding cobordism groups of singular maps (with a shift of the
dimension).
Hence the Hurewicz map for this spectrum induces a
rational isomorphism from the cobordism groups of singular maps to the
homology groups of Kazarian's space.

This allows us to extend quite a few classical theorems of cobordism theory to
cobordisms of singular maps. In particular we extend (at least modulo torsion):

\begin{enumerate}[\indent a)]
\item
  The Pontrjagin - Thom theorem claiming that the characteristic numbers
  determine the cobordism class of a manifold.
\item
  The definition of Conner and Floyd of the characteristic numbers of bordism
  classes and the statement that these numbers form a complete set of invariants.
We give a complete and computable set of obstructions to elimination of
  singularities by cobordism.

\item For any set of corank one singularities we determine explicitly the
  ranks of the corresponding  cobordism groups.

\item Give some general results beyond Morin maps, and also on the torsion
  groups.

\item Show a ``Postnikov like tower'' that produces the classifying space for 
cobordisms of singular maps as an iterated fibration. This gives a spectral
sequence starting with the cobordism groups of {\it immersions} (of the singularity
strata) and converging to the cobordism group of {\it singular} maps.  

\end{enumerate}

History: To our knowledge
cobordisms of singular maps were considered first by U. Koschorke \cite{Ko} and  the author \cite{Sz1}.
In 1984 Eliashberg in [E] proved a Pontrjagin - Thom type theorem for these cobordisms.
For a long time there were hardly any more results except the paper of Chess~\cite{Ch} 
considering a similar but different problem.   In the last decade
works of Ando \cite{A}, Saeki~\cite{S}, \cite{I--S},
Ikegami~\cite{I}, Kalm\'ar~\cite{Ka}
dealt with some cases of negative codimension maps.
The subject of Sadykov's recent preprint~\cite{Sa} is closely related to our one,
 but applies a different approach, it generalises Ando's paper and gives a
 proof for the Pontrjagin - Thom construction for singular maps, different
 from our one.

Acknowledgement: 
Endre Szab\'o spent countless hours discussing with the author the material of
the present paper. He made a significant contribution to this paper by
suggesting that the ``key bundles'' (Definition 107.) might be bundles. 
The  author gave talks on the subject of this paper once a week during a semester
on a seminar from February to May of 2005, where the
discussion with the members of the audience helped in clarifying the
ideas. My particular thanks are due to
Endre Szab\'o and  Andr\'as N\'emethi, further to G\'abor Lippner, Tam\'as Terpai 
and Boldizs\'ar Kalm\'ar.

\part{Preliminaries}
\section{Cobordism groups}
\begin{defi}\label{susp}
By {\it {singularity}} we mean an equivalence class of germs $(\BR^c, 0)
\longrightarrow (\BR^{c + k}, 0)$ of fixed codimension $k > 0$, where
the equivalence is generated by $\mathcal A$-equivalence
(=~left-right equivalence) and suspension:
\footnote {also called {\it {trivial extension}}, see \cite{K3}.}
\[
\Sigma \eta(u, t) = \bigl(\eta(u), t\bigr),
\]
i.e.\ $\eta : (\BR^c, 0) \to (\BR^{c + k}, 0)$ is equivalent to
\[
\Sigma \eta = \eta \times \text{\rm id}_{\BR^1} : (\BR^c \times \BR^1, 0)
\to (\BR^{c + k} \times \BR^1, 0).
\]
For a germ $\eta$ we denote by $[\eta]$ its equivalence class,
and will call it also the singularity class of $\eta$.
We say that the germ  $\eta : (\BR^c, 0) \to (\BR^{c + k}, 0)$ has an isolated
singularity at the origin if at no point near the origin the germ of $\eta$ is
equivalent to that at the origin.
Such an $\eta$ will be called the {\it {root}} of its singularity
class~$[\eta]$ and will be denoted by $\eta^0$. The root of $[\eta]$ is characterised also by having the
smallest possible dimension $c$ in its class.
Given a smooth map $f: M^n \to P^{n + k}$ we say that $x \in
M^n$ is an $[\eta]$-point (or simply an $\eta$-point) if the
germ of $f$ at $x$ belongs to~$[\eta]$. The set of such points in $M^n$ will
be denoted by $\eta(f).$ This is submanifold (possibly non closed subset) in
$M^n$
of dimension $n-c.$ By this reason $c$ will be called the {\it {codimension}} of $\eta.$
\end{defi}

In this paper we shall always consider stable singularities (see \cite{Ar2}),
of maps of a fixed positive codimension $k$. For these type of singularities
there is a natural hierarchy.
We say that $\eta$ is more complicated than $\xi$, and write
$\xi < \eta$, if in any neighbourhood of an $\eta$-point there
is a $\xi$-point.

\begin{exam}\

\begin{enumerate}[1)]
 \item 
  The class of germs of regular (i.e. non-singular) maps $\Sigma^0$
 is the lowest ``singularity.''
\item
 The  class of a  fold germ (i.e. a stable 
 $\Sigma^{1,0}$\ germ) is lower than that of a cusp germ (i.e. of the stable
 $\Sigma^{1,1}$\ germ).

\end{enumerate}

\end{exam}
\begin{defi}\
\begin{enumerate}[\indent a)]
\item Let $\tau$ be a set  of stable singularity classes of codimension~$k > 0$.
Given two smooth manifolds $M^n$ and $P^{n+k}$ of dimensions $n$ and $n+k$ respectively
we say that
a proper smooth map $f: M^n \to P^{n+k}$ is a
{\it {$\tau$-map}} if at any point $x \in M^n$ the germ of $f$ 
belongs to a singularity class from~$\tau$.

\item If $M^n$ has non-empty boundary, then $P^{n+k}$ must have it too, and a $\tau$-map
$f : M^n \to P^{n+k}$ must map $\partial M^n$ into $\partial P^{n+k}$, moreover
$\partial M^n$ and $\partial P^{n+k}$ have collar neighbourhoods
$\partial M^n \times [0, {\varepsilon})$ and $\partial P^{n+k} \times [0,
{\varepsilon})$ such that $f\big|_{\partial M^n \times [0,
{\varepsilon})}$ has the form $g \times \id_{[0,{\varepsilon})}$,
where $g : \partial M^n \to \partial P^{n+k}$ is a $\tau$-map, $\id_{[0,
{\varepsilon})}$ is the identity
map of~$[0, {\varepsilon})$.
\end{enumerate}

\end{defi}
Next we define 
manifolds with corners and then 
 $\tau$-maps when both the source and the target
are manifolds with corners.

\begin{defi}
We extend the definition of a smooth manifold with boundary to
the definition of a {\it smooth manifold with corners}  by
allowing neighbourhoods of the type 
$\BR^{n-m} \times [0, {\varepsilon})^{m}$ for any $m \ge 0.$
The numbers $n$ and $m$ will be called the dimension of the manifold and the
codimension of the corner respectively.
The positive codimension part (where $m > 0$) will be called {\it {the
    boundary}}, its complement is the {\it {interior part}} of the manifold
with corners. These parts will be denoted by $\partial$ and $
int$ respectively.
Note that the boundary itself is an $n-1$-dimensional manifold with corners.

Now if $M^n$ and $P^{n+k}$ are manifolds with corners then 
a proper smooth map $f: M^n \to P^{n+k}$ is a
{\it $\tau$-map} if a neighbourhood of
the corner of codimension $m$ of $M$ is mapped into that of $P^{n+k}$ by a  map
of the local form $g \times \id_{[0,{\varepsilon})^m}$, where
$g: \BR^{n-m} \to \BR^{n+k-m}$ is a $\tau$-map.
So near a corner of codimension $m$ the map looks like an $m$-tuple suspension
of a $\tau$-map, see Definition~\ref{susp}.
\end{defi}

\begin{defi}
  For an arbitrary (possibly non-compact) oriented smooth 
  $(n~+~k)$-dimensional manifold~$P^{n + k}$ 
we will denote by  $\text{\rm
    Cob}_\tau(P^{n + k})$ the set of oriented
  $\tau$-cobordism classes of $\tau$-maps of
  $n$-dimensional manifolds in $P^{n + k}$.
  Namely: Given two $\tau$-maps $f_0 : M^n_0 \to P^{n + k}$ and
  $f_1 : M^n_1 \to P^{n + k}$ (where $M^n_0$, $M^n_1$ are closed,
  oriented, smooth $n$-manifolds) they are said to be
  $\tau$-cobordant if there exist
  \begin{enumerate}[\indent 1)]
  \item
    a compact, oriented $(n + 1)$-manifold $W^{n + 1}$ such that 
    
    $\partial
    W^{n+1} = -M^n_0 \sqcup M^n_1$, and
  \item
    a $\tau$-map $F : W^{n + 1} \to P^{n+k} \times [0,1]$ s.t.
    \[
    F^{-1} \bigl(P^{n+k} \times \{ i \}\bigr) = M^n_i \quad \text{ and }
    \quad
    F\big|_{M^n_i} = f_i \quad \text{for } \ i = 0,1.
    \]
  \end{enumerate}
  The $\tau$-cobordism class of a $\tau$-map $f$ will be denoted
  by $[f].$
\end{defi}

Note that for any two $n$-manifolds $N^n_i$,\ $i = 0, 1$, and any two $\tau$-maps $g_i : N^n_i \to P^{n + k}$ 
the disjoint union $g_0 \sqcup g_1 : N^n_0 \sqcup N^n_1 \to P^{n +
k}$ is also a $\tau$-map.
(The images of $g_0$ and $g_1$ may intersect.)
The disjoint union defines an associative and commutative
operation on $\Cob_\tau(P^{n+k})$ with a null-element represented
by the map of the empty set.

\begin{rema}
In this paper we shall deal exclusively with oriented cobordisms of singular
maps and so - typically - we shall not indicate that in our notation.
(Each construction, definition, theorem of our paper has its unoriented counterpart, but we
do not include them here.)
\end{rema}

\begin{defi}\label{simpl}

The notion $\Cob_\tau(P)$ of cobordism set of $\tau$-maps into $P$ can be
extended to the case when $P$ is a finite simplicial complex.
Let $P$  consist of a single simplex first.
Note that a simplex is a manifold with corners in a natural way and also a
cylinder over a simplex is that.
Since we have defined the $\tau$-maps for manifolds with corners, the
definition of $\Cob_\tau(P)$ can be given just repeating the definition for
the manifold case (replacing the word ``manifold'' by ``manifold with corners''). 
When $P$ consists of several simplices, then 
a $\tau$-map into $P$ is a set of $\tau$-maps into each simplex of (possibly
empty) manifolds
with corners, and with natural identifications of the corners and their maps
corresponding to the identifications of the faces of the simplices.
(Note that into each open simplex of dimension $i$ we map manifolds of
dimension $i-k$.)
The cobordism will be a $\tau$-map into the cylinder $P\times [0,1].$ 
\end{defi}

\begin{rema}
It is not obvious but true that for any finite simplicial complex $P$ the set
 $\Cob_\tau(P)$  with
this operation forms a group, i.e.\ there is an
inverse element. 
This fact follows  of course from the homotopy representation of this set
- the generalised Pontrjagin - Thom construction - proved in a number of
papers
\cite{R--Sz}, \cite{A}, \cite{Sa}, see also the present paper, but
here we give an explicit geometric description of the inverse element.

\end{rema}

\subsection{On the inverse element}

\begin{enumerate}[\bf a)]
\item
  If $P^{n+k} = \BR^{n + k}$, then composing a map $f : M^n \to
  \BR^{n + k}$ with a reflection $r$ in a hyperplane of 
$\BR^{n + k}$
 and an orientation
  reversing diffeomorphism $T:M^n \to M^n$ of the one
  gets the inverse element: $-[f] = [r \circ f \circ T]$.

\item
  If $P^{n + k}$ has the form $N^{n + k - 1} \times \BR^1$,
  then still one can reflect in a ``hyperplane'', i.e.\ in a
  submanifold $N^{n+k-1} \times \{c\}$, $c \in \BR^1$ and composed it with the
  orientation reversing map $T$ of the source we get the inverse.
\item
  For an arbitrary manifold $P^{n + k}$ we need the notion of ``framed
  $\tau$-map'':
\end{enumerate}
\begin{defi}
A {\it {framed $\tau$-map}} from an $n$-manifold $M^n$ into an
$(n + k + 1)$-manifold $Q^{n + k + 1}$ is a map germ along
$M^n \times \{0\} \approx M^n$ of a $\tau$-map
$M^n \times \BR^1 \supset
\mathcal N^{n+1} \overset{\wt f}{\longrightarrow} Q^{n + k + 1}$
defined in a neighbourhood $\mathcal N^{n+1}$ of $M^n \times \{0\}$, such that
around each point $p = (x, 0) \in M^{n} \times \{0\}$ and its image
$q = \wt f(p)$ there exist coordinate neighbourhoods
$W^{n+1}_x, \wt W_q^{n+k+1}$ such that
 $W_x^{n+1} = U_x
\times V_0$, where $U_x \subset M^n$, $V_0 \subset \BR^1$, $U_x \approx
\BR^n$, $V_0 \approx \BR^1$,  and  $\wt W_q^{n+k+1} \approx \BR^{n + k + 1} =
\BR^{n + k} \times \BR^1$ in which $\wt f$ has the form: $g \times
\id_{\BR^1}$, where $g : \BR^n \to \BR^{n + k}$ is a $\tau$-map and 
$\text{id}_{\BR^1}$ is the identity map $\BR^1 \to \BR^1.$
(Note
that the $\BR^1$-direction in the source is always directed along
the second factor of $M^n \times \BR^1$, while in the target there is
no such a fixed direction.)
\end{defi}
\begin{rema}\
  \begin{enumerate}[1)]
  \item
    The framed $\tau$-map is an extension of the notion of
    ``immersion with a normal vector field'' to the case of maps
    with singularities.
  \item
    Analogously it can be defined the notion of $\ell$-framed
    $\tau$-maps as germs at $M^n \times 0$ of maps $M^n \times \BR^\ell
    \supset \mathcal N^{n+\ell} \to Q^{n + k + \ell}$ etc.
  \end{enumerate}
\end{rema}
\begin{exam}
If $f : M^n \to P^{n + k}$ is a $\tau$-map, then the germ of the
map $f \times \id_{\BR^1} : M^n \times \BR^1 \to P^{n+k} \times \BR^1$ along $M^n
\times \{0\}$ gives a framed $\tau$-map of $M^n$ into $P^{n + k}
\times \BR^1$.
\end{exam}
\begin{defi}\label{framed}
The {\it {cobordism set of framed $\tau$-maps}} of closed, (oriented) smooth
$n$-mani\-folds into an (oriented) smooth $(n + k + 1)$-dimensional manifold $Q^{n + k + 1}$ can be defined by an
obvious modification of the definition of $\tau$-cobordism.
The set of cobordism classes will be denoted by
$\Cob_{\tau\oplus 1} (Q^{n+k+1}).$ The cobordism set of $\ell$-framed $\tau$-maps will be denoted
by $\Cob_{\tau\oplus \ell} (Q^{n+k+1}).$
\end{defi}

\begin{prop}
  \label{cob-tau-n=cob-tau-n+l}
Given  oriented $n$ and $n+k$-dimensional manifolds $M^n$ and  $P^{n+k}$
associating to
a $\tau$-map $f: M^n \to P^{n+k}$  the $\ell$- framed $\tau$-map
  $f \times \id_{\BR^{\ell}}$ one obtains a one-to-one
  correspondence
  $$\Cob_\tau(P^{n + k})
  \overset{1\text{-}1}{\longrightarrow} \Cob_{\tau \oplus \ell}
  (P^{n+k}\times \BR^\ell)$$
  for any natural number $\ell$. (Recall, that $k > 0$.)
\end{prop}

\noindent
\emph{Proof} will be given later, see Remark~\ref{ell-framed}.
It will follow from a stratified version of the so-called
``Compression Theorem'' (see \cite{Ro--Sa} and also
\cite{Gr}). The manifold $P^{n+k}$ can be replaced by any finite 
simplicial complex. \hfill $\square$
\medskip

Now for an element in $\Cob_\tau(P^{n + k})$ its inverse can be
obtained using the one - to - one correspondence in the Proposition, since
in $\Cob_{\tau\oplus 1}(P^{n+k} \times \BR^1)$ we have a target manifold
of the type~b).
More explicitly for a $\tau$-map $f : M^n \to P^{n + k}$ a
representative of the inverse cobordism class $-[f]$ can be
obtained as follows:
Consider the germ of $f \times  \id_{\BR^1} : M^n \times \BR^1 \to P^{n+k}
\times \BR^1$ at $M^n \times 0$. This can be visualised as the map
$\wt f : M^n \overset{f}{\longrightarrow} P^{n+k} \hookrightarrow P^{n+k} \times
\BR^1$ with an ``upward'' directed vector field~$\uparrow$. (The
lines $\{p\} \times \BR^1 \subset P^{n+k} \times \BR^1, \ p \in P^{n+k},$ with the positive direction
 we call upward directed
lines.) Now according to b) we reflect this map $\wt f$ in the
``hyperplane'' $P^{n+k} \times \{0\}$ in order to obtain a representative of the
inverse cobordism class $-[\wt f]$. Since the image of $\wt f$ lies in
$P^{n+k}\times \{0\}$, this reduces to taking the same
map $\wt f$ with a downward directed vector field~$\downarrow$.
The `` Stratified Compression Theorem'' (see Theorem~\ref{SCT}) will claim that
there is an isotopy $\varphi_t$ of
$P^{n+k} \times \BR^1$
such that $\varphi_0 $ is the identity map, and the differential $d(\varphi_1)$
maps the downward directed vector field along $\wt f$ into the upward directed
one along $\varphi_1 \circ \wt f.$
 Now $\varphi_1 \circ \wt f \circ T  : M^n \to P^{n+k}
\times \BR^1$ can be composed with a projection to $P^{n+k}$, and this gives a $\tau$-map
$g : M^n \to P^{n+k}$ such that $[g] = -[f]$ in $\Cob_\tau(
P^{n+k})$.

Finally if $P$ is an arbitrary finite simplicial complex, then we will use
the relative version of the Stratified Compression Theorem to get an 
``isotopy'' of $P$
moving each simplex within itself and turning the downward directed
vector-field into an upward directed one.

Another way of reducing the case when $P$ is a finite simplicial complex to
that when it is a manifold is the following:

A finite simplicial complex is the deformation retract of its very small
neighbourhood in a big dimensional Euclidean space. This neighbourhood can
be chosen so that its closure be a smooth manifold with boundary.
Now having a $\tau$-map in the simplicial complex $P$ it extends to a
$\tau$-map into a small closed neighbourhood $U$ of $P$ in the Euclidean
space, where $U$ is a manifold with boundary. Applying the procedure 
constructing the inverse element in $\Cob_\tau(U)$ and then restricting this
element to $P$ we obtain the inverse element of the class of the original 
$\tau$-map.

In order to show that $\Cob_\tau(\ )$ is a contravariant functor we have to
define the induced maps. This will be done if we show how to  pull
back a   $\tau$-map.

\begin{prop}\label{pullback} (Pulling back a $\tau$-map.)

$\Cob_\tau(P)$ is a contravariant functor on the category of finite simplicial
complexes and their simplicial maps. (Actually it can be extended to their
continuous maps.)
\end{prop}

\begin{proof}
\begin{enumerate}[1)]\
\item
Let $P^{n+k}$ and $P^{n'+k}$ be first smooth manifolds of dimensions $n+k$ and $n'+k$ 
respectively.
Let $g: P^{n'+k} \to P^{n+k}$ be a continuous map. We show that $g$ induces a well-defined
homomorphism  $g^*: \Cob_\tau(P^{n+k}) \to \Cob_\tau(P^{n'+k}).$ Let $f: M^n \to P^{n+k}$
be a $\tau$-map. The map $g$ can be approximated by a smooth map (homotopic to
$g$) transverse to the image of each singularity stratum  of $f$, we will denote
this transverse map again by $g.$ Let us define the ``pull-back'' map 
$g^*(f)$ as follows.
Consider the
 map $\wt f: \wt M^n \to P^{n'+k}$, where $$\wt M^n = \{(p', m)\in P^{n'+k}\times M^n |
    g(p') = f(m)\},\ \ \wt f (p',m) = p'$$ and put $g^*(f) = \wt f.$
    The map $g^*(f)$
    will be a $\tau$-map. (This was formulated without proof in \cite{K3},
    page 71. The proof of this statement is contained implicitly in the 
    paper \cite{R--Sz}. Below we shall prove it.)
    Further the cobordism class of the obtained map $\wt f$ will be shown to
    be  
    independent
    of the chosen transverse approximation of $g.$ Moreover we show that
    homotopic maps induce the same homomorphisms, i.e $g \cong h
    \Longrightarrow g^* = h^*.$
   
\begin{clai}\ 

\begin{enumerate}[\indent a)]
\item
$\wt M^n$ is a smooth manifold.
\item 
$\wt f$ is a $\tau$-map.
\end{enumerate}  
\end{clai}

\begin{proof}
It was shown in the proof of Theorem 1 in \cite{R--Sz} (see also 
\cite{Sz8} and \cite{Sz9} for Morin maps) that in a neighbourhood of the $\eta$-stratum the map
$f$ has local form
$$  
\eta \times \id_{\BR^u}: \BR^c\times \BR^u \to \BR^{c+k}\times \BR^u
$$

Let us call these local coordinates {\it canonical.} 

Here the subspaces $\{0\} \times R^u$ are directed along the singularity 
stratum and its image, while the complementary subspaces $\BR^c \times\{0\}$
and   
$\BR^{c+k} \times\{0\}$ are directed along  the normal fibres of the stratum
$\eta(f)$ in the source and that of its image $\wt \eta(f) = f(\eta(f))$ in the target
respectively.
Since the map $g$ is transverse to the stratum $\wt \eta(f)$ for any point 
$p'\in P^{n'+k}$ such that $g(p') = f(x),$ where  $x \in \eta(f),$ there is a local coordinate
system $\BR^{c+k}\times\BR^v$ centred around $p'$, such that 
in these local coordinates and in the canonical local coordinates around
$g(p')$ the map $g$ has the form
$\id_{\BR^{c+k}} \times h: \BR^{c+k}\times\BR^v \to \BR^{c+k}\times\BR^u $,
where $h: \BR^v \to \BR^u$ is a smooth map.
Now it is clear that in the given local coordinates the map $\wt f = g^*(f):
\BR^c \times \BR^v \to \BR^c \times \BR^u$
has the form   $\id_{\BR^c}\times h.$
In particular $\wt M^n$ is a smooth manifold (because the arbitrarily chosen point $(x,p')$ 
has a Euclidean neighbourhood and the change of coordinates is smooth
since they are such in $M^n$ and $P^{n'+k}$). Further $\wt f$ has in these local
coordinates a normal form equivalent to $\eta$ around the point $(x, p').$
Hence the pulled back map is a $\tau$-map.
\end{proof}

\item
All this works also when $P^{n+k}$ and $P^{n'+k}$ are not manifolds but those with
corners. Then the map $g$ must map corners of codimension $r$ in $P^{n'+k}$ into those in
$P^{n+k}.$ (Of course in this case $\wt M^n$ will be a manifold with corners.)

\item
Finally let $P$ and $P'$ be any finite simplicial complexes. We have defined 
the groups $\Cob_\tau(P)$ and $\Cob_\tau(P')$ for this case too. Recall that
each simplex is considered as a manifold with corners, and
 $\tau$-maps in this
case are those into each simplex in a compatible way.  
Each simplicial map can be considered  as a map of manifolds with corners. Hence 
$g^*: \Cob_\tau(P) \to \Cob_\tau(P')$ is defined.

\item
If $P$ is a simplicial complex and $P_1$ is its subdivision, then there are
natural isomorphisms
$\Cob_\tau(P) \to \Cob_\tau(P_1)$ and $\Cob_\tau(P_1) \to \Cob_\tau(P).$

The first is obtained by taking transversal intersections of the image of the
$\tau$-maps with the new faces in $P_1.$
The second map is obtained by noticing that at the points of the omitted faces
the $\tau$-maps can be attached along their corners and get an interior point
of the image.
\item
If $g$ is any continuous map of the simplicial complex $P'$ into the
simplicial complex $P$, then we can choose an appropriate subdivision of $P'$
and replace $g$ by a homotopic simplicial map.
\item
For homotopic maps $g \cong h$ from $P'$ to $P$ the induced maps $g^*$ and $h^*$
from
$ \Cob_\tau(P) \to \Cob_\tau(P')$ will coincide.

Indeed, given a $\tau$-map in $P$ the homotopy between $g$ and $h$ also can be made transverse to each stratum $\wt \eta(f).$
Pulling back the map $f$ by this transverse homotopy  we obtain a
$\tau$-cobordism between the  two pull-backs $g^*(f)$ and $h^*(f).$

\end{enumerate}
\end{proof}

\begin{cor}

Given two finite simplicial complexes $P$ and $P'$ denote by $[P',P]$
the set of homotopy classes of continuous maps $P' \to P.$
Then a map arises
$$
[P',P] \times \Cob_\tau(P) \to \Cob_\tau(P')$$

that is additive in its second factor.
\end{cor}

\begin{rema}\ \label{deloop}
  \begin{enumerate}[a)]

  \item
    Above we described the extension of the group-valued functor $\Cob_\tau(-)$ from smooth
    manifolds to (finite) simplicial complexes. 
    This extended functor satisfies the conditions of the Brown
    representability theorem.
    Hence there is a space $X_\tau$ such that: $\Cob_\tau(P) = [P,
    X_\tau]$ for any finite simplicial complex~$P$, see \cite{Sz3}.
    (Here the brackets $[\ ,\ ]$ denote the set of homotopy
    classes. The space $X_\tau$ has been constructed in\cite{R--Sz}
    explicitly for more general types of maps.)
 \footnote {The definition of $\tau$-maps in \cite{R--Sz} was different from
  that here. It included both local 
restrictions (as here) and global ones (these are not considered here).
In the presence of global restrictions the set of $\tau$-maps does not form a
group, the classifying space is not an $H$-space and finally there is no such
a nice connection with the  Kazarian space that we are going to prove here. In general the
situation is much more difficult in the presence of global restrictions.
We can illustrate this by comparing the classification of embeddings with that
of immersions (these are actually special cases of maps with global
restrictions and maps with only local restrictions.) While immersion theory is
completely reduced to algebraic topology, embedding theory is  
hopeless in general.} 
    
\item
    Since $\Cob_\tau(P)$ is actually a group, the space $X_\tau$
    will be an $H$-group according to an addendum to Brown's theorem.

  \item
    The space $X_\tau$ is a loop space, i.e.\  there is a space
    $Z_\tau$ such that $X_\tau$ is homotopy equivalent to $\Omega Z_\tau$. 
    Indeed, if we denote by $Z_\tau$ the classifying space of framed
    $\tau$-maps, then for any finite simplicial complex $P$
    the group of framed $\tau$-maps into the suspension $SP$ is
    isomorphic to the group $\Cob_\tau(P)$
    (by Proposition~\ref{cob-tau-n=cob-tau-n+l}).
    Hence
    \[
    [P, X_\tau] = [SP, Z_\tau] \quad \text{for any }\ P.
    \]
    On the other hand $[SP, Z_\tau] = [P, \Omega Z_\tau]$ and so
    $X_\tau \cong \Omega Z_\tau$
    since the representing space of a functor is homotopically unique.
    (Note that this shows again that $X_\tau$ is an $H$-group - since it is a
    loop-space - and
    $\Cob_\tau(P)$ is a group.)
  \item
    The above procedure can be iterated and we obtain that
    $X_\tau$ is an infinite loop space, i.e.\ for each $\ell$ there
    is a space $Z^\ell_\tau$ such that $X_\tau = \Omega^\ell Z^\ell_\tau$.
    Here $Z^\ell_\tau$ will be the classifying space of
    $\ell$-framed $\tau$-maps.
  \end{enumerate}
\end{rema}
In the introduction we promised to establish a relationship
between the classifying space of cobordism of singular maps with
prescribed singularities -- this is the space $X_\tau$ -- and
Kazarian's space (which we will denote by $\Kaz_\tau$).
In the next section we recall the two definitions of this space
$\Kaz_\tau.$ The first is  described by Kazarian in \cite{K1}. (The referee informed me, that
it is due to Thom, \cite{HK}.)
The second was communicated to me - without proof - by Kazarian \cite{K2}. 
We show the equivalence of these two definitions, this equivalence will be
used in the concrete computations.

\part{Kazarian's space}
\section{Definition derived from the Borel construction.}
\subsection{Unstable version}

Here we construct a ``preliminary'', ``unstable'' version of Kazarian's space.\footnote {``unstable'' means that in
this subsection we do not consider yet the
identification of the germs with their suspensions, 
see Definition~\ref{susp}.}
Let $\tau$ be a set of stable singularity classes such that each of
them is $K$-determined for a big enough number $K$.
Put $V_n = J^K_0(\BR^n, \BR^{n + k})$ the space of $K$-jets of germs
 $(\BR^n, 0) \to (\BR^{n + k}, 0)$ at the origin.
Let us denote by $G_n$ the group $J^K_0\,\bigl(\text{\rm Diff}^+(\BR^n,
0)\bigr) \times J^K_0\, \bigl(\text{\rm Diff}^+(\BR^{n + k}, 0)\bigr)$
i.e.\ the group of $K$-jets of orientation preserving
diffeo\-morphism-germs of $\BR^n$ and $\BR^{n + k}$ at the origin.
Those $K$-jets in $V_n$ which belong to the same singularity
class, form a $G_n$-orbit.
In particular the $K$-jets in the singularity classes lying in
$\tau$ form a $G_n$-invariant subset in $V_n$, which we denote
by $V^n_\tau$.
Put $\Kaz_\tau^{\prime}(n) = V^n_\tau \underset{G_n}{\times} E G_n$.
Let us denote by $\pi'$ the fibration $\pi' : \Kaz_\tau^{\prime} (n)
\overset{V^n_\tau}{\longrightarrow} BG_n$.
If $[\eta] \in \tau$ is a singularity class from $\tau$, then we
denote by $B\eta(n)$ the corresponding subset in $\Kaz_\tau^{\prime
}(n)$.
(For any fibre $V^n_\tau$ the intersection $V^n_\tau \cap B
\eta(n)$ coincides with the set  of germs $\mathcal A_+$-equivalent to
$\eta.$ Here $\mathcal A_+$-equivalent means $\mathcal A$-equivalent through
orientation preserving diffeomorphisms.)

\begin{lem}[\cite{HK}, \cite{K1}]\
\label{Kaz}
  \begin{enumerate}[a)]
  \item
    Let $f: M^n \to P^{n + k}$ be any $\tau$-map.
    Let us denote by $\eta(f)$ the subset of $M^n$, where the germ
    of $f$ is in~$[\eta]$.
    Then there is a map $\kappa_f' : M^n \to \Kaz_\tau^{\prime} (n)$ such that
    $(\kappa_f')^{-1}(B\eta(n)) = \eta(f)$ for any $[\eta] \in \tau.$
  \item
    Moreover this map $\kappa_f'$ has the following extra property.
    First we note that the space $BG_n$ is homotopy equivalent to
    $B\,SO_n \times B\,SO_{n + k}$. Let us lift the universal vector bundles over 
$B\,SO_n$ and $B\,SO_{n+k}$ to $BG_n.$
    The extra property of the map $\kappa_f'$ is that, 
    the composition $\pi' \circ \kappa_f'$ pulls back the bundles $TM^n$
    and $f^* TP^{n+k}$ from these lifted bundles 
    respectively, where $TM^n$ and $TP^{n+k}$ denote the tangent
    bundles over $M^n$ and $P^{n+k}$, respectively.
  \end{enumerate}
\end{lem}
\begin {rema}
Actually Kazarian has constructed this map
$\kappa_f' : M^n \to \Kaz_\tau^{\prime} (n)$, we recall his construction in
the sketch of the proof below. By the Kazarian map $\kappa_f'$ we shall mean
this particular construction (or its stabilised version described later).
\end{rema}

\begin {rema}

The space $\Kaz^{\prime}_\tau (n)$ is not yet the final definition of the
Kazarian space. It still depends on the dimension $n.$ The final definition
will be obtained in section~\ref{stab}, after some stabilisation.
\end{rema} 

\begin{defi}
\label{n-normal}

Pull back the universal bundles $\gamma^{SO}_n \to B\,SO_n$ and
$\gamma^{SO}_{n+k} \to B\,SO_{n+k}$ to the space $\Kaz^{\prime}_\tau (n)$ and consider
their formal difference as a virtual $k$-dimensional vector bundle. This
bundle will be denoted by $\nu^k_\tau(n)$ (or sometimes simply by $\nu(n)$)
and called the {\it {$n$-universal virtual normal bundle for $\tau$-maps}}.
\end{defi}

\begin{proof}[Sketch of the proof, \cite{K1}]
Let us fix Riemannian metrics on $M^n$ and $P^{n+k}$ and denote by $T^{\varepsilon} M^n$
the set of tangent vectors of length less than $\varepsilon.$
There is a fiberwise map $\widehat f : T^{\varepsilon} M^n \to f^* TP^{n+k}$
defined on a neighbourhood of the zero section of $TM^n$ such that
$\exp_P^{n+k} \circ \widehat f \circ \exp^{-1}_M = f$, where $\exp$
denotes the exponential map.
Let $V^n_\tau(M^n)$ denote the subspace of $J^K_0 \, (TM^n, f^* TP^{n+k})
\to M^n$ corresponding to the singularities in~$\tau$.
So $V^n_\tau(M^n)$ is the total space of a fibration over $M^n$
with fibre $V^n_\tau$ and structure group~$G_n$.
The map $\widehat f$ defines a section $\wt f$ of this
bundle.
The bundle $V^n_\tau(M^n) \overset{V^n_\tau}{\longrightarrow} M^n$
can be induced from the bundle $K'_\tau(n)
\overset{V^n_\tau}{\longrightarrow} BG_n$.
The obtained map $V^n_\tau(M^n) \longrightarrow K'_\tau(n)$
composed with the section~$\wt f$ gives the map
$\kappa_f' : M^n \to \Kaz_\tau^{\prime} (n)$.
Clearly $(\kappa_f')^{-1}(B\eta(n)) = \eta(f)$.
\end{proof}

\subsection{Addendum to the Borel construction.}

Let $G$ be a compact group Lie, let $V$ be a contractible smooth manifold with smooth (left) $G$-action 
on it, and let $BV= EG\times_G V$ be the Borel construction applied to $V.$
Further let $\Sigma$ be an orbit
of $G$ in $V$, let $x$ be a point in $\Sigma$, and let $G_x$ be the stabiliser of $x.$
(Hence $\Sigma = G/G_x.$) Finally let $S_x$ be a small transverse slice of $\Sigma.$
We can choose a Riemannian metric on $V$ such that the $G$-action will be
isometric. Choose $S_x$ orthogonal to $\Sigma.$ Then $G_x$ acts on the tangent
space $TS_x$ orthogonally. Let us denote by $\rho_x$ the arising representation of $G_x$
on $TS_x.$

\begin{lem}\label{n-Kaz}\

\begin{enumerate}[\indent a)]

\item 
The subset $B\Sigma = EG\times_G \Sigma$ in $BV = EG\times_G V$ is the
classifying space $BG_x$ of the stabiliser group $G_x.$ 
\item 
A neighbourhood of $B\Sigma$ in $BV$ can be identified with the
universal bundle over $BG_x$ associated with the representation $\rho_x.$
\end{enumerate}
\end{lem}

\begin{proof}\

\begin{enumerate}[\indent a)]

\item 
This part was proved by Kazarian \cite{K1},\cite{K3}. Below our proof for part b) will
show it too. 

\item 
Let $\mathcal N$ be a tubular $G$-invariant neighbourhood of $\Sigma$ (in $V$)
fibred over $\Sigma$ with fibre $S_{gx}$ over $gx \in \Sigma.$
Note that for any $g \in G$ we obtain a bijection from $S_x$ to $S_{gx}.$
Nevertheless the fibration $\mathcal N \to \Sigma$ might be non-trivial
since different elements $g, g'$ of $G$ might give different bijections onto the same
fibre $S_{gx}$ if $gx = g'x.$

Let us lift the bundle $\mathcal N \to \Sigma = G/G_x$ to $G$ pulling it
back by the projection $G \to G/G_x.$ Let $\wt {\mathcal N}$ be the obtained bundle over $G.$
We claim that the bundle $\mathcal N \to G$ is trivial.
Indeed, let $\wt S_g$ be the fibre over $g.$ Then by the definition of the
pull-back bundle there is an identification of $\wt S_g$ with $S_{[g]}$,
where $[g] \in G/G_x$ is the coset of $g.$ ($[g] \in G/G_x$ is identified with
the point $gx \in \Sigma.$) Since each element $g \in G$
gives a bijection $S_x \to S_{gx}$, we obtain a canonical identification 
$\wt S_e \to \wt S_g$ of the fibre over the neutral element $e \in G$
with the fibre over  $g \in G.$
So we get a $G$-action on $\wt {\mathcal N}$ and the obtained decomposition
$\wt {\mathcal N} = G \times \wt S_e$ is $G$-equivariant (with trivial
$G$-action on $\wt S_e.$)

The fibration $\wt {\mathcal N} \to \mathcal N$ with fibre $G_x$
is  $G$-equivariant. The well-known Borel construction is a functor from the
category of $G$-spaces and $G$-equivariant maps to the category of topological
spaces and their maps. Applying the Borel construction to the equivariant map
$\wt {\mathcal N} \to \mathcal N$
we obtain a map 
$\wt {\mathcal N} \times _G EG \to \mathcal N \times _G~EG.$
Here the domain is $ \wt {\mathcal N} \times _G EG = G \times _G \wt S_e
\times EG = \wt S_e \times EG.$
The target  $\mathcal N \times _G EG$ is a neighbourhood of $B\Sigma$ in $BV.$
The map itself is the quotient map factorising out by the diagonal
$G_x$-action on $\wt S_e \times EG.$ 
Hence the neighbourhood $B\cal N$ of $B\Sigma$ in $BV$ can be identified with 
$\wt S_e \times_{G_x} EG$, where $\wt S_e$ is a ball, and $G_x$ acts on it by
the representation $\rho_x.$ 

\end{enumerate}

\end{proof}

Later (in Lemma \ref{transvslice}) we shall need the following simple Remark.
\begin{rema}\label{slice}

Let $V, \Sigma, G$ be as above. Let $F\subset V$ be transverse to $\Sigma.$
We identify $V$ with a fibre of the Borel construction $BV.$
Then $F$ is transverse to $B\Sigma$ in $BV.$
\end{rema}

\subsection{Stabilisation}\label{stab}
Until this point in the construction of the space $K'_\tau(n)$ we
did not use that forming the equivalence classes of germs we
identified each germ with its suspension. Now we are going to use
it.
Note that adding an arbitrary smooth vector bundle $F$ over
$M^n$ to both $TM^n$ and $f^* TP^{n+k}$, and replacing the map 
$\widehat f: T^{\varepsilon}M \to f^*(TP)$
by $\widehat f \oplus \id_F$, where $\id_F$ is the identity map of $F,$
we obtain the same stratification on $M^n$ by singularity
types. In other words at $x \in M^n$ the $K$-jet of $(\widehat f
\oplus \id_F)_x : T^{\varepsilon}_x M^n \oplus F_x \to T_{f(x)} P^{n+k} \oplus F_x$ is
equivalent to the $(\dim F_x)$-tuple suspension of the $K$-jet
of $\widehat f : T_x M^n \to T_{f(x)} P^{n+k}$. 
For $x\in M^n$ the following equivalences hold:
 $x \in \eta(f)
\Leftrightarrow x \in \eta(\widehat f) \Leftrightarrow x \in
\eta(\widehat f\oplus \id_F)$.
Now choose the bundle $F$ to be equal to the (stable) normal
bundle of $M^n$ in $\BR^{N}$ for a big enough $N \ (N \ge 2n+2$: $F = \nu_M$.
 $TM^n \oplus \nu_M \approx {\varepsilon}^N$).
Consider the fiberwise polynomial maps of degree at most $K$
and with zero constant term
 from ${\varepsilon}^N = M^n
\times \BR^N$ into $f^* TP^{n+k} \oplus \nu_M$, i.e.\ the bundle
$J^K_0({\varepsilon}^N, f^* TP^{n+k} \oplus \nu_M) \to M^n$ with fibre
$J^K_0(\BR^N, \BR^{N + k})$ (= polynomial maps from $(\BR^N, 0)$
to $(\BR^{N + k}, 0)$ of degree at most $K$).
Again in each fibre we can take only those polynomial maps which
have at the origin a singularity type belonging to~$\tau$, i.e.\
the subspace~$V^N_\tau$.
The bundle $J^K_0({\varepsilon}^N, f^* TP^{n+k} \oplus \nu_M) \to M^n$ can be
induced from the bundle $J^K_0({\varepsilon}^N, \gamma^{SO}_{N + k}) \to
B\,SO_{N + k}$.
\begin{defi}
Denote by $\Kaz_\tau(N)$ the subspace in $J^K_0({\varepsilon}^N,
\gamma^{SO}_{N + k})$ such that the intersection of $\Kaz_\tau(N)$
with any fibre $J^K_0(\BR^N, \BR^{N + k})$ is~$V^N_\tau$.
\end{defi}
The inclusion $SO_N \subset SO_{N + 1}$ induces a map $\Kaz_\tau(N)
\to \Kaz_\tau(N + 1)$.
Finally put $\Kaz_\tau = \lim\limits_{N \to \infty} \Kaz_\tau(N)$ and
call this space {\it {Kazarian space for $\tau$-maps}}.
\begin{defi}
\label{nu}
The lift of the $n$-universal virtual normal bundle $\nu(n)$ (see Definition~\ref{n-normal}) to the space
$\Kaz_\tau$ will be denoted by $\nu = \nu^k_\tau$ and called the {\it {universal
virtual normal bundle of  $\tau$-maps.}}
\end{defi} 

\begin{rema}\label{felbont}
Note that there is a partition of $\Kaz_\tau$ according to the
singularity classes in $\tau$, i.e. $\Kaz_\tau = \bigcup\limits_{[\eta]
\in \tau} B\eta$ and the analogue of  Lemma~\ref{Kaz}  holds, see the Appendix.
\end{rema}

\section{Definition of the Kazarian space by gluing}

For this second definition we will need the notion of ``global
normal form'' of a singularity type $[\eta]$ (see \cite{Sz8},
\cite{Sz9}, \cite{R--Sz}), hence we recall it now.
Let $\eta: (\BR^c, 0) \to (\BR^{c + k}, 0)$ be the root of the
singularity type~$[\eta]$.
(Recall that this means that $c$ is the minimal possible
dimension for the representatives of~$[\eta]$.)
Let $G^O_\eta$ denote the maximal compact
subgroup of the
automorphism group of $\eta$
(in the sense of J\"anich, see \cite{J}, \cite{W}, \cite{R--Sz}).
Let $\lambda_\eta$ and $\wt \lambda_\eta$ be the
representations of $G^O_\eta$ in the source and the target of
$\eta$, respectively
(i.e.\ for $\forall g \in G^O_\eta$, $\wt\lambda_\eta(g)
\circ \eta \circ \lambda^{-1}_\eta(g) = \eta$. 
Since we will consider always oriented source and target manifolds we will need only the
``orientation preserving'' subgroup of $G^O_\eta$, i.e. the group $G_\eta$
defined as follows:
$G_\eta = \{g \in G_\eta^O | \text {det} \lambda (g) \cdot \text{det}
 \wt \lambda (g) > 0\}$).
\begin{defi}\label{glob}
  The {\it {``global normal form''}} of $\eta$ consists of
  \begin{enumerate}[\indent 1)]
  \item
    the vector bundle $\xi_\eta = EG_\eta
    \underset{\lambda_\eta}{\times}\BR^c
    \overset{\BR^c}{\longrightarrow} BG_\eta$,
  \item
    the vector bundle $\wt\xi_\eta = EG_\eta
    \underset{\wt\lambda_\eta}{\times} \BR^{c + k}
    \overset{\BR^{c + k}}{\longrightarrow} B G_\eta$,
  \item
    the (non-linear) fiberwise map $\Phi_\eta : \xi_\eta
    \longrightarrow \wt\xi_\eta$ having the form
    
    $(\id_{EG_\eta} \times~\eta)\big/G_\eta$, where $\id_{EG_\eta}$ denotes
    the identity map of the contractible space $EG_\eta.$
  \end{enumerate}
\end{defi}
\noindent
{\bf Explanation:}
The map $\id_{EG_\eta} \times \eta$ is a $G_\eta$-equivariant
map from $EG_\eta \times \BR^c$ into $EG_\eta \times \BR^{c + k}$.
The map $\Phi_\eta$ is obtained  quotiening out by the action
of~$G_\eta$.
\begin{prop}[\cite{Sz7}, \cite{R--Sz}]
\label{pro:1a}
Let $M^n$ and $P^{n+k}$ be oriented manifolds of dimensions $n$ and $n+k$ respectively.
Suppose that $f: M^n \to P^{n + k}$ is a proper $\tau$-map,
and $f\big|_{\eta(f)}$ is an embedding into $P^{n+k}$.
Let us denote by $T$ and $\wt T$ the tubular
neighbourhoods of $\eta(f)$ in $M^n$ and $f(\eta(f))$ in $P^{n+k}$,
respectively.
Then there is a commutative diagram
\begin{equation}
\label{snormal}
\CD
T @> f\big|_T >> \wt T \\
@VVV @VVV \\
\xi_{\eta} @>\Phi_\eta >> \wt\xi_{\eta}
\endCD
\end{equation}
where the vertical maps are vector bundle morphisms, which are
isomorphisms on the fibres.\hfill $\square$
\end{prop}
\begin{rema}
This proposition shows that the universal map $\Phi_\eta$ describes the
$\tau$-map $f$ in
the tubular neighbourhood of the whole stratum $\eta(f)$, that
is why it was called a ``global normal form''.

\end{rema}

\begin{rema}
The relative version of the above Proposition also holds, with the same proof
as in \cite{R--Sz}. This means the following: Suppose that the vertical arrows of 
\eqref{snormal} are
given on $T|_A$ and $\wt T|_A$ for some closed subset $A \subset \eta(f)$,
so that the diagram is commutative, where the maps are defined, 
and the pair 
$(\eta(f), A)$ is a relative CW complex.
Then these given vertical arrows 
can be extended to have the commutative diagram \eqref{snormal}.
\end{rema}

Now we give another description of the Kazarian space by a glueing
procedure. In order to avoid the confusion of the two definitions we denote
this space by $gl\Kaz_\tau.$ 
{\footnote {The space $gl \Kaz_\tau$ will not be homeomorphic to the space
$\Kaz_\tau$, but it will be homotopically equivalent. Moreover
$gl \Kaz_\tau$ is a deformation retract in $\Kaz_\tau.$}}

Let $[\eta]$ be a top singularity type in $\tau$ and put $\tau ' = \tau
\setminus \{[\eta]\}.$ Suppose that the space $gl\Kaz_{\tau'}$ has been constructed
and it
has been shown that $gl \Kaz_{\tau'}$ is homotopically equivalent to the
space $\Kaz_\tau'.$ Let us denote by $v':gl\Kaz_{\tau'} \to \Kaz_{\tau'}$
the constructed homotopy equivalence and by $u'$ we denote its homotopy
 inverse. 
(The starting step of the induction will be given at the end of the proof.)
Let $BG_{\eta,fin}$ be a finite dimensional manifold approximating the
classifying space $BG_\eta.$ (For example on the Stiefel manifold $V_c(R^N)$
of $c$-tuples of vectors in $R^N$ there is a free $G_\eta$-action, and this
manifold is 
$N-c-1$-connected, where $c$ is the codimension of the singularity $\eta.$
So we can put $BG_{\eta,fin} = V_c(\BR^N)/\lambda_\eta(G_\eta).)$

Let us denote the restrictions of the bundles  $\xi_\eta$ and $\wt\xi_\eta$
to $BG_{\eta, fin}$by  $\xi_{\eta,fin}$ and $\wt \xi_{\eta,fin}$ respectively, and their disc
bundles by $D\xi_{\eta,fin}$ and $D\wt \xi_{\eta,fin}.$
Finally let $\Phi = \Phi_{\eta,fin}$ denote the restriction of $\Phi_\eta$ to 
$D\xi_{\eta,fin}$, see Diagram \eqref{snormal} . (The disc bundles are chosen so that
   $D\xi_{\eta,fin} = \Phi_\eta^{-1}(\wt D\xi_{\eta,fin}).$)
Moreover this map $\kappa_f'$ has the following property.
By Remark~\ref{felbont}  
 the map $\Phi: D\xi_{\eta,fin} \to \wt D\xi_{\eta,fin} $ defines a map $\kappa_\Phi: D\xi_{\eta,fin}
\to \Kaz_\tau.$ This map carries a point $x \in D\xi_{\eta,fin}$
into such a point of $\Kaz_\tau$ that corresponds to a singularity of the same
type as the germ of $\Phi$ at $x.$ Hence the preimage of the $\eta$-stratum
$B \eta \subset \Kaz_\tau$ is the zero section of the bundle 
$D\xi_{\eta,fin}.$ 
Finally let us denote by $\kappa_\Phi^S$ the restriction of the map
$\kappa_\Phi$ to the sphere bundle $S(\xi_{\eta,fin}) = \partial D(\xi_{\eta,fin}).$

\begin{rema}\label{puback}

It follows from 
Remark\ref{felbont}
that $\kappa_\Phi$ pulls back from the
universal virtual normal bundle $\nu_\tau$ the virtual bundle
$\pi_{D(\xi_\eta)}^*(\wt \xi_\eta - \xi_\eta)$ where 
$\pi_{D(\xi_\eta)}$ is the projection $D(\xi_\eta) \to BG_\eta.$ 

\end{rema}

\begin{defi}
We define a finite dimensional approximation of the  space $gl\Kaz _\tau$ as follows:
$$gl \Kaz _{\tau,fin} \overset{def} = gl\Kaz _{\tau',fin}
 \underset{u' \circ \kappa^S_\Phi} \bigcup  D\xi_{\eta,fin}
= gl\Kaz _{\tau'}\sqcup D\xi_{\eta,fin}\Bigm/ x \sim u' \circ \kappa^S_\Phi
(x),$$
for\ all $\ x \in S(\xi_{\eta, fin}),$ i.e. in the disjoint union of 
$gl\Kaz_\tau'$ and $D\xi_{\eta,fin}$ we identify each $x \in 
S(\xi_{\eta, fin})$ with its image $u' \circ \kappa^S_\Phi(x) \in
gl\Kaz_{\tau'}.$ Recall that $u': \Kaz_{\tau'} \to gl\Kaz_{\tau'}$
is a homotopy equivalence, that exists by the assumption of the induction.
Its homotopy inverse is $v'.$

\end{defi}

\begin{defi}
Let  $v: gl\Kaz_\tau \to \Kaz_\tau$ be a continuous map with the following properties:
\begin{enumerate}[\indent a)]

\item
$v\big|_{gl\Kaz_{\tau'}} \cong v'$
\item
$v\big|{D(\xi_\eta}) = \kappa_\Phi$

Such a map exists, because on the sphere bundle $S(\xi_\eta)$ the two maps
on the right hand sides of 
(a) and (b) are homotopic. 
\end{enumerate}
\end{defi}
\begin{clai}\label{v}

The map $v$ is a homotopy equivalence.
\end{clai}
\begin{proof}
It is enough to show that the map
$$
\hat \kappa_\Phi: (D(\xi_\eta), S(\xi_\eta)) \to (\Kaz_\tau, \Kaz_{\tau'})
$$
defined by $\kappa_\Phi$
induces isomorphisms of the homology groups.

\begin{lem}
\label{transvslice}
$\kappa_\Phi$ maps each fibre of the bundle $D(\xi_\eta)$ into a slice
transverse to $B\eta$ in $\Kaz_\tau.$
\end{lem}
\begin{proof}
Postponing the proof of this Lemma we prove Claim~\ref{v}.
After an isotopy it can be supposed that $\kappa_\Phi$ maps each fibre
of $D(\xi_\eta)$ into a fibre of the normal bundle $\nu_\eta$ of the stratum
$B\eta$ in $\Kaz_\tau.$
Let $g: BG_\eta \to B\eta$ be the map defined by the map $\kappa_\Phi$
mapping the zero section of $\xi_\eta$ into that of $\nu_\eta.$
Hence $g^*(\nu_\eta) = \xi_\eta.$
On the other hand $\xi_\eta$ is the universal bundle with structure group
$G_\eta$, and $\nu_\eta$ admits the group $G_\eta$ as structure group. Hence
there is a map $h: B\eta \to BG_\eta$ such that $h^*(\xi_\eta) = \nu_\eta.$
We obtain that the composition map $h\circ g: BG_\eta \to BG_\eta$ is
homotopic to the identity map (because $(h\circ g)^*(\xi_\eta) = \xi_\eta$ and
the inducing map is homotopically unique).
Since the homology groups of $BG_\eta = B\eta$ are finitely generated Abelian
groups in each dimension we obtain that the homomorphisms induced by the maps 
$h$ and $g$ are isomorphisms.
Let us denote by $\hat g$ the map of the corresponding Thom spaces induced by
$g$, i.e. $\hat g : T\xi_\eta \to T\nu_\eta.$ By the Thom isomorphism $\hat g$
also induces isomorphism of the homology groups, since so does $g.$
But the homology groups of the Thom spaces $T\xi_\eta$ and $T\nu_\eta$ are
canonically isomorphic to those of the pairs 
$(D(\xi_\eta), S(\xi_\eta))$ and  $(\Kaz_\tau, \Kaz_{\tau'})$ respectively,
and identifying the corresponding homology groups using these canonical
isomorphisms the homomorphism $\hat g_*$ turns into $\hat \kappa_{\Phi*}$.
Now by the five lemma the map $v$ induces isomorphism of the homology groups.
\end{proof}
\end{proof}

\begin{proof} of Lemma~\ref{transvslice}.
Let us use  Remark~\ref{slice} with the substitutions:
\begin{enumerate}[\indent a)]

\item
$ V = J^K_0(\BR^N, \BR^{N + k}),$ i.e.  
$V$ is the space of polynomial maps  $\mathcal P =
 \mathcal P^K(\BR^N,\BR^{N+k})$  from $\BR^N$ into $\BR^{N+k}$
of degree (at most) $K$ having zero constant term.
\item
$\Sigma = [\eta] \cap \mathcal P.$
This means the following: Let $\eta: \BR^c \to \BR^{c+k}$
be a polynomial map of degree (at most) $K$ giving  the root of the stable
singularity class $[\eta].$ Then $\Sigma$ consists of all the polynomial maps from
$\mathcal P$
equivalent to $\eta$ in the sense of Definition~\ref{susp}.
\end{enumerate}

The map 
$\kappa_\Phi$ restricted to a fibre $\BR^c$ can be described as follows:
$a \in \BR^c$ is mapped into the polynomial map $f_a$ in $\mathcal P$
such that:
if $x \in \BR^c$, $y \in \BR^{N-c}$ and $\BR^N = \BR^c \times \BR^{N-c}$, then 
$f_a(x, y) = ((\eta (x-a) - \eta(a)), y).$ 
By Remark~\ref{slice} it is enough to show that the image of $\BR^c$ is
transverse to $[\eta]\cap \mathcal P$ in $\mathcal P.$
If $\mathcal U$ is a small enough neighbourhood of $f_0$ in $\mathcal P$
then for each $g \in \mathcal U$ a linear subspace $L_g$ can be defined in
$\BR^N$ such that $L_g$ contains $\text{ker}dg$ and $\text{dim} L_g =
\text{dim}df_0$. We denote by $\ell$ the dimensions of these spaces $L_g.$
Again if $\mathcal U$ is small enough then all the spaces $L_g$ can be 
projected isomorphically onto $L_{f_0}$, and so they can be identified with 
$\BR^{\ell}$ in a continuous way.

We define a map 
$$
\Pi: \mathcal U \to   
\mathcal P^K(\BR^{\ell},\BR^{N+k})$$
by taking the restriction of each $g \in \mathcal U$ to $L_g.$
Note that $\Pi(g)$ is either the genotype of $g$ or a trivial extension
(i.e. a multiple suspension) of it.
{\footnote{Recall that the genotype of the germ $g$ can be defined as its
  restriction to the kernel of its differential.}} 

 Further clearly $\Pi$ is a submersion onto
its image.
A theorem of Mather says that 
two stable germs $f$ and $g$ are $\mathcal A$ equivalent if and only if their
genotypes are contact equivalent. From this it follows that the preimage of
the  contact 
orbit of $\Pi(f_0)$ is the $\mathcal A$ orbit of $f_0.$

$\Pi(\kappa_\Phi(\BR^c))$ is the deformation of the genotype of $f_0$, which
is the same as the deformation of the genotype of $\eta$ given by $\eta.$ Again by a theorem of Mather
such a deformation is transverse to the contact orbit of the genotype if the 
germ itself (in this case $\eta$) is $\mathcal A$ stable.
Hence $\Pi(\kappa_\Phi(\BR^c))$ is a transverse slice to the contact orbit of 
the genotype of $f_0.$  But then the lift of this slice
$\kappa_\Phi(\BR^c))$
 is transverse to the preimage of the contact orbit,  which is precisely the
 $\mathcal A$ orbit of $f.$

The start of the induction:
Recall that  the simplest ``singularity type'' $[\eta_0]$ is always
the class of maps of maximal rank $\Sigma^0$.

We show that for $\tau = \{[\eta_0]\}$ both spaces $\Kaz_\tau$ and  $gl\Kaz_\tau$
are homotopy equivalent to the Grasmann manifold  $B\,SO(k).$ We denote these
spaces here by $\Kaz_0$ and $gl\Kaz_0$ respectively.

  \begin{enumerate}[\bf 1)]
  \item
    By the gluing procedure we have to consider the root of
    codimension $k$ maps  of maximal rank, this is the germ
    $(0, 0) \hookrightarrow
    (\BR^k, 0)$.
    The automorphism group of this germ is the group of germs of all
    diffeomorphisms $\text{\rm Diff}(\BR^k, 0)$.
    Its maximal compact subgroup is $G^O_{\eta_0} = O(k)$.
    The oriented version of this group $G_{\eta_0}$ is $SO(k)$.
    The representations $\lambda$ and $\wt\lambda$ on the source
    $(0)$ and on the target $(\BR^k)$ respectively are the trivial map $\lambda : SO(k)
    \to \{1\}$ and the standard inclusion $\wt\lambda : SO(k) \to O(k).$
    Now the space $gl \Kaz_0(k)$ is $E\, SO(k)
    \underset{SO(k)}{\times} \{0\} = B\, SO(k)$.
  \item
    By the Borel construction we obtain the following:
    Let us denote by $\mathcal P^K_{\eta_0}$ the space of polynomial maps
    $\BR^N \rightarrowtail\BR^{N + k}$
    of
    degree (at most) $K$ having differentials of maximal rank at
    the origin. The group $K$-jets of germs of diffeomorphisms of
    $\BR^{N+k}$ acts on this space. Associating to each polynomial map and to
    each germ of diffeomorphisms their
    linear parts at the origin we obtain an equivariant homotopy equivalence
from $\mathcal P^K_{\eta_0}$
    to the Stiefel manifold
    $V_N(\BR^{N  + k})$ provided with the action of the linear group
    $GL^+(\BR^{N+k})$. Replacing the linear group by the homotopically
    equivalent orthogonal group we obtain
    that 
    $$
    \Kaz_{0} =\underset{N \to \infty}{\text {lim}} V_N(\BR^{N + k}) 
    \underset{SO(N + k)}{\times} E\,SO(N + k).$$
    Since $V_N(\BR^{N + k}) = SO(N + k) / SO(k)$ we conclude
    $$\Kaz_{0}(k) =\underset{N\to \infty}{\text{lim}} E\, SO(N + k) \bigm/
    SO(k) = B\, SO(k).$$
  \end{enumerate}
\end{proof}
From now on we will use the notation $\Kaz_\tau$ for the space $gl\Kaz_\tau.$
In particular we will use the decomposition 
$\Kaz_\tau = \underset{i} \bigcup D\xi_{\eta_i}.$

\section{The Kazarian spectral sequence}

Let $\eta_0 < \eta_1 < \dots$ be a linear ordering of the
singularities in $\tau$ compatible with the hierarchy of singularities.
Let $\tau_i$ be the set $\bigl\{\eta_j \mid j \leq i \bigr\}$
and denote $\Kaz_{\tau_i}$ by $\Kaz_i$, $G_{\eta_i}$ by $G_i$,
$\xi_{\eta_i}$ by $\xi_i$ etc.
Let us consider the homological spectral sequence arising from
the filtration $\Kaz_0 \subset \Kaz_1 \subset \dots \subset \Kaz_\tau$.
Then $E^1_{p,q} = H_{p + q} (T\xi_q)$, where $\xi_q =
\xi_{\eta_q}$ and $T\xi_q$ is the Thom space. By the Thom isomorphism
\[
E^1_{p,q} \approx H_{p + q - c_q} \bigl(BG_q; \wt \Z(\xi_q)\bigr),
\]
where $c_q$ is the dimension of the bundle $\xi_q$, the
coefficients $\wt \Z(\xi_q)$ are the integers twisted by the
orientation of~$\xi_q$.
This spectral sequence converges to the group $H_*(\Kaz_\tau)$,
i.e.\ $\bigoplus\limits_{p + q = n} E^\infty_{p,q}$ is
associated to $H_n(\Kaz_\tau).
\footnote{Kazarian's ordering is different from ours.}$

\part{Pontrjagin - Thom construction for $\tau$-maps.}

In our proof of Pontrjagin - Thom construction for $\tau$-maps
we follow the method of Wells \cite{We} reducing the  cobordisms of immersions to
those of embeddings with many normal vector fields by lifting  immersions
into a high dimensional space.
 (This method  is completely different from 
that of \cite{R--Sz}, not to mention \cite{A} or \cite{Sa}.)

\section{Definition of $\tau$-embedding}

\begin{defi}
  Let $\tau$ be a set of stable singularity classes of codimension $k>0$.
  A quintuple $(M^n, Q^q, e, V, \mathcal F)$ is
  called a{\it { $\tau$-embedding}} if
  \begin{enumerate}[\indent 1)]
  \item $M^n$ and $Q^q$ are manifolds (possibly with non-empty
    boundaries) of dimensions $n$ and $q$ respectively.
  \item $e : M^n \hookrightarrow Q^q$ is an embedding.
    (If $M^n$ has boundary, then $Q^q$ has it too,
    $e^{-1}(\partial Q^q) = \partial M^n$,
    and $e$ is transverse to the boundary of~$Q^q$.)
  \item $V = \{v_1, \dots, v_N\}$ is a finite set of linearly
    independent vector fields along $e(M^n)$ (i.e.\ sections of the
    bundle $TQ^q\big|_{e(M^n)}$). Here $N = q - n - k.$ We shall identify $V$ with the $N$-dimensional
    trivialised sub-bundle of $TQ^q|_{e(M^n)}$ generated by the vector fields. 
    We require that the
    vector fields
    on the boundary $e(\partial M^n)$ be
    tangent to~$\partial Q^q.$
    \item $\mathcal F$ is a foliation of dimension $N$ on a neighbourhood of $e(M^n)$, and it is
      tangent to $V$ along $e(M^n).$
    \item In a neighbourhood of any point $p \in M^n$ the
      composition of the embedding $e$ with the local projection along
      the leaves of the foliation $\mathcal F$ gives a local map (onto a small
      $n+k$-dimensional slice transverse to the leaves), having at
      $p$ a
      singularity belonging to~$\tau$.
    \end{enumerate}

\end{defi}

\begin{rema}\
  \begin{enumerate}[1)]
  \item
   Note that $N = \text{dim}\ Q^q - \text{dim} M^n - k.$ When the codimension 
   of
   the $\tau$-embedding is not clear from the context, then we have to include $N$ in the
   notation. In this case we shall write $(\tau,N)$-embedding for the above
   defined $\tau$-embedding.
  \item
    The quintuple $(M^n, Q^q, e, V, \mathcal F)$ often will be denoted simply by
    $e$, and we shall say that $e$ is a $\tau$-embedding,
  \end{enumerate}
\end{rema}

\begin{exam}
Let $f: M^n \to P^{n + k}$ be a $\tau$-map, let $g : M^n
\hookrightarrow \BR^N$ be an embedding and $Q^q = P^{n+k} \times \BR^N, \ q = n+k+N$. Put
$e = f \times g$ and $V = \{v_1, \dots, v_N\}$, where $v_1,
\dots, v_N$ are the vector fields arising from a basis in~$\BR^N.$
Let $\mathcal F$ be the vertical foliation with the leaves $\{x\}\times \BR^N.$
Then the quintuple $(M^n, Q^q, e, V, \mathcal F)$ is a $\tau$-embedding.
\end{exam}
\begin{defi}
\label{vert}
In $P^{n+k} \times \BR^N$ we shall call the above defined vector fields
$v_1, \dots, v_N$ \it {vertical}.
\end{defi}
\begin{rema}
If $(M^n, Q, e, V, \mathcal F)$ is a $\tau$-embedding, where $Q = P^{n + k}
\times \BR^N$ and $V = \{v_1, \dots, v_N\}$ are the vertical basis
vector fields, $\mathcal F$ is formed by the leaves $\{p\}\times R^N$, for
$p\in P^{n+k}$, then the composition $f = \pi \circ e$ is a 
$\tau$-map $f :
M^n \to P^{n + k}$
, where 
$\pi: P^{n+k} \times \BR^N \to P^{n+k}$ is the projection. 

\end{rema}
\begin{rema}
Given a $\tau$-embedding $(M^n, Q, e, V, \mathcal F)$ a
stratification arises on the manifold $M^n$;  $M^n = \cup \{\eta(\pi \circ e) \mid [\eta] \in \tau \}$
according to the singularity types of the composition of $e$ with
the local projection $\pi$ along the leaves. Sometimes we shall write $\eta(e)$
instead of $\eta(\pi \circ e).$
\end{rema}

\section{The Compression Theorem for stratified manifolds.}

\begin{defi}(Gromov, \cite{Gr})
Let $\mathcal S$ be a finite stratification on a compact
manifold $M^n$:  $M^n = S_1 \cup \dots \cup S_i$, \ $S_j$ are the
strata, $S_j \subset \overline S_{j - 1}$.
Let $e$ be an embedding of $M^n$ in a manifold $Q^q, \ q = n+k+N$ (with the usual
assumption on the boundary) and let $v_1, \dots, v_N$ be
linearly independent vector fields along $e(M^n)$ (at $e(\partial M^n)$ they are tangent
to $\partial Q^q$).
For $x \in S_j$ denote by $V_x$ the vector space generated by
the vectors $v_1, \dots, v_N$ at $e(x)$.
So $V_x \subset T_{e(x)} Q^q$. Let $de$ be the differential of the embedding $e.$
We say that the vector fields $v_1, \dots, v_N$ are
{\it {transverse to the stratification}} $\mathcal S$ if for any $x \in
M^n$ the space $V_x$ and the $de$-image $de(T_xS_j)$ of the tangent space
 of the stratum $S_j$ at $x$ intersect only in the
origin of the space $T_{e(x)}Q^q$.
\end{defi}
\begin{thm}[Stratified Compression Theorem, SCT]
  \label{SCT}
  \
  \begin{enumerate}[a)]
  \item
    Let $M^n$ be a compact smooth manifold, 
    $P^{n + k}$ be an arbitrary smooth
    $(n + k)$-manifold, 
    $e : M^n \hookrightarrow P^{n+k} \times \BR^N$ an embedding.
    Let $\mathcal S$ be a finite stratification of $M^n : M^n = S_1 \cup
    \dots \cup S_i$.
    Further let $v_1, \dots, v_\ell$, $\ell \leq N,$ be vector fields
    along $e$ and transverse to the stratification~$\mathcal S$.
    Then there exists an isotopy $\varphi_t : P^{n+k} \times \BR^N \to
    P^{n+k} \times \BR^N$, $t \in [0,1]$, such that $\varphi_0$ 
    is the  identity map, and the vector fields $d\varphi_1(v_1),\dots,
    d\varphi_1(v_\ell)$ are vertical vector fields along the
    embedding~$\varphi_1 \circ e$, (see Definition~\ref{vert}).
  \item
    (Relative version)
    Suppose we are in the previous conditions and
    $K \subset M^n$ is a compact in $M^n$ such that on a
    neighbourhood of $K$ the vector fields $v_1, \dots, v_\ell$ are vertical.
    Then the isotopy $\varphi_t$ can be chosen so that it keeps fixed
    the points of a neighbourhood of~$K$.
  \end{enumerate}
\end{thm}
\begin{proof} of $b).$
Induction on $i$ ($=$ the number of strata).
For $i = 1$ this is the usual Compression Theorem for
immersions, see \cite{Ro--Sa}.
It can be supposed that the vertical vector fields are given not only along
the submanifold $e(M^n)$ but also along a tubular neighbourhood of $e(M^n)$  in
$P^{n+k}\times R^N$ and are tangent to the leaves of the foliation.
(Such an extension always exists and it is homotopically unique.)
Now suppose that the vector fields are vertical already on a
neighbourhood $U$ (in $P^{n+k} \times \BR^N$) of the union of the $e$-images of the
strata: $S_i \cup S_{i - 1} \cup \dots \cup S_{i - j}$.
(Recall that $\dim S_i < \dim S_{i - 1} < \dots < \dim S_{i - j}
< \dots < \dim S_1 = n$.)
Let us choose neighbourhoods $U'$, $U''$ of the union $S_i \cup
\dots \cup S_{i - j}\cup K$  such that $U'' \subset  \overline U''
\subset U' \subset \overline U' \subset U$.
Put $M^* = M^n\setminus U''$ and $K^* = (M^n \setminus U'') \cap \overline U'$.
Then $M^*$ is stratified by the strata $S_1 \cup \dots \cup S_{i - j
- 1}$ intersected with~$M^*$.
The vector fields are made vertical on $U$, in particular on~$K^*$.
By the relative version of the inductional assumption we can
make the vector fields vertical on $M^*$ (or even on a
neighbourhood of it in $P^{n+k} \times \BR^N$) by an isotopy fixed on
$K$.
So we get vector fields vertical on the whole~$M^n$.
\end{proof}

Now we show that starting with a $(\tau, N)$-embedding in $P^{n+k}\times \BR^N$
obtained as a result of the isotopy described in the previous Lemma
(i.e the vector fields are vertical along $e(M^n)$) a further isotopy
(fixed on $e(M^n)$) can be
applied that will turn the foliation vertical in a
neighbourhood of $e(M^n).$ (Recall that the foliation is called vertical if its
leaves coincide with the sets $\{p\}\times\BR^N$ in a neighbourhood of $e(M^n).)$

\begin{lem}\label{vertfolia}
Let $(M^n,P^{n+k}\times \BR^N,e, V,\mathcal F)$ be a $(\tau,N)$-embedding
such that the vector fields $v_1, \dots, v_N$ along $e(M^n)$ are vertical, i.e.
coincide with the standard vector fields.
Then there is an isotopy $\psi_t$ of $P^{n+k}\times \BR^N$ such that
$\psi_0$ is the identity map, $\psi_1(\mathcal
F) = \mathcal F'$, where $\mathcal F'$ denotes the vertical foliation
formed by the leaves $\{x\}\times \BR^N, \ x \in P^{n+k},$ 
in a neighbourhood of $e(M^n).$
Moreover, $d\psi_1(v_j) = v'_j$ for $j = 1,\dots, N$, 
where $\{v_1', \dots v_N'\}$ 
denote the standard vertical vector fields arising from the standard
basis of $\BR^N$  .
\end{lem}

\begin{defi}
Let us call the subsets $P^{n+k}\times \{y\}$ in $P^{n+k} \times \BR^N$,
 where $y
\in \BR^N,$ {\it { horizontal sections}}.
\end{defi}

\begin{proof}
We start with two general remarks:
\begin{enumerate}[\bf 1)]
\item
  If $L$ and $L'$ are two subsets of $P^{n+k} \times \BR^N$ such that
  they intersect each horizontal section $P^{n+k} \times \{y\}$ for $y \in U
  \subset \BR^N$ in exactly one point 
  , then a one-to-one map $L \to
  L'$ arises.
\item
  If $A$ is a subset in $P ^{n+k}\times \BR^N$,
  and for each $a = (p, y)$, $a \in A,$ subsets $L_a$ and $L'_a$
  are given such that $a \in L_a$, $a \in L'_a$ and $L_a$, $L'_a$
  intersect each horizontal section $P^{n+k} \times \{y'\}$ for $y'$ close
  to $y$ in a unique point, then a family of maps $L_a \to L'_a$ arises. Suppose
  that $L_{a_1} \cap L_{a_2} = \emptyset$ and $L'_{a_1} \cap L'_{a_2} =
  \emptyset$ for $a_1 \neq a_2$.
  Then a one-to-one continuous map $\alpha: \bigcup\limits_{a
    \in A} L_a \to \bigcup\limits_{a \in A} L'_a$ arises by associating to each
  other the points on the same horizontal section.
  If $A$, $L_a$, $L'_a$ are smooth submanifolds so that $\{ \cup
  L_a \mid a \in A\}$ and $\{\cup L'_a \mid a \in A\}$ are also smooth
  submanifolds, then the map $\alpha$ is a smooth map.

\medskip

  Now let $V^\perp$ be the sub-bundle in $T(P^{n+k} \times \BR^N)\big|_{e(S_i)}$ formed by the
  orthogonal complement of $V\big|_{e(S_i)} \oplus Te(S_i).$ (Recall that
  $S_i$ is the smallest dimensional stratum.) 
  Let
  $\exp(V^\perp)$ be the image of $V^\perp$ at the exponential map.
  Choose for $A$ a small neighbourhood of $e(S_i)$ in $\exp(V^\perp).$
  If $A$ is chosen sufficiently small, then the local leaves $L_a$
  and $L'_a$ of the foliations $\mathcal F$ and
  $\mathcal F'$ around $a \in A$
  will satisfy the conditions that for $a_1 \ne a_2$ 
  the intersections
  $ L_{a_1} \cap
  L_{a_2}$ and $L'_{a_1} \cap L'_{a_2}$ are empty, and so a 
  map $\alpha : \{\cup L_a \mid a \in A\} \to \{ \cup L'_a \mid a
  \in A\}$ arises.
  Note that $U = \{ \cup L_a \mid a \in A\}$ and $U'
  = \{ \cup L'_a \mid a \in A\}$ are neighbourhoods of the stratum
  $e(S_i)$ in $P^{n+k} \times \BR^N$.
  The map $\alpha$ will
  have differential equal to the identity at the points $a \in
  e(S_i)$. Hence $\alpha$ is a diffeomorphism in a neighbourhood of
  $e(S_i)$ and $\alpha(\mathcal F) = \mathcal F'$ around $e(S_i)$.
  Joining the points $x$ and $\alpha(x)$ by a minimal
  geodesic curve in the horizontal section they both belong to, we see
  that $\alpha$ can be embedded into an
  isotopy $\alpha_t, (\ 0 \le t \le 1)$ such that $\alpha_0 =$ identity,
  $\alpha_1 = \alpha$.
  Next we repeat the same procedure around $S_{i - 1}$.
  Note that where the foliations $\mathcal F$ and $\mathcal F'$
  already coincide the above procedure gives  automatically the identity
  diffeomorphism. And so on.
  In $i$ steps we obtain a diffeomorphism (moreover an isotopy)
  mapping $\mathcal F$ into $\mathcal F'$ around $e(M^n)$. Q.E.D.
\end{enumerate}
\end{proof}

\section{Cobordism of $\tau$-embeddings}

\begin{defi}
The {\it {cobordism of $\tau$-embeddings}} $(M^n, Q^q, e, V,\mathcal F)$ of
$n$-manifolds in the manifold $Q^q$ can be defined in an obvious way.
The cobordism is a $\tau$-embedding in $Q^q \times[0,1].$
Let us denote the set of equivalence classes by $\text{\rm
  Emb}_{(\tau,N)}(Q^q)$, where $N = q - n - k$.
\end{defi}

The next theorem shows that the computation of the cobordism group of
$\tau$-maps can be reduced to that of $\tau$-embeddings.

\begin{thm}
\label{tau-embeddings}
Let $P^{n+k}$ be a fixed oriented manifold of dimension $n+k$, where $k>0$ and consider the
group $\Cob _\tau(P^{n+k})$ of $\tau$-maps of oriented $n$-manifolds  in
$P^{n+k}$, where 
$\tau$ is a family of stable singularities, such that for each natural
number $c$
there are finitely many elements in $\tau$ of codimension $c.$  
If $N$ is sufficiently big, then
\[
\aligned
\Cob_\tau(P^{n+k}) &\approx \Emb_{(\tau,N)}(P^{n+k} \times \BR^N).
\endaligned
\]
\end{thm}
\begin{proof}\
  \begin{enumerate}[1)]
  \item
    We define a map $\alpha : \Cob_\tau(P^{n+k}) \to \Emb_{(\tau,N)}(P^{n+k} \times
    \BR^N)$ as follows:
    Let $f : M^n \to P^{n + k}$ be a $\tau$-map, and let $g : M^n
    \hookrightarrow \BR^N$ be an arbitrary embedding.  Put $e = f
    \times g$.
    Then $e : M^n \hookrightarrow P^{n+k} \times \BR^N$ is an embedding and
    choosing for $\{v_1, \dots, v_N\}$ the standard vertical
    vector fields (defined by the standard basis in $\BR^N$)
and for $\mathcal F$ the vertical foliations with leaves $\{p\}\times R^N,
p\in P^{n+k}$,
 we obtain a 
$\tau$-embedding $(M^n, P^{n+k} \times \BR^N,
    e, V, \mathcal F)$, where $V = \{v_1, \dots, v_N\}.$
    The cobordism class of this $\tau$-embedding does not depend on the
    choice of $g.$ Indeed, if $g'$ is another embedding $g' : M^n
    \hookrightarrow \BR^N$, then $g$ and $g'$ can be joined by
    an isotopy $g_t$, $t \in [0,1]$, $g_0 = g$, $g_1 = g'$ and
    putting $G(x,t) = (g_t(x), t) \in \BR^N \times [0,1]$
    for $x \in M^n$, $t \in [0,1]$ the $\tau$-embedding $E = f\times G : M^n \times
    [0,1] \to P^{n+k} \times \BR^N \times [0,1]$, (where $V$ is chosen to be the
    vertical foliation
    for each $P^{n+k} \times \BR^N \times \{t\}$),
    gives a cobordism joining $f \times g$ with $f\times g'$.

    Now let us choose any other representative $f'$ of the same
    $\tau$-cobordism class $[f] \in \Cob_\tau(P^{n+k})$.
    Then there exist an $(n + 1)$-manifold
    $W^{n + 1}$ such that $\partial W = -M^n \cup M'$, and a $\tau$-map
    \[
    F:(W; M^n; M') \to (P^{n+k}\times [0,1]; P^{n+k} \times \{0\}; P^{n+k}\times \{1\}).
    \]
    such that $F\big|_M^n$ and $F\big|_{M'}$ coincide with $f$ and $f'$ respectively.
    Further there is an embedding $G^W : W^{n + 1} \hookrightarrow
    \BR^N$ such that $G^W\big|_{M^n} = g$ and $G^W\big|_{M'} = g'$.
    Then $E^W = F\times G^W : W^{n + 1} \to P^{n+k} \times [0,1] \times \BR^N$
    is a $\tau$-embedding giving a cobordism between $e = f \times
    g$ and $e' = f' \times g'$.
    Hence the map $\alpha : \Cob_\tau(P^{n+k}) \to \Emb_{(\tau,N)}(P^{n+k}\times \BR^N)$
    is well defined.
  \item
    Next we define a map $\beta : \Emb_{(\tau,N)}(P^{n+k}\times \BR^N) \to
    \Cob_\tau(P^{n+k})$.
    Let $(M^n, P^{n+k} \times \BR^N, e, V,\mathcal F)$ be a $\tau$-embedding.
    By Theorem~\ref{SCT} and Lemma~\ref{vertfolia} the foliation $\mathcal F$
    can be made vertical by an isotopy $\varphi_t, \ \varphi_0 =
    \text{identity},$ and $\varphi_1 \circ e$ is a $(\tau,N)$-embedding with
    vertical leaves. 
    Projecting $\varphi_1 \circ e : M^n \hookrightarrow P^{n+k} \times \BR^N$
    to $P^{n+k}$ we obtain a $\tau$-map $f : M^n \to P^{n+k}$.
    We have to show that
    \begin{enumerate}[\indent a)]
    \item the cobordism class of $[f]$ is independent of the
      choice of the isotopy $\varphi_t$, and
    \item if $e'$ is another $\tau$-embedding cobordant to $e$,
      then applying the above procedure to $e'$ we obtain a
      $\tau$-map~$f'$ \ $\tau$-cobordant to $f$.
    \end{enumerate}
    Both of these statements follow from the relative version of the
    Stratified Compression Theorem.
    First we prove b). Suppose that $e$ and $e'$ are cobordant
    $\tau$-embeddings joined by a $\tau$-embedding~$E$.
    Apply isotopies $\varphi_t$ and $\varphi'_t$ turning vertical the
    leaves of the foliations of the $\tau$-embeddings $e$ and $e'$ 
    respectively.
    $\varphi_t$ and $\varphi'_t$ are defined on $P^{n+k} \times \BR^N \times\{0\}$
    and $P^{n+k} \times \BR^N \times \{1\}$, respectively.
    They can be extended to an isotopy $\Phi_t$ of the whole
    cylinder $P^{n+k} \times \BR^N \times [0,1]$.
    As a result we obtain a $\tau$-embedding $\Phi_1 \circ E$
    with boundary such that its foliation is vertical on
    the boundary.
    Now applying the relative version of Theorem~\ref{SCT}
    and Lemma~\ref{vertfolia} we obtain a
    $\tau$-embedding cobordism in $P^{n+k} \times \BR^N \times [0,1]$ with
    vertical leaves everywhere.
    Projecting it to $P^{n+k}\times [0,1]$ we obtain a $\tau$-cobordism
    between the $\tau$-maps
    $\pi \circ \varphi_1' \circ e'$ and $\pi \circ \varphi_1 \circ e$
    , where $\pi: P^{n+k} \times \BR^N \to P^{n+k}$ is the projection.
    That proves b).
    In order to prove a) it is enough to put $e' = e$ in the previous
    argument.
  \end{enumerate}
\end{proof}

\begin{rema}\label{ell-framed}
The same proof shows the following:
\begin{enumerate}
\item
If $N$ is sufficiently big (compared with $n,k$ and $\ell$) , then
\[
\aligned
\Cob_{\tau\oplus \ell}(P^{n+k}) &\approx \Emb_{(\tau,N+\ell)}(P^{n+k} \times
\BR^N), \ \  \text{and}
\endaligned
\]

\item
$\Emb_{(\tau,N)}(P^{n+k} \times \BR^N)$ does not
depend on $N$, moreover the  embedding $\BR^N \subset \BR^{N + 1}$
induces an isomorphism $$\Emb_{(\tau,N)}(P^{n+k} \times \BR^N) \approx
\Emb_{(\tau, N+1)}(P^{n+k} \times \BR^{N + 1}).$$
\end{enumerate}
\end{rema}

\begin{defi}

For any natural numbers $n, N, k$, any manifold $Q^{n+k+N}$ and
set of stable singularities $\tau$ of codimension $k$ maps
put
$$\Emb_{(\tau,N)}^{\text{Stable}} (Q^{n+k+N})\overset{\text {def}} = \underset{L\to \infty}{\text {lim}}\ \Emb_{(\tau,N+L)} (Q^{n+k+N} \times \BR^{L}).$$

\end{defi}

Our conclusion from this section is that in order to prove
the Pontrjagin--Thom construction  
for $\tau$-maps it is enough to do that for $(\tau, N)$-embeddings with big $N$. 
The key tool for the latter will be the normal form of a $\tau$-embeddings.

\section{Normal forms of $\tau$-embeddings of bounded source dimensions.}

\subsection{Normal form around a stratum.}

Fix a natural number $n.$ In this subsection
we shall consider $(\tau, N)$-embeddings with source
manifolds of dimensions less or equal to $n.$
Let $\eta$ be a top singularity stratum in $\tau.$
Choose a finite dimensional approximation $BG_{\eta,fin}$ of $BG_\eta$ such
that the pair $(BG_\eta, BG_{\eta, fin})$ is $(n+1)$-connected.
(Hence for the purpose of classification of $G_\eta$-bundles with base space
 of dimensions less
or equal to $n$ the approximation $BG_{\eta,fin}$ is fine enough.)
It can be supposed that $BG_{\eta, fin}$ is a smooth manifold, let us denote
its dimension by $\mu(n).$

Let us denote by $\xi_{fin}$ and $\wt \xi_{fin}$ the
restrictions of the bundles $\xi_\eta \to BG_\eta$ and $\wt \xi_\eta \to
BG_\eta$  to $BG_{\eta, fin}$, the corresponding restriction of the map $\Phi_\eta:\xi_\eta
\to \wt \xi_\eta$  will be denoted by $\Phi_{\eta, fin}: \xi_{\eta, fin} \to
\wt \xi_{\eta, fin}.$

Put $N(n) = 2 \mu(n) + 2n + 2 .$
Let $N$ be any number such that $N \ge N(n)$, and let $j_\eta: \xi_{\eta, fin} \to
\BR^N$ be a smooth embedding. (The number $N(n)$ was chosen so that 
for $N \ge N(n)$ such an embedding exists and it is unique up to isotopy.) 
Put $\Phi^{lift}_{\eta,fin} = (\Phi_{\eta, fin}, j_\eta): \xi_{fin} \to \wt
\xi_{\eta, fin} \times \BR^N.$

This map $\Phi^{lift}_{\eta,fin}$ is a $(\tau, N)$-embedding. 
(The foliation is the
vertical one with leaves parallel to the Euclidean factor.) 

\begin{defi}
The map $\Phi^{lift}_{\eta,fin}$ will be called the {\it {normal form}} of $(\tau,
N)$-embeddings around the stratum $\eta$ and in dimension $\le n.$
\end{defi}

The reason for this name is the following:
\begin{thm}
\label{normform}
Let $e: M^n \hookrightarrow Q^q$ be a $(\tau, N)$-embedding, with $N \ge
N(n).$
Let $\eta(e)$ be the stratum of $\eta$-points of $e$, and $\wt \eta(e) =
e(\eta(e))$ its image. Then there exist tubular neighbourhoods $T$ and $\wt T$
of $\eta(e)$ and $\wt \eta(e)$ respectively, such that $T = e^{-1}(\wt T)$ and 
the
restriction map $e \big|_T: T \to \wt T$ can be induced from the map
$\Phi^{lift}_{\eta,fin},$ i.e. there is a commutative diagram:

\begin{equation}
\label{gnormal}
\begin{CD}
T@> e\big|T >> \wt T\\
@V g VV @VV \wt g V\\
\xi_{\eta,fin} @> \Phi_{\eta,fin}^{lift} >> \wt \xi_{\eta,fin}\times \BR^N
\end{CD}
\end{equation}

Here $g$ and $\wt g$ are fiberwise maps, isomorphic on each fibre.
\end{thm}

\begin{defi} Given two smooth disc bundles $T$ and $\wt T$ over the same smooth
  manifold (possibly with boundary) and a fiberwise map $\mathcal E: T \to \wt T$, we
  shall say that this map is provided with a 
{\it {$G_\eta \times \id_{\BR^N}$-structure}}, if fiberwise maps $g$ and $\wt g$ are given,
 such that diagram~\ref{gnormal} holds (substituting $\mathcal E$ for $e\big|T$) .
Note that a
{$G_\eta \times \id_{\BR^N}$-structure} defines on the target bundle $\wt T$
an $N$-dimensional foliation.
\end{defi}
\begin{rema}
The fiberwise homotopy classes of the maps $g$ and $\wt g$ are
unique.
\end{rema}
\begin{thm}
\label{normformrel}
~(Relative version of the normal form around a stratum) Suppose that the maps $g$ and $\wt g$ are given
over a subset $A \subset \eta (e)$ and they satisfy the commutativity
expressed by diagram \eqref{gnormal}. Then they can be extended to the maps
$g$ and $\wt g.$
\end{thm}
\begin{proof}
Theorems~\ref{normform} and \ref{normformrel} follow immediately
from Proposition~\ref{pro:1a}. Nevertheless we give details in
subsection~\ref{PTbiz}.
  
\end{proof}

\subsection{Global normal form of $(\tau, N)$-embeddings.}
\label{GNF}

Let $\tau$ be a finite sequence of singularities
$\tau = \{[\eta_0] < [\eta_1] \dots < [\eta_r]\}$, and let $n$ be any natural
number.
We shall simplify the notation by writing 
$\xi_i, \xi_{i,fin}, \wt\xi_i, \wt \xi_{i,fin}$ and $j_i$ instead of 
$\xi_{\eta_i},  \xi_{\eta_i, fin},
\wt\xi_{\eta_i}, \wt \xi_{\eta_i,fin}$ and $j_{\eta_i}$ respectively.
(Recall that these objects are:
\begin{enumerate}[$\bullet$]
\item $\xi_i$ - is the universal normal bundle in the source manifold 
of the stratum $\eta_i$
\item $\xi_{i, fin}$ - is its final dimensional (manifold) approximation)
(``good enough'' for $n$-dimensional base spaces). 
\item $\wt \xi_i$ - is the universal normal bundle of the $\eta_i$-stratum in
  the target,
\item $\wt \xi_{i,fin}$ - is its finite dimensional approximation,
\item $j_i$ is the embedding of $\xi_{i,fin}$ into $R^N$, unique up to isotopy.)
\end{enumerate} 

\begin{prop}
 There is a natural number $N(n)$ depending on $n$, such that for
any $N \ge N(n)$ a $(\tau,N)$-embedding $(M^n,e,Q^q, V, \mathcal F)$,\ $(N = q -
n - k)$ 
has the following structure:

There are subsets $T_i,\  \text{for}\ i = 1, 2, \dots, r$ in \ $M^n$
and subsets $\wt T_i,\  \text{for}\ i = 1, 2, \dots, r$ in \ $Q$
such that:
\begin{enumerate}[a)]
\item
$ M^n = T_r \cup T_{r-1} \cup \dots \cup T_0$ and
\item
$ Q \supset \wt T = \wt T_r \cup \wt T_{r-1} \dots \cup \wt T_0,$
where $\wt T$ is a tubular neighbourhood of $e(M^n)$ in $Q,$

\item

$\wt T_i$ is an open set in $Q \setminus \underset{j > i}\cup \wt T_j,$
it is a tubular neighbourhood of $\wt \eta_i(e) \setminus \underset{j > i}\cup
\wt T_j,
$ where the sets $T_{r+1}$ and $\wt T_{r+1}$ are empty.
\item
$T_i = e^{-1}((\wt T_i),$
$T_i$ is an open set in $M^n \setminus \underset{j > i}\cup  T_j,$
it is a tubular neighbourhood of
 $ \eta_i(e) \setminus \underset{j > i}\cup T_j,$
\item
the map
$
e\big|_{T_i}: T_i \to \wt T_i
$
has a $G_i \times \id_{\BR^N}$-structure (i.e. there is a diagram like (\ref{gnormal}).)
\item

{\it Compatibility condition}

There exist $T_i'$ and $\wt T_i'$  slightly bigger neighbourhoods of $T_i$ and
$\wt T_i$ for $i= 1, 2,
\dots , r$  such that $e^{-1}(\wt T_i') = T_i'$
and the restriction maps $e\big|T_i'$  still have $G_i\times \id_{\BR^N}$-structures
extending those of $e\big|T_i.$
Then the leaves of the foliations on the overlapping parts of the enlarged
neighbourhoods $\wt T'_i$ coincide for all indices $i.$
\end{enumerate}
\end{prop}

\begin{proof}
Backward induction on the indices of the strata. Apply first Theorem~\ref{normform}
to the map $T'_r \to \wt T'_r$, and then apply the relative version 
Theorem~\ref{normformrel} to 

$T'_i~\setminus \underset {j > i}\cup T_j \to \wt T'_i
\setminus \underset {j > i}\cup \wt T_j$ modulo the map
$T'_i \cap (\underset {j > i}\cup T'_j) \to \wt T'_i \cap
(\underset {j > i}\cup \wt T'_j).$
\end{proof}

\subsection{Pulling back a $(\tau, N)$-embedding}

The following is an extension of Thom's lemmas, (see \cite{T}, Theorem I.5 and
Lemma IV. 5') its proof is omitted.

\begin{lem}\label{tube}
\begin{enumerate}[\indent a)]
\item
Let $f: V^n \to M^p$ be a $C^{\infty}$ smooth map, $N^{p-m}$ a compact
submanifold in $M^p$, let $T_N$ be a tubular neighbourhood of $N^{p-m}.$ 
Suppose that $f$ is transverse to $N^{p-m}$, and hence the preimage 
$L^{n-p} = f^{-1}(N^{p-m})$ is a submanifold in $V^n$, let us denote its
tubular neighbourhood by $T_L.$ We can suppose that $T_N$ and $T_L$ are chosen so
that $T_L = f^{-1}(T_N).$
Note that $T_L$ and $T_N$ are ball bundles over $L$ and $N$ respectively with
fibres $D^m.$ Let $\pi_L: T_L \to L$ and $\pi_N: T_N \to N$ be the fibering
maps. Let $D'^m \subset D^m$ be a smaller (concentric) ball in $D^m$
and let  $\pi'_L: T'_L \to L$ and $\pi'_N: T'_N \to N$ be the smaller tubular
neighbourhoods  and associated fibre maps with fibre $D'^m.$  

Then there is a diffeomorphism $A$ of $T_N$ onto itself, isotopic to the identity and fixed
near the boundary of $T_N$ such that the map $ A\circ f$
maps the fibres of $\pi'_L$ onto those of $\pi'_N$ diffeomorphically.

\item
(Relative version)
Let $K$be a compact subset of $N$, and $U$ its open neighbourhood in $N.$ Then
if the isotopy $A$ is given on $\pi^{-1}(U)$ 
then it can be chosen so for the whole tubular neighbourhood $T_N$
that it coincide with the given one on $\pi^{-1}(U')$, where $U'$ is
a  neighbourhood of $K$ (possibly smaller than $U).$

\end{enumerate}

\end{lem}

Now let $(M^n, Q^q, e, V,\mathcal F)$ be a $(\tau, N)$-embedding, $(N = 
q - n - k$). Let $\wt T'_i$ be the slightly enlarged sets $\wt T_i$
in $Q^q$, described in
the global normal form of a $(\tau, N)$-embedding and corresponding to the
$(\tau, N)$-embedding $e.$

Let $Q'$be a manifold of dimension $q'$ and
let $\varphi: Q' \to Q^q$ be a map transverse to all the
strata of $e$. Apply the relative version of Lemma~\ref{tube}
by a backward induction to the sets $\wt T_i$, starting with $i=r$ and 
keeping the isotopy constructed in the $i$-th step fixed on the set 
$\wt T'_i \setminus \wt T_i$ during the $i-1$-th step.  
Finally we obtain an
isotopy $A$ on $Q^q$ such that the map
$\varphi^* = A \circ \varphi$ maps diffeomorphically the fibres of the tubular
neighbourhoods of the preimages of the strata to those of the tubular neighbourhoods 
$\wt T_i.$ 
 
Now we obtain a $(\tau, N)$-embedding $(M', Q', e', V',\mathcal F')$
just taking the preimages of the corresponding objects
in $(M^n, Q^q, e, V,\mathcal F).$
In particular the leaves of the foliation $\mathcal F$ are the preimages of 
those of $\mathcal F$ under $\varphi^*.$ This does make sense because $\varphi^*$
maps the fibres of the tubular neighbourhoods diffeomorphically.
This  fact implies also that the local projection along the leaves of the manifold
$M' =(\varphi^*)^{-1}(e(M^n))$ will have only singularities from $\tau.$
The
 cobordism class of the obtained pulled back
$(\tau, N)$-embedding 
$(M', Q', e', V',\mathcal F')$ in 
$\text{\rm Emb}_{(\tau,N)}(Q')$
does not depend on the 
choices made during the construction.
We can summarise the pulling back procedure as follows.

\begin{prop}
$\tau$-embeddings have the following functor-like property.
Let $\varphi: Q' \to Q^q$ be a smooth map, transverse to all the strata of a
given $(\tau,N)$-embedding $(M^n, e, Q^q, V, \mathcal F).$
Then $\varphi^{-1}(M^n)\subset Q'$ will be a $(\tau,N)$-embedding with respect to the
 foliation $\mathcal F' = \varphi^{-1}\circ A^{-1}(\mathcal F)$, where $A$ is the
 diffeomorphism of $Q^q$ described in the application of Lemma~\ref{tube} to
 pulling back a $(\tau,N)$-embedding.
\end{prop}

\begin{defi}
Let $e: M^n \hookrightarrow Q$ be a $\tau$-embedding and $\varphi: Q' \to Q$ be
any map. Let $\varphi'$ be any map homotopic to $\varphi$ and transverse
to all the strata of $e(M^n).$ The cobordism class of $\tau$-embeddings
of the pulled back $\tau$-embedding $(\varphi')^{-1}(M^n)$ is called the
cobordism class \it {associated} to the map $\varphi.$
\end{defi}

\begin{rema}

Clearly homotopic maps $\varphi_0: Q' \to Q$ and $\varphi_1: Q' \to Q$
associate to a given $(\tau,N)$-embedding $e: M^n \hookrightarrow Q$ the same cobordism
class in $Emb_{(\tau,N)}(Q').$

Indeed, approximations of the maps $\varphi_0$ and $\varphi_1$ transverse to
the strata
can be joined by a homotopy transverse to the strata. Now the preimage of
$e(M^n)$ under the homotopy will provide a cobordism of $\tau$-embeddings.
Similarly it can be shown, that the associated cobordism class will be the
same for two $\tau$-embeddings that are cobordant.
 Hence a map arises
$$
[Q', Q^q] \times \Emb_{(\tau,N)}(Q^q) \to \Emb_{(\tau,N)}(Q')
$$
 where $\text{dim}\ Q' - \text{dim}\ Q^q = n' - n$, and
$[Q',Q^q]$ denotes the set of homotopy classes of maps, and $n'$ denotes the
dimension of the source manifold of the $(\tau,N)$ embedding into $Q'$.
\end{rema}
\noindent{
\section{Pontrjagin--Thom construction for $\tau$-embeddings}}

Here we shall show that for any $q$-dimensional manifold $Q^q$ the group

$\Emb^{\text{\rm
Stable}}_{(\tau,N)}(Q^q)$ admits a homotopical description: it is
isomorphic to the group of stable homotopy classes of maps $\{Q^q,
S^NV_\tau\}$ from $Q$ to the $N$-th suspension of a certain
(virtual) space~$V_\tau$. (Here $N = q - n - k.$)
The problem is that we can construct the ``classifying space'' $V_\tau$ only in some
stable sense. For giving the precise meaning of this last sentence we
introduce the notion of the virtual complex.

\subsection{Virtual complexes}

Virtual complexes will be
-- roughly speaking -- CW complexes in which the gluing maps of
the cells are given only as stable maps.
In our case this
notion will arise as a result of a gluing procedure using stable maps.
So we give first the definition of such a gluing procedure using
usual, i.e. non-stable maps. (So the next definition is well-known, we give it
only to make clear the analogy and the differences with the definition
following it.)

\begin{defi}
  Suppose that a sequence of CW pairs $(A_j, B_j) , j =0,  1,2, \dots$ is
  given
  and further  attaching maps $\rho_j$ are defined on $B_j$ with values in
  inductively constructed spaces $X_{j-1}$ (see below), then we can obtain a
  limit space $X = \text{lim}_{j\to \infty} X_j$ by the following well-known 
  glueing procedure:
  \begin{enumerate}[\indent 1)]
    \setcounter{enumi}{-1}
  \item
    $X_0 \overset{def} = A_0.$ Let a map $\varrho_1 : B_1 \to X_0$ be given.
  \item
    $X_1 \overset{def} = A_1 \underset{\varrho_1} 
\cup X_0 = A_1 \sqcup X_0 \Bigm/ b \sim \varrho_1(b), \ \forall b \in B_1. $
i.e. in the disjoint union $A_1 \sqcup X_0$ we identify each point $b \in B_1$
with its image $\rho_1(b) \in X_0.$ 

Let a map
    $\varrho_2 : B_2 \to X_1$ be given.
  \item
    $X_2 \overset{def} = A_2 \underset{\varrho_2} \cup X_1 =
   A_2 \sqcup X_1 \Bigm/ b \sim \varrho_2(b), \ \forall b \in B_2.
 $ Let a map $\varrho_3 :
    B_3 \to X_2$ be given, etc.
  \end{enumerate}

   Let us denote by $c_j$ the connectivity of the pair $(A_j,B_j)$,
   and suppose that the sequence $c_j$ goes to infinity as
  $j \to \infty.$ For simplicity we suppose that this sequence $(c_j)$ of
  connectivities is monotonic.
  Put $X = \displaystyle\lim_{\longrightarrow} X_j.$
\end{defi}

Clearly the
$n$-homotopy type of $X$ is well-defined for any $n.$

\medskip

\begin{exam} Any CW complex having in each dimension finitely many cells is
obtained in this way if for each $j= 1, 2, \dots$ we take $(A_j,B_j) = (D^{c_j+1}, S^{c_j})$
 a closed ball of
dimension $c_j+1$ and its
boundary.
\end{exam}
\begin{defi}[(Virtual complex)]
  Suppose that we are almost in the same situation as before with the only
  difference that {\it {the attaching maps $\rho_j$ are given only in stable sense.}}
  More precisely:
  \begin{enumerate}[\indent 1)]
    \setcounter{enumi}{-1}
  \item
    A sequence of CW pairs $(A_j, B_j) , j =0,  1,2, \dots$ is
    given. Put $X_0 = A_0.$
  \item
    A stable map $\varrho_1: B_1 \not\to X_0$ is given
\footnote{The sign $\not \to$ means ``stable map''.}
, i.e. $\exists
    n_1$ such that a map $S^{n_1}B_1 \to S^{n_1}X_0$ is given. We shall denote
    this map by $S^{n_1} \varrho_1$, although $\varrho_1$ may not exist.
    Put
    $$
    S^{n_1} X_1 \overset{def}= S^{n_1}A_1 \sqcup S^{n_1}X_0 \Bigm/{b \sim
      S^{n_1}\varrho_1(b)}
    $$
    \ \ $\forall b \in S^{n_1} B_1 \subset S^{n_1} A_1$
    Note, that $X_1$ may not exist.
  \item
    A stable map $\varrho_2: B_2 \not\to X_1$ is given, i.e. $\exists
    n_2 > n_1$ such that a map $S^{n_2}\varrho_2 :S^{n_2}B_2 \to S^{n_2}X_1 (= S^{n_2 - n_1}
    (S^{n_1} X_1))$ is given. We shall denote
    this map by $S^{n_2} \varrho_2$, although $\varrho_2$ may not exist.
    Put
    $$
    S^{n_2} X_2 \overset{def} = S^{n_2}A_2 \sqcup S^{n_2}X_1
    \Bigm/{b \sim S^{n_2}\varrho_2(b)}
    $$
    \ \ $\forall b \in S^{n_2} B_2 \subset S^{n_2} A_2$
    Note again, that $X_2$ may not exist.
  \item
    And so on:
$$
S^{n_i} X_i \overset{def} = S^{n_i}A_i \sqcup S^{n_i}X_{i-1}
    \Bigm/{b \sim S^{n_i}\varrho_i(b)}
$$

Here $n_i > n_{i-1}$

  \end{enumerate}
  Again suppose that the connectivity $c_j$ of $(A_j,B_j)$ tends to infinity
  monotonically. Does   the limit space
  $X = \text{lim}\ X_j$
  make sense even if the spaces $X_j$ may not exist ?
  We show that  $X$ exists in the following stable sense.

  For each $m$ there is a natural number $N(m)$ such that for $N> N(m)$
  the $(N+m)$-homotopy type of the space $S^NX$ is well-defined. Indeed, let $i$
  be such that $ c_i > m$, and let $N(m)$ be $ n_{c_i}.$ Then the
  $S^N$ suspensions exist for all the virtual spaces
  $X_0, X_1, \dots, X_i$ and the maps $S^N\varrho_i$ exist for
  $ j = 0, 1, \dots,
  i.$ We claim that the $(N+m)$-th homotopy type of $S^N X_i$ can be taken for
  that of $S^NX$. Indeed,  if $i' > i$, and $N' > N$ is 
   such that the suspension $S^{N'} X_j$ exits for each $ j= 0, 1, \dots,
  i'$,
  then the $ (N' + m)$-th homotopy type of $S^{N'} X_{i'}$ coincides with
  the$(N' - N)$-th suspension of the $(N+m)$-th homotopy type of $S^NX_i .$

  We will say that $X$ exists as {\it {a virtual complex}}, and will denote it
  by $X = \overset {\bullet} \cup_{\varrho_j} \ A_j.$ The dot will
  indicate that the union is
  taken in virtual (i.e. stable) sense.
\end{defi}
\begin{rema} For a virtual complex $X$ its {\it {de-suspension}} can be defined
  as follows: It is a virtual complex $ Y = \text{lim} Y_j$, where
$S^{n_{j+1}}Y_j = S^{n_j}X_j.$ (Note the $Y_j$ may not exist.)
\end{rema} 

\begin{rema}
 Note that $\Gamma X =\Omega^{\infty}S^{\infty} X$ is
  a space well-defined up to homotopy if $X$ is a virtual complex. Indeed we say that
  the  $m$-homotopy type of $\Gamma X$ is that of $\Omega^NS^N X_i$, where $N$
  and $i$ are as above.
\end{rema}

\begin{defi}
\label{skelet}

There is an obvious modification of the above definition when the maps
$S^{n_j} B_j \to S^{n_j} X_{j - 1}$ are defined only on the
($n_j$-suspension of the) $m_j$ skeleton of $B_j$, where $m_j
\geq m_{j - 1}$, $m_j \to \infty$ as $j \to \infty$.
Then the space $\Gamma X$ still can be defined:
\[
\Gamma X = \lim\limits_{j \to \infty} \Omega^{n_j} S^{n_j}
\sk_{m_j} X_j.
\]
\end{defi}

\begin{rema}\label{spect}
It is a natural feeling of the reader that virtual complex is just another
name for the notion of the spectrum. In the present Remark we discuss the relation
of these two notions.
In order to obtain a spectrum from a virtual complex one has to choose
finite dimensional approximations of the blocks $(A_i)$ attached stably to the
previously obtained spaces $(X_{i-1})$.
Two different choices of the sequence of approximations give different
spectra, but they are equivalent in a natural sense, in particular they define
the same extraordinary cohomology theories.
So the advantage of the notion of the virtual complex is that one does not
have to choose data (the sequence of approximations) that finally turn out to
be unimportant.
On the other hand a virtual complex has a filtration, an essential one, that reflects
the steps of its construction ($X_i$ in the definition).
So one can say, that a virtual complex defines an equivalence class of spectra
without pointing out a particular spectrum of this equivalence class.
\end{rema}
 
\subsection{Construction of the virtual complex $V_\tau$}

Recall that $\tau = \bigl\{ [\eta_0] < [\eta_1] < \dots <
[\eta_i]\bigr\}$, where each $[\eta_j]$ is a stable, codimension
$k$ map-germ class for $j = 0, \dots, i$.
They are ordered according to their hierarchy,
(more precisely in a way compatible with the partial ordering given by their hierarchy).
Further let the
set $\tau$ be ``closed'' in the sense that $[\xi] \in \tau$ and
$[\zeta] < [\xi]$ imply $\zeta \in \tau$.
(This assumption implies that
$[\eta_0]$ is the class of germs of maximal rank.)
In Section~2 we introduced for $\eta = \eta_j$ the notations:
$G_j$, $\xi_j$, $\wt\xi_j$, for the objects $G_{\eta _j}$,
$\xi_{\eta _j}$, $\wt\xi_{\eta _j}$.
For $\tau_j = \{\eta_0 < \dots < \eta_j\}$, $j \leq i$, we
denoted the space $\Kaz_{\tau_j}$ simply by $\Kaz_j$.
Recall that $\Kaz_i = \cup \bigl\{ D(\xi_j) \mid j \leq i\bigr\}$
for an appropriate gluing.

In this subsection
we shall show that the target disc bundles $D(\wt\xi_j)$ can be
attached to each other in ``virtual sense'' and give a virtual
complex $V_\tau.$ More precisely we will attach  finite dimensional
approximations of these disc bundles, (see Definition~\ref{skelet}).

For a non-compact space $Q$ we denote by $\overset{\bullet} Q$ the one point
compactification. We put $\overset{\bullet} Q = Q$ if $Q$ is compact.

\begin{thm}\label{PT}
[Construction of the virtual complex $V_\tau$ and the Pontrjagin--Thom construction for $\tau$-embeddings]
  There is a virtual complex $$V_\tau = \overset {\bullet}
\cup \{D(\wt\xi_{\eta})\big|
  \eta \in \tau \}$$ and for each natural number $n$ there is a natural number
  $N(n)$ such that, for  $N \geq N(n)$ and for any manifold $Q$ of dimension $N+n+k$
  the following statement holds:
$$
\Emb_{(\tau,N)}(Q) \approx \{\overset {\bullet} Q, S^N V_\tau\}.
$$

Here the brackets $\{ \,\ \}$ denote stable homotopy classes, i.e. $$\{X,Y\} =
\underset{N \to \infty} {\text{lim}} [S^NX, S^NY].$$
Note that the stable homotopy classes $\{X,Y\}$ make sense also if $X$ is a
space and $Y$ is a virtual complex.
\end{thm}

\begin{add}
The isomorphism holds also if $Q$ is a  manifold with corners or  if it is a
finite simplicial complex.
\end{add}

\begin{add}
(Naturality properties)
For any manifold $Q$ (possibly with corners)
we shall construct homomorphisms
$\alpha:
\Emb_{(\tau,N)}(Q) \to \{\overset {\bullet} Q, S^N V_\tau\}
$
and
$ \beta:
\{\overset {\bullet} Q, S^N V_\tau\} \to \Emb_{(\tau,N)}(Q)
$
that are each others' inverses. Further
the isomorphisms $\alpha$ and $\beta$ commute with the homomorphisms $f^*$ 
induced
by a map $f: (Q', \partial Q') \to (Q, \partial Q)$:
$f$ induces a map $\Emb_{(\tau,N)}(Q) \overset {f^*}\longrightarrow \Emb_{(\tau,N)}(Q')$
by pull-back and a map $\{\overset {\bullet}Q, S^NV_\tau\}\overset {f^*}\longrightarrow
 \{\overset{\bullet}Q', S^NV_\tau\} $ by composition.
\end{add}

The isomorphisms $\alpha$ and $\beta$
will be constructed by an induction on $\tau$ together with the virtual
complex $V_\tau.$

\begin{defi}[of $V_\tau$ and the homomorphism $\alpha.$]
Given a natural number $n$ let $BG_{\eta,fin}$ be such a finite dimensional
approximation of $BG_\eta$, that the pair $(BG_\eta, BG_{\eta,fin})$ is 
$(n+1)$-connected. Let $\mu(n)$ be the dimension of $BG_{\eta,fin}.$
Let $V_{\tau',fin}$ be such a finite dimensional
approximation of $V_\tau'$ that for some $N \geq 2\mu(n) + 2n + 2$ the 
suspension $S^N V_{\tau'}$ exists and the pair $(V_{\tau'}, V_{\tau',fin})$
is $\mu(n)+n+1$-connected.
Then $S^N V_{\tau',fin}$ classifies 
$(\tau',N)$-embeddings with source
dimension not greater than $\mu(n)+n.$
The map 
$$\Phi_{\eta,fin}^{lift}:\xi_{\eta,fin} \to \wt \xi_{\eta,fin}\times R^N$$ is a
$(\tau,N)$-embedding
and let us denote by $\Phi_{S,fin}^{lift}$ or shortly by $\Phi_S$
the associated map of spherical bundles.
More precisely we consider a small spherical bundle 
$S(\wt \xi_{\eta,fin})$ of the vector bundle
$\wt
  \xi_{\eta,fin}$ and its transverse preimage by $\Phi_{\eta,fin}$, and denote
  this preimage by $\mathcal S(\xi_{\eta,fin}).$
The mapping $\Phi_{\eta,fin}^{lift}$ maps $\mathcal S(\xi_{\eta,fin})$ into
$S(\wt \xi_{\eta,fin})\times R^N$. This is the map that we shall denote by $\Phi_S$
and call the associated sphere bundle map.

$\Phi_S$ is a
$(\tau',N)$-embedding with source
dimension not greater than $\mu(n)+n.$ Let us denote by $[\Phi_{S,fin}^{lift}]$
the cobordism class of this $\tau'$-embedding. By the inductional assumption 
a map $\alpha([\Phi_{S,fin}^{lift}]) 
: S\bigl(\wt\xi_{\text{\rm fin}} \times D^N \bigr) 
\to S^N V_{\tau',\text{\rm fin}}$ arises. ($D^N$ is a ball of
sufficiently big radius, so that $\text{image \ of}\ \Phi_S$
is contained in  
$S\bigl(\wt\xi_{\text{\rm fin}} \times D^N \bigr)$
Now put $S^N V_{\tau,\text{\rm fin}} \defo S^N V_{\tau',\text{\rm fin}}
\bigcup
_{\varrho_{\text{\rm fin}}}
D\bigl(\wt\xi_{\text{\rm fin}}\bigr) \times D^N$, where
the attaching map ${\varrho_{\text{\rm fin}}}$ is the map
$\alpha([\Phi_{S,fin}^{lift}]).$ 
This will be a finite-dimensional approximation of the $N$-th
suspension of the virtual space~$V_\tau$.
Taking higher and higher dimensional approximations of
$BG_\eta$ and of the bundles $\xi_\eta$ and $\wt\xi_\eta$ we
obtain better and better approximations of the virtual
complex~$V_\tau$.
\end{defi}
\begin{rema}
\label{univ}

The attaching map
$\rho_{\text{fin}}$
can be chosen so that the subsets 
$$\Phi_{\eta_j,fin}^{lift}(D\xi_{j, fin})
\subset
D\wt \xi_{j,fin} \times D^N$$ for $j = 1, 2, \dots , i$ will be attached 
to each other and they will form a finite dimensional approximation of the
Kazarian space $\Kaz_{\tau} = \cup D\xi_j$, (constructed by the glueing
procedure).
The obtained finite dimensional approximation will be denoted by
 $\Kaz_{\tau, fin} = \cup D\xi_{i,fin}.$
In the process of the proof we will construct the map $\rho_{\text{fin}}$
explicitly. It will also follow that
the inclusion
 $\Kaz_{\tau, fin} \subset S^NV_{\tau, fin}$
restricted to the bundle $D\xi_{i,fin}$
coincides with the map $\Phi_{\eta,fin}^{lift}:D\xi_{i,fin} \to D\wt \xi_{i,fin}
\times \BR^N,$ composed with the map  $D\wt \xi_{i,fin}\times \BR^N \to S^NV_{\tau,fin}.$

\end{rema}

\subsection{Proof of Theorem~\ref{PT} (Pontrjagin - Thom
  construction).} 
\label{PTbiz}

{\bf {Start of the induction}}

Suppose that $i = 0$ i.e.\ $\tau = \{\eta_0\}$.
Recall that $[\eta_0]$ is the class of germs of
maximal rank.
In this case a $\tau$-embedding is a usual embedding $e : M^n
\hookrightarrow Q^{n + k + N}$ with normal vector fields $v_1,
\dots, v_N$ (defined on a neighbourhood of $e(M^n)$).
Choose $V_\tau$ to be the Thom space $M SO(k) = T\gamma^{SO}_k$.
(In this case $V_\tau$ exists as a space, not just as a virtual complex).
The $N$-th suspension $S^N MSO(k)$ is the Thom space of the
bundle $\gamma^{SO}_k \oplus \varepsilon^N$.
Now the usual Pontrjagin--Thom construction provides all the
statements of the Theorem.

{\bf {Induction step}}

Let $\tau$ be equal to $\tau' \cup \{\eta\}$, where $\eta$ is a top
singularity type in~$\tau$.
Suppose that for $\tau'$ we have constructed the virtual complex
$V_{\tau'}$ and all the statements of the Theorem including the  Addenda
are proved
with $\tau'$ substituted for $\tau.$

Let $e : M^n \hookrightarrow Q^q$ be a
$(\tau, N)$-embedding ($q = n+k+N$) and let $\eta(e)$ be the stratum of $\eta$-points of~$e$.
Then $\eta(e)$ is a proper submanifold in~$M^n$
(this means that it is the image of a smooth embedding, it is a closed subset,
and  corners of
codimension $i$ in $\eta(e)$ - and only they -  lie in those of $M^n$).
Let $\wt T$ be a small tubular neighbourhood  of $e(\eta(e))$ in
$Q^q$ and let $T$ be  that of $\eta(e)$ in~$M^n$.
Since the $N$-dimensional foliation $\mathcal F$
is transverse to the stratum $e(\eta(e))$
the tubular neighbourhood $\wt T$ can be
decomposed into the product $\wt T^\perp \times D^N$. The leaves of the
foliation $\mathcal F$ are the sets $x \times D^N,$ for $x \in \wt T^\perp.$ Here $\wt
T^\perp$ can be identified with a
$D^{c + k}$-bundle over $e(\eta(e))$,
its fibres are the images under the exponential map of the small
$\varepsilon$-balls
orthogonal to $V_x \oplus T_x(e(\eta(e))$ in $T_xQ^q$ for $x \in e(\eta(e)).$
(Recall that $c+k$ is the dimension of the target of the root of $\eta.)$
Let us denote by $\pi$ the projection $\wt T \to \wt T^\perp$
along the leaves of $\mathcal F$ (i.e. the projection of $\wt T = \wt T^\perp
\times D^N$ onto $\wt T^\perp$).
Let us consider the map $e^\perp_\eta = \pi \circ e\big|_T : T
\to \wt T^\perp$.
This is a proper $\tau$-map, having $\eta$ points at the
zero-section of the bundle $T \to \eta(e)$.
By  Proposition~\ref{pro:1a}
(about the normal form around a stratum) there is a
commutative diagram
\begin{equation}
\label{perp}
\begin{CD}
T@> e^\perp_\eta >> \wt T^\perp\\
@V g VV @VV \wt g V\\
\xi_{\eta} @> \Phi_{\eta} >> \wt \xi_{\eta}
\end{CD}
\end{equation}

 From now on we will write
simply $\xi$ and $\wt \xi$ omitting the subindex $\eta.$
Using the fact that the bundles $T$ and $\wt T$  have finite-dimensional base spaces of
dimension $\leq n$ , we can
replace the bundles $\xi$ and $\wt\xi$ by their
finite-dimensional approximations $\xi_{\text{\rm fin}}$ and
$\wt\xi_{\text{\rm fin}}$ over a finite-dimensional
approximation of $BG_\eta$ denoted by
$BG_{\eta,\text{\rm fin}}$ so that \eqref{perp} still holds.
Let us denote its dimension by $\mu(n).$
We have seen that 
depending on the dimension $\mu(n)$ of this approximation there
is a natural number $\wh N = \wh N(\mu(n))$ such that there is an
embedding (unique up to isotopy) $\xi_{\text{\rm fin}} \to
D^{\wh N}$ of the total space of the bundle $\xi_{fin}$.
(One can choose for  $\wh N(\mu(n))$ the number
$N(n) = 2 \mu(n) + 2n + 2 .$)

For $N \geq N(n)$, there is an embedding $j :
\xi_{\text{\rm fin}} \to D^N$ (unique up to isotopy).
We  denoted by $\Phi^{\text{\rm lift}}_{\eta,fin}$ the embedding $(\Phi_{\eta,fin},
j): \xi_{\text{\rm fin}} \to \wt\xi_{\text{\rm fin}} \times D^N$.
\[
\gathered
/\hskip-1pt\raisebox{6.85pt}{$\overset{\textstyle e\big|_T}
{\hbox to39mm{\hrulefill}}$}\hskip-.8pt\raisebox{.5pt}{$\searrow$}\quad\ \\[-3mm]
\begin{CD}
T @> e^\perp_\eta >> \wt T^\perp @<
<< \wt T\\
@V g VV @VV \wt g V @VV \ell V\\
\xi_{\text{\rm fin}} @>\Phi_{\eta} >> \wt\xi_{\text{\rm fin}} @<<<
\wt\xi_{\text{\rm fin}} \times D^N
\end{CD}\\
\backslash\hskip-1pt\raisebox{-2.5pt}{$\underset{\textstyle\Phi_{\eta,fin}^{\text{\rm
lift}}}{\hbox to39mm{\hrulefill}}$}\hskip-.7pt\raisebox{-.7pt}{$\nearrow$}\quad
\endgathered \qquad\qquad \parbox{4cm}{There is a map $\ell : \wt T \to
\wt \xi_{\text{\rm fin}} \times D^N$ that makes this
diagram commutative.}
\] 

The composition of the map $\ell : \wt T \to D(\wt\xi) \times D^N$
with the inclusion $D(\wt \xi) \times D^N \to S^N V_{\tau,fin}$ together with the
(stable) map $Q^q \setminus \wt T  \to S^N V_{\tau',fin}$ 
corresponding to the
$\tau'$-embedding $e\big|_{M^n \setminus T} : M^n\setminus T \to Q^q
\setminus \wt T$ (where $T = e^{-1}(\wt T)$), gives a stable map
$\varphi : Q^q  \to S^N V_\tau$.
(The two maps obtained on the common boundary $\partial \wt T =
\partial (Q^q \setminus\text {int}\ \wt T)$ are homotopic, because these are
maps corresponding to the same $(\tau',N)$-embedding  $\partial T \subset
\partial \wt T.$
Then it can be achieved by a deformation that these maps coincide.)

The proof of the statement that stably cobordant $\tau$-embeddings
give stably homotopic maps is quite similar.
Thus the virtual complex $V_\tau$ and the  map
$\alpha$ have been constructed.

Next we construct the map 
$$
\beta:\{\overset{\bullet}Q^q,S^NV_{\tau}\} \to
\Emb^{stable}_{(\tau,N)}(Q^q)
$$

According to Remark~\ref{univ} given a map
$\overset{\bullet}Q^q \to S^NV_{\tau,fin}$ we can
pull back the  $(\tau,N)$-embedding $\Kaz_{\tau,fin} \subset S^NV_{\tau,fin}$ and then we obtain
a $(\tau, N)$-embedding unique up to cobordism and this
cobordism class depends only on the homotopy class of the given map.

In general (when the map
$\overset{\bullet}Q^q \to S^NV_{\tau,fin}$ is given only stably, then)
 we have to replace $Q^q$ and $N$ by $Q^q\times R^L$ and $N+L$ respectively
 for a big enough $L.$

So we obtain a map

$$
\beta_{\fin}:\{\overset{\bullet}Q^q,S^NV_{\tau,fin}\} \to
\Emb^{stable}_{(\tau,N)}(Q^q)
$$
It remained to note that the natural map
$$
\{\overset{\bullet}Q^q, S^NV_{\tau,fin}\} \to
\{\overset{\bullet}Q^q, S^NV_{\tau}\}
$$
is an isomorphism and then $\beta$ is the composition of $\beta_{fin}$ with
this isomorphism.
It is an isomorphism because the pair $(BG_i, BG_{i, fin})$ is
$(n+1)$-connected, hence the pair $(S^NV_{\tau}, S^NV_{\tau,fin})$ is
$(N+n+1)$-connected, so we can apply the generalised Freudenthal theorem.

\bigskip

Clearly the maps $\alpha$
and $\beta$ are each others inverses by construction.

The naturality properties also follow clearly from the constructions 
of the maps.

\subsection{Summary}

We have shown that given manifolds $P^{n+k}$ and $Q^q$ we have:
\begin{enumerate}[\indent a)]
\item $\Cob_\tau(P^{n+k}) = \Emb_{(\tau,N)}(P^{n+k} \times \BR^N)$ for
  any sufficiently large $N$.
\item For any sufficiently large $N$ it holds:
$\Emb_{(\tau,N)}^{\text{\rm Stable}}(P^{n+k}) \approx \Emb_{(\tau,N)} (P^{n+k} \times \BR^N)$
\item $\Emb_{(\tau,N)}^{\text{\rm Stable}}
(Q^q) = \{\overset {\bullet}Q^q, S^N V_\tau\}$, for any natural number $N$,
where
  $V_\tau = \bigcup\limits^\bullet \bigl\{D\wt\xi_\eta \mid [\eta]
  \in \tau \bigr\}$ is the  virtual complex described above.
\end{enumerate}
\begin{cor}
$\Cob_\tau(P^{n+k}) \approx \bigl\{ \overset{\bullet}{P^{n+k}}, V_\tau\bigr\}$.
\end{cor}
\begin{rema}\
  \begin{enumerate}[1)]
  \item
    Recall that in \cite{R--Sz} a space $X_\tau$ was constructed
    such that $\Cob_\tau(P^{n+k}) = \bigl[ \overset{\bullet}{P^{n+k}}, X_\tau\bigr]$.
  \item
    We have seen that the functors $\Cob_\tau(-)$ and $\Emb_\tau(- \times \BR^N)$ can
    be extended from manifolds to simplicial complexes.
    (See Definition~\ref{simpl} and also \cite{Sz3}).
    The Corollary above will hold for any simplicial complex $P$.
  \end{enumerate}
\end{rema}
We obtain the following
\begin{prop}\label{X=GV}
$X_\tau \cong \Gamma V_\tau$, where $\Gamma(\ ) = \Omega^\infty
S^\infty(\ )$.
\end{prop}
\begin{proof}
\[
\aligned
\Cob_\tau (P) &\approx \bigl\{ \overset{\bullet}{P},
V_\tau\bigr\} = \bigl[ \overset{\bullet}{P}, \Gamma V_\tau\bigr]\\
\text{and } \ \Cob_\tau(P) &= \bigl[\overset{\bullet}{P}, X_\tau\bigr].
\endaligned
\]
Hence both the space $X_\tau$ and the space $\Gamma V_\tau$ are
classifying spaces for the functor $\Cob_\tau(-)$.
But by Brown's theorem the classifying space is homotopically unique.
\end{proof}

\part{The Kazarian conjecture}
\section{Formulation of the conjecture}

Maxim Kazarian formulated a conjecture relating the space
$\Kaz_\tau$ (which we called Kazarian space) with the space
$X_\tau$, classifying $\tau$-cobordisms.

Below $\cong_\Q$ means
``rationally homotopically equivalent'' and $\Kaz^+_\tau$ denotes the disjoint union of
$\Kaz_\tau$ with a point.

\begin{conj}[Kazarian, \cite{K2}]
$$X_\tau \cong_\Q \Gamma S^k(\Kaz^+_\tau).$$
\end{conj}
This means that $X_\tau$ gives the stable rational homotopy type of
the $k$-th suspension of $\Kaz_\tau^+.$
Here we shall formulate and prove a statement giving a
relationship between $X_\tau$ and $\Kaz_\tau$ that gives a homotopy
equivalence -- not just a rational one -- after modifying the
right-hand side properly.
Notice that $S^k(\Kaz_\tau^+)$ is the Thom space of the trivial
$k$-dimensional
vector bundle over $\Kaz_\tau.$
Roughly speaking we replace the $k$-th suspension
$S^k(\Kaz^+_\tau)$ by the Thom space of the $k$-dimensional vector
bundle $\nu^k$ over $\Kaz_\tau$, see Definition~\ref{nu}.
But actually $\nu^k$ is only a virtual vector bundle.
\begin{rema}
Recall that a {\it {virtual vector bundle}} $(\alpha, \beta)$ over a
space $A$ is an equivalence class of formal differences of
vector bundles. By dimension of such a virtual vector bundle $(\alpha, \beta)$
one means the difference
$\dim(\alpha, \beta) = \dim \text{fiber of } \alpha - \dim
\text{fiber of } \beta$.
If $A$ is compact, then any virtual vector bundle admits a
representation of the form $(\alpha, \varepsilon^N)$, where
$\varepsilon^N$ is the trivial $N$-bundle.
If $A$ is a union of compact subspaces $A_1 \subset A_2 \subset
\dots$, $A = \bigcup A_i$, then a virtual bundle of
dimension~$k$ over $A$ is represented by a sequence $\bigl(\alpha_i,
\varepsilon^{N_i}\bigr)$ of pairs, where $\alpha_i$ is a vector
bundle over $A_i$ of dimension $N_i + k$, \ $\varepsilon^{N_i}$ is
the trivial $N_i$-dimensional bundle and $\alpha_i\big|_{A_{i -
1}} = \alpha_{i - 1} \oplus \varepsilon^{{N_i} - N_{i - 1}}$.
\end{rema}
\begin{defi}
By Thom space of a virtual bundle $(\alpha, \varepsilon^N)$ we
mean the virtual complex $S^{(-N)} T \alpha$.
Here $T\alpha$ is the Thom space of the bundle $\alpha$, and
$S^{(-N)}$ is the inverse of the $N$-th suspension.
Note that the $(-N)$-th suspension is a meaningful
functor from virtual complexes into virtual complexes.
\end{defi}
\begin{prop}\label{V}
There exists a virtual $k$-dimensional vector bundle $\nu^k$
over $\Kaz_\tau$ such that $\Gamma V_\tau {\cong} \Gamma T\nu^k$
\footnote{$\cong$ stands for ``homotopy equivalent''}
\end{prop}

The next corollary is a refinement of Kazarian's conjecture.

\begin{cor}
\label{refine}
$X_\tau \cong \Gamma T\nu^k$.\hfill $\square$
\begin{proof} It follows immediately from Propositions~\ref{V} and
\ref{X=GV}. 
\end{proof}
\end{cor}
\begin{proof} of Proposition~\ref{V}.
We have seen that the bundles $D\wt\xi_{\eta, \text{\rm fin}}
\times D^N$ attached properly to each other give the space $S^N
V_{\tau,\text{\rm fin}}$.
Here the subindex fin indicates that we consider certain
finite-dimensional approximations.
From now on we will omit this subindex.
The attaching maps are defined by the embeddings $D\xi_\eta
\subset D\wt\xi_\eta \times D^N$.
When we attach the bundle $D\wt\xi_\eta \times D^N$ to $S^N
V_{\tau'}$, then the bundle $D\xi_\eta$ is attached to
$\Kaz_{\tau'} \subset S^N V_{\tau'}$.
So we obtain that $\Kaz_\tau \subset S^N V_\tau$.
Recall that $\nu^k$ is the virtual vector bundle over $\Kaz_\tau$
which is defined by the following restrictions:
$\nu^k \big|_{D\xi_\eta} = \pi^* \wt\xi_\eta - \pi^* \xi_\eta$,
where $\pi : D\xi_\eta \to BG_\eta$ is the disc-bundle of the
bundle~$\xi_\eta$.
Note that for those finite-dimensional approximations where the embedding
$D\xi_{\eta,\text{\rm fin}} \hookrightarrow D
\wt\xi_{\eta,\text{\rm fin}}$
$\times D^N$ exists the virtual
bundle $\nu^k \oplus \varepsilon^N$ can be realized as a
(non-virtual) vector bundle (it can be regarded as the ``normal
bundle'' of the embedding $\Kaz_{\tau, \text{\rm fin}}
\hookrightarrow S^N V_{\tau,\fin}$).
Both the space $\Kaz_\tau$ and the virtual complex $V_\tau$ are
stratified by their constructions: 
$$\Kaz_\tau = \bigcup \bigl\{ D
\xi_{\eta_i} \mid [\eta_i] \in \tau \bigr\}; \  \  \ V_\tau =
\overset{\bullet}{\bigcup} \bigl\{D \wt\xi_{\eta_i} \mid [\eta_i] \in
\tau \bigr\}.$$
The inclusion
$\Kaz_\tau \subset S^N V_\tau$ respects this stratification:
\[
(\Kaz_\tau, \Kaz_{\tau'}) \subset (S^N V_\tau, S^N V_{\tau'}).
\]
Hence it induces  map of the quotient spaces
\[
i_{\eta}: T\xi_\eta \hookrightarrow S^N T\wt \xi_\eta.
\]
(Here $\eta$ is a highest singularity occurring in $\tau$-maps with at most
$n$-dimensional source manifolds, and  $\tau'$ is the set of singularities in
$\tau$, strictly lower  than $\eta$.)
The base space of both bundles $\xi_\eta$ and $\wt \xi_\eta$ is the space
$BG\eta$ and $i_{\eta}$ is the identity on it.

Let $U$ be a
tubular neighbourhood of $\Kaz_\tau$ in $S^N V_\tau$. Then the
map
\begin{equation}
\label{quotient}
S^N V_\tau \to S^N V_\tau / (S^N V_\tau \setminus U) = S^N T\nu^k
\end{equation}
induces maps of the quotients
$$
S^NV_\tau/S^NV_{\tau'} =
S^N T\wt\xi_\eta = T(\wt\xi_{\eta} \oplus \varepsilon^N)\to$$
$$ \to S^N T\wt\xi_\eta \big/ (S^N T\wt\xi
\setminus U) = S^N T \bigl(\nu^k \big|_{D\xi_\eta}\bigr) \big/
S^N T(\nu^k \big|_{\partial D\xi_\eta}\bigr) = T(\nu^k \oplus \xi_\eta \oplus \varepsilon^N).
$$
These maps map the zero section of $\wt \xi_\eta \oplus \varepsilon^N$
into that of $\nu \oplus \xi_\eta \oplus \varepsilon^N$ by the
identity homeomorphism.
Therefore they induce isomorphisms in the cohomologies.
Hence the map \eqref{quotient} also induces isomorphism on cohomologies.
\end{proof}
\begin{rema}

Note that in (4) we actually have some finite dimensional approximation
$V_{\tau, fin}$ of the virtual complex $V_\tau$ and the natural number $N$
is chosen so that 
\begin{enumerate}[a)]

\item $S^N V_{\tau,fin}$ exists as a (non-virtual) space,

\item over the corresponding finite dimensional approximation of $\Kaz_\tau$ the virtual
bundle $\nu^k$ has the form $\alpha^{N+k} - \varepsilon^N$.
\end{enumerate}

 So the Thom
space $S^N T\nu^k$ (as a virtual complex) has the corresponding finite dimensional
approximation $T\alpha^{N+k}.$ 
We have shown that $S^NV_{\tau,fin} = T\alpha^{N+k}.$
Therefore the virtual complexes $V_\tau$ and $T\nu^k$ have sequences of finite
dimensional approximations and natural numbers $N_j$
such that the $N_j$-th suspensions of the $j$-th approximations exist and are
homotopically equivalent. We express this  by saying that {\it {the 
virtual complexes 
$V_\tau$ and $T\nu^k$ are homotopically equivalent}} and write: 
$V_\tau \overset{\bullet}{\cong} T\nu^k$.
In particular $\Gamma V_\tau \cong \Gamma T\nu^k.$
\end{rema}
\begin{nota}
For any space $A$ let us denote by $SP\,A$ the infinite symmetric
product of~$A$.
\end{nota}
\begin{rema}
It is known, (see e.g.\ \cite{H} page 472) that for any base pointed, connected $CW$-complex $A$ the
space $SP\,A$ is (homotopy equivalent to) a product of Eilenberg--MacLane spaces.
Since by the Dold--Thom theorem $\pi_i(SP\,A) = H_i(A)$, (see [H]), we obtain
\begin{equation}
\label{SP}
SP\,A = \prod\limits_i K(H_i(A), i).
\end{equation}
\end{rema}
\begin{lem}
\label{homol}
Let $A$ and $B$ be arbitrary spaces such that $H_*(A) \approx H_*(B)$.
Then $SP\,A$ is weakly homotopy equivalent to~$SP\,B$.
\end{lem}
\begin{proof}
Obvious from \eqref{SP}.
\end{proof}

\begin{rema}
We do not require that $A$ and $B$ were simply connected,
neither that the isomorphism $H_*(A) \approx H_*(B)$ were
induced by a map of the spaces.
\end{rema}

\begin{rema}\
  \begin{enumerate}[1)]
  \item
    In particular $\Omega (SP (SA)) = SP A.$
  \item
    A co-fibration $B \subset A \to A/B$ gives a quasi-fibration 
   (see \cite{H}) of the
   corresponding infinite
   symmetric products $SPB \to SPA \to SP(A/B)$ for which the homotopy exact
   sequence holds. Hence if the pair $(A,B)$ is $\ell$-connected, then the
   pair $(SPA, SPB)$ is also at least $\ell$-connected.
 \item
   Therefore the functor $SP(\ )$ can be defined for a virtual complex $V$ 
   as well as
   follows: 
Let $(m_j, n_j)$ be a sequence of pairs of natural numbers from
Definition~\ref{skelet}, i.e. $S^{n_j} sk_{m_j}V$ exists, 
$n_j \to \infty, m_j \to \infty.$ Now put
 $$SP(V) =\underset {j \to \infty} {\text {lim}} \Omega^{n_j} SP S^{n_j} sk_{m_j}V.$$
  
 \end{enumerate}
\end{rema}

\section{Proof of the Kazarian conjecture}

By the Thom isomorphism we have
\[
\overline H^{* + k} (T\nu^k) \approx H^*(\Kaz_\tau) \approx
\overline H^{* + k} \bigl(S^k (\Kaz^+_\tau)\bigr).
\]
We obtain the following equalities:
\[
X_\tau = \Gamma V_\tau = \Gamma T\nu^k \ \text{ and } \
SP\,T \nu^k = SP\bigl(S^k(\Kaz^+_\tau)\bigr).
\]
Here equality means homotopy equivalence.
\begin{lem}
\label{Hurew}
For any (virtual) complex $A$ there is a natural map $h^A : \Gamma
A \to SP\,A$ that induces a rational homotopy equivalence.
Moreover the map $h^A_\#$ induced by $h^A$ in the homotopy groups is
the stable Hurewicz map $h^A_\# : \pi_i (\Gamma A) =
\pi_i^s(A) \to H_i(A) = \pi_i(SP\,A)$.
\end{lem}
\begin{proof}
In \cite{B--E} a model has been constructed for the
space $\Gamma A$ for any connected complex~$A$:
\[
\Gamma A = \coprod_i W S(i) \underset{S(i)}{\times} (A \times
\dots \times A) \Bigm/ \sim
\]
where $W S(i)$ is a contractible space with a free $S(i)$-action
on it, and $\sim$ denotes gluing by some natural equivalences.
($S(i)$ is the symmetric group on $i$ elements.)
The projections 
$$WS(i) \underset{S(i)}{\times} (A\times \dots
\times A) \to (A \times \dots \times A) \Bigm/S(i)$$ are
consistent for different $i$-s with the gluing in
$\Gamma(A)$ and $SP(A)$.
So they define a map $h : \Gamma(A) \to SP(A)$ and the induced
map $h^A_\# : \pi _i(\Gamma(A)) = \pi^s_i(A) \to \pi^s_i(SP(A))
= H_i(A)$ is the stable Hurewicz map.
By Serre's theorem this is a rational isomorphism (see Theorem 18.3 in \cite{M--S}). This proves
the lemma when $A$ is a non-virtual space.
If $A$ is a virtual complex, then we replace it by its
appropriate suspension (over each finite skeleton) and repeat
the previous proof using that $\Omega^N \Gamma S^N A = \Gamma A$
and $\Omega^N SP\, S^N A = SP\,A$.
\end{proof}
Now we finish the proof of the Kazarian conjecture.

$$
X_\tau(k) \overset{(1)}{\cong} \Gamma T\nu^k \overset{(2)}{\cong}_\Q SPT\nu^k
\overset{(3)}{=} SPT\varepsilon^k \overset{(4)}{=} SPS^k(\Kaz_\tau^+)
\overset{(5)}{\cong}_\Q \Gamma S^k(\Kaz_\tau^+)
$$

Here

(1) is Corollary~\ref{refine},

(2) and (5) follow from Lemma~\ref{Hurew},

(3) comes from Lemma~\ref{homol},

(4) is obvious.

\hfill $\square$

\begin{rema}
It follows from the above proof that the following modification of the
Kazarian conjecture also holds:
\footnote{Kazarian recently has told me that this was the first version of his conjecture.}
$$ X_\tau \cong  _\Q SP(S^k(\Kaz_\tau^+)).$$
\end{rema}

\section{Corollaries of the Kazarian conjecture}

Here we compute the rational cobordism groups of $\tau$-maps.
Further we define characteristic classes of $\tau$-maps,
show analogues of the Pontrjagin--Thom theorem and those of the Conner--Floyd
theorem claiming that the characteristic numbers completely
determine the cobordism (respectively bordism) class.

\begin{prop}
\label{pro:0}
$\Cob_\tau(n, k) \otimes \Q \approx H_n(\Kaz_\tau; \Q)$.
\end{prop}
\begin{proof}
$\Cob_\tau(n, k) \otimes \Q = \pi_{n + k} (X_\tau) \otimes \Q =
\pi_{n + k} (\Gamma V_\tau) \otimes \Q = \pi^s_{n + k}
(V_\tau) \otimes \Q = \pi^s_{n + k} (T\nu^k) \otimes \Q
\overset{h_\Q}{\longrightarrow} H_{n + k} (T\nu^k; \Q) =
H_n(\Kaz_\tau; \Q)$.
Here $h_\Q$ is the tensor product of the stable Hurewicz homomorphism
with the identity map
$\id_\Q : \Q \to \Q$, hence it is an isomorphism.
\end{proof}
So the stable homotopy groups of the virtual complex $V_\tau$
give the $\tau$-cobordism groups, while the homology groups of
$V_\tau$ are isomorphic to those of $\Kaz_\tau$ after a shift of
the dimensions by~$k$.
Now we are going to give a more explicit description of the arising
map $\Cob_\tau(n,k) \to H_n(\Kaz_\tau)$, denote this map by~$\varphi$.
\begin{prop}
\label{pro:1b}
Let $f : M^n \to \BR^{n + k}$ be a $\tau$-map, and let $\kappa_f : M^n
\to \Kaz_\tau$ be the corresponding Kazarian map, lifting the stable normal
map $\nu_f : M^n \to BSO.$
Then $\varphi([f]) = (\kappa_f)_* [M]$.
\end{prop}
\begin{proof} is standard, see \cite{Sz10}, page 34. \end{proof}

\begin{rema}
This Proposition can be considered as a (rational) analogue for
$\tau$-maps of the Pontrjagin--Thom theorem saying that two
manifolds are cobordant if and only if their characteristic
numbers coincide.
To make the analogy  more transparent we give the definition
of characteristic numbers of a $\tau$-map and reformulate the
Proposition.
\end{rema}
\begin{defi}
If $f : M^n \to \BR^{n + k}$ is a $\tau$-map, and $x \in H^n
(\Kaz_\tau)$ is any $n$-dimensional cohomology class of $\Kaz_\tau$,
then the number $\left<\kappa^*_f(x), [M^n]\right>$ will be called
the $x$-character\-is\-tic number of~$f$.
\end{defi}
\begin{prop}
There is an isomorphism
\[
\aligned
\Cob_\tau (n,k) \otimes \Q
&\overset{\approx}{\longrightarrow}
\Hom \bigl(H^n(\Kaz_\tau);\Q\bigr)\\
\text{given by the formula }
\psi: [f]\  &\longmapsto \bigl(x \longmapsto \left< \kappa^*_f(x), [M^n]\right>\bigr).
\endaligned
\]
\end{prop}
\begin{proof}
$\bigl<\kappa^*_f(x), [M^n]\bigr> = \left<x, \kappa_{f*}[M^n]\right>$.
The set of these numbers for $\forall x \in H^n(\Kaz_\tau)$ define the class
$\kappa_{f*}[M^n]$. Hence the map $\psi$ is injective.
Its source and target have equal dimensions over $Q,$ hence $\psi$ is an isomorphism.

\end{proof}
\begin{exam}
\label{exa:1a}
Let $\tau$ be the set of all possible singularities of codimension $k$ germs.
Then $\Cob_\tau(n,k) = \Omega_n(\BR^{n + k}) = \Omega_n$,
where $\Omega_n(X)$ denotes the $n$-th bordism group of a space
$X$ and $\Omega_n$ is the cobordism group of oriented
$n$-dimensional manifolds.
In this case $\Kaz_\tau \cong B\,SO = \underset {N \to \infty} {\text {lim}} B\,SO(N)$.
(The fibre of the map $\Kaz_{\tau} \to BSO$ is the space of all
polynomial maps
and that is contractible.)
So in this case our proposition says:
\[
\Omega_n \otimes \Q \approx \Hom(H^n(B\,SO);\Q) \approx
H_n(B\,SO;\Q),
\]
and that is a well-known theorem of Thom.
\footnote{This case is not quite covered by our consideration, since we
  allowed in $\tau$ only stable singularities.}
\end{exam}
\begin{exam}
\label{exa:2a}
Let $\tau$ be equal to  $\{\Sigma^0\}$ i.e.\ in this case a $\tau$-map is the same as
 an immersion. Then $\Kaz_\tau$ is the $V_N(\BR^{N + k})$ bundle
associated to the bundle $\gamma^{SO}_{N + k} \to B\, SO(N + k)$,
for $N \to \infty$.

But this is just $B\,SO(k)$.
Hence the Proposition says that the rational cobordism group of
oriented immersions $\Cob_\tau(n, k) \otimes \Q = \text{\rm
Imm}^{SO}(n, k) \otimes \Q$ is isomorphic to $H_n(B\, SO(k);\Q)$.
(Compare with \cite{B}.)
\end{exam}
Next we consider the case when the target is an arbitrary
oriented $(n + k)$-manifold~$P^{n + k}$.
\begin{prop}
\label{pro:1c}
$\Cob_\tau(P^{n + k}) \otimes \Q \approx \bigoplus\limits_j H^j
(P^{n+k}; H_{j - k} \bigl(\Kaz_\tau; \Q)\bigr)$.
\end{prop}
\begin{proof}
We will use the sign ``$A \underset{ \Q}{\approx} B$'' meaning
that $A \otimes  \Q \approx B \otimes  \Q$, where $A$ and $B$ are Abelian groups.
(Further we recall that for any space $X$ it holds $SP\, X =
\prod\limits_j K(H_j(X), j)$ and there is a weak rational
homotopy equivalence $h^X: \Gamma X \to SP\, X$, where $\Gamma X =
\Omega^\infty S^\infty X$.)
\[
\Cob_\tau(P^{n + k}) = [P^{n+k}, X_\tau] = [P^{n+k}, \Gamma V_\tau] = [P^{n+k},
\Gamma T\nu^k]
\underset{\overset{\approx}{\Q}}{\overset{h^{T\nu^k}}{\longrightarrow}}
[P^{n+k}, SP\, T\nu^k] =
\]
\[
= \Bigl[ P^{n+k}, \prod\limits_j K\bigl(H_j (T\nu^k), j\bigr)\Bigr] =
\Bigl[ P^{n+k}, \prod\limits_j K \bigl(H_{j - k} (\Kaz_\tau), j\bigr)
\Bigr] =
\]
\[= \bigoplus\limits_j H^j \bigl(P^{n+k}; H_{j - k}
(\Kaz_\tau)\bigr).\] \end{proof}
\begin{defi}
We will call the cohomology classes in $H^*(\Kaz_\tau)$ 
{\it {characteristic classes of $\tau$-maps.}}
If $f$ is a $\tau$-map, $\kappa_f : M^n \to \Kaz_\tau$ is the
corresponding Kazarian map and $x \in H^*(\Kaz_\tau)$, then
$\kappa^*_f(x)$ will be called the {\it {$x$-characteristic class of the
$\tau$-map~$f$.}}
\end{defi}

\begin{defi}\label{char}
Let $f : M^n \to P^{n + k}$ be a $\tau$-map, representing a
cobordism class $[f] \in \Cob_\tau(P^{n + k})$.
Let $x$ be a cohomology class in $H^*(\Kaz_\tau)$ and $z \in H^*(P^{n+k})$.
Then the number $\bigl< \kappa^*_f(x) \cup f^*(z), [M^n]\bigr>$ will be
called the {\it {$(x, z)$-characteristic number of the $\tau$-map $f$.}}
\end{defi}
\begin{prop}
Let $f : M^n \to P^{n + k}$ be a $\tau$-map.
Then the $\tau$-characteristic numbers $\left<\kappa^*_f(x) \cup
f^*(z), [M^n] \right>$ determine the rational $\tau$-cobordism
class of $f$, i.e.\ $[f] \otimes 1_ {\Q} \in \Cob_\tau(P^{n + k})
\otimes  \Q$, and they are well-defined invariants of this class.
Moreover, if for any $x \in H^*(\Kaz_\tau)$ and $z \in H^*(P^{n+k})$ a
rational number $r(x, z)$ is given, then there is a unique
element in $\Cob_\tau(P^{n+k}) \otimes  \Q$ such that the numbers $r(x,
z)$ give the $(x, z)$-characteristic numbers of this element.
\end{prop}

\begin{rema}
Following the suggestion of the referee I omit the proof here (but refer to
\cite{Sz10}, page 37, where one can
find also an application of these characteristic numbers, namely the {\bf
ring structure} on the direct sum of all
the cobordism groups of fold maps, sum for all dimensions and
all codimensions, has been computed, \cite{Sz10}, Part IX.).
This will be published elsewhere, see \cite{L--Sz}.
\end{rema}

\part{Concrete computations}
Here we compute explicitly the rational cobordism
groups of
\begin{enumerate}[\indent a)]
\item
  Morin maps
    with trivial (and trivialised) kernel bundle
  (these are the so called prim maps),
\item
  arbitrary Morin maps,
  \item
  Morin maps with at most $\Sigma^{1_r}$ singularities

\end{enumerate}
 
\section {Morin maps.}

Let $\tau$ be the set of all stable singularity classes of corank $\le 1,$
i.e. $\tau = \{ \Sigma^0, \Sigma^{1,0}, \dots, \Sigma^{1_i,0}, \dots\}.$

Here $1_i$ stands for $i$ digits $1,$ and any Thom- Boardman symbol stands
for the corresponding stable singularity class.
Such a $\tau$-map is called a Morin map. For this $\tau$ the group
$\Cob_{\tau}(n,k)$ will be denoted by $\Cob_{Morin}(n,k).$
In this section we prove the following:

\begin{thm}
\label{Cob-Morin}
$
\Cob_{Morin}(n,k)\otimes \Q =
\begin{cases} H_n(BO(k);\Q)\ \text{if} \ k \ \text{is \ even} \\
H_n(BO(k+1);\Q)\ \text{if} \ k \ \text{is \ odd}
\end{cases}$

\end{thm}

\begin{rema}

Recall that we always consider \underline{oriented} cobordism group of
singular maps. Surprisingly these are related to the \underline{unoriented}
Grasmann manifolds.

\end{rema}

\begin{lem}\label{HK-Morin}
$\Kaz_{Morin}(k) \cong_\Q \begin{cases} BO(k)\ \text{if} \ k \ \text{is \ even} \\
BO(k+1)\ \text{if} \ k \ \text{is \ odd}
\end{cases}$
\end{lem}

Since $\Cob_{\tau}(n,k)\otimes \Q \approx H_n(\Kaz_\tau;\Q)$
 Lemma~\ref{HK-Morin} implies
Theorem~\ref{Cob-Morin}.
The proof of the Lemma consists of computing Kazarian's spectral
  sequence for Morin maps. The main tool will be the notion of ``prim maps''.

\begin{defi} We say that a Morin map $f:M^n \to P^{n+k}$ is a {\it {prim map}}
(the word prim is an abbreviation for {\it {pr}}ojected {\it{im}}mersion) if the line bundle
over the singularity set formed by the kernels of the differential is trivial
and trivialised. Such a map always can be decomposed as $\pi \circ g$, where
$g: M^n \to P^{n+k} \times \BR^1$ is an immersion, and $\pi: P^{n+k} \times \BR^1
\to \BR^1$ is the projection. Here the immersion $g$ is unique up to regular
homotopy (if we require that $\frac{\partial (\pi \circ g)}{\partial \underline
  v} > 0$, where $\underline v$ is the positive direction of $\text{ker}df$).
\end{defi}

\begin{rema}~\label{involution}
Let us denote by $\Kaz_{prim}(k)$ the Kazarian space for all prim maps of
codimension $k.$ Alike to the decomposition
$\Kaz_{Morin}(k) =\cup_i D\xi_i$, where
$\xi_i$ is the universal normal bundle in the source for the singularity
stratum $\Sigma^{1_i}$, there is a decomposition
$\Kaz_{prim}(k) = \cup_i D\bar \xi_i$, where $D\bar \xi_i \to D\xi_i$ is a double
covering for each $i > 0,$
\footnote{This double cover is associated to the kernel line bundle of a
  Morin map. Note that for a prim map the kernel line bundle is 
always trivial.} and it is the identity of $BSO(k) = D\bar \xi_0 =
D\xi_0$ for $i= 0.$ Hence there is a $\Z_2$ action on $\Kaz_{prim}(k)$ such
that the quotient is $\Kaz_{Morin}(k).$

Prim maps are easy to handle because the automorphism groups (in the case of
oriented source and target manifolds) are all isomorphic to $SO(k)$, and
the corresponding bundles $\bar \xi_i$ (the universal normal bundle of the
$\Sigma^{1_i}$-stratum in the source) are the bundles $i\cdot \gamma_k \oplus
\varepsilon^i$, where $\varepsilon^i$ is the $i$-dimensional trivial bundle. 
\end{rema}

\begin{lem}
\label{lem:2b}
Let $\ol K$ be a CW complex with a cellular $\Z_2$-action and let
$K$ be the quotient $\ol K/\Z_2$, and let $p : \ol K \to K$ be
the quotient map.
Suppose that $K$ is simply connected.
Then $H^*(K; \Q)$ is equal to  $H^*(\ol K; \Q)^{\Z_2} $ the $\Z_2$-invariant
part of $H^*(\ol K; \Q),$ and the homomorphism $p^* : H^*(K; \Q) \to H^*(\ol K; \Q)$
is  injective and its image is $H^*(\ol K; \Q)^{\Z_2} $.
\end{lem}
\begin{proof}
Let us consider the map $q : \ol K \underset{\Z_2}{\times}
S^\infty \to K = \ol K/\Z_2$.
Note that $q^{-1}(x) = S^{\infty}$ if $x$ is not a fix point of
the $\Z_2$-action, and $q^{-1} (x) = RP^\infty$ otherwise.
Recall that there exists a Leray spectral sequence (see \cite{Go}):
\[
E^{i,j}_2 = H^i\left(K; \mathcal H^j\left(q^{-1}(x); \Q\right)\right)
\Longrightarrow H^* \bigl(\ol K \underset{\Z_2}{\times} S^\infty;
\Q \bigr).
\]
Here $\mathcal H^j\left(q^{-1}(
x); \Q\right)$ is the group of local
coefficients at the point $x \in K.$
Since $H^j\bigl(q^{-1}(x); \Q\bigr) = 0$ for $j > 0$ and
$H^0\bigl(q^{-1}(x); \Q
\bigr) = \Q$ we have that
\[
\aligned
E^{i,j}_2 &= 0 \ \text{ for } \ j > 0,\\
 \text{and } \ E^{i,0}_2 &= H^i(K ;\Q).
\endaligned
\]
So the Leray spectral sequence degenerates and gives that
\begin{equation}
\label{Leray}
H^i \Bigl( \ol K \underset{\Z_2}{\times} S^\infty; \Q) = H^i (K; \Q).
\end{equation}
On the other hand, there is a fibration $\ol K
\underset{\Z_2}{\times} S^\infty \overset{\ol K}{\longrightarrow}
RP^\infty$, and now we consider the spectral sequence of this
fibration:
$E^{i,j}_2 = H^i\bigl(RP^{\infty}; \mathcal H^j(\ol K; \Q)\bigr)$
with twisted coefficients.
The fundamental group $\pi_1(RP^\infty) = \Z_2$ acts on $\mathcal
H^*(\ol K; \Q)$.
Let us denote by $\mathcal H_+$ and $\mathcal H_-$, respectively
the eigenspaces for the eigenvalues $+1$ and $-1$,
respectively of the $\Z_2$-action.
It is well known that
\[
\aligned
H^*(RP^\infty; \mathcal H^*_-) &= 0\\
\text{and } \
H^i(RP^\infty; \mathcal H^j_+) &= \begin{cases}
\mathcal H^j_+ &\text{if } i = 0\\
0 &\text{if } i \neq 0.
\end{cases}
\endaligned
\]
Hence $E^{i,j}_2 = \begin{cases}
0 &\text{if } i \neq 0 \\
H^j(\ol K; \Q)^{\Z_2} &\text{if } i = 0.
\end{cases}$.

Hence
\begin{equation}
\label{invol}
H^j\Bigl(\ol K\underset{\Z_2}{\times} S^\infty; \Q\Bigr) = H^j
(\ol K; \Q)^{\Z_2} .
\end{equation}
Now this last equality combined with \eqref{Leray} implies the Lemma.
\end{proof}

\begin{rema}

If $K'$ is another space with a $\Z_2$-action and $f : K' \to \ol K$
is an equivariant map, then the induced map $f_*$ in the
homologies respects the decomposition according to the
eigenvalues $+1$ and $-1$ of the $\Z_2$-action.
This implies that if

\begin{equation}
\label{primfiltr}
\ol K_0 \subset \ol K_1 \subset \dots
\subset \ol K_i \subset \dots \subset \ol K
\end{equation}

is a filtration by $\Z_2$-invariant subspaces, then the whole
spectral sequence decomposes according to the eigenvalues $+1$
and $-1$. (Indeed, apply the present Remark to $K' = \ol K_i.$)

In particular, if $K_i = \ol K_i / \Z_2$, then the cohomological
spectral sequence of the filtration
\begin{equation}
\label{filtr}
K_0 \subset K_1 \subset \dots \subset K_i \subset  \dots \subset
K
\end{equation}
can be identified with the $\Z_2$-invariant part of the
cohomological spectral sequence of~\eqref{primfiltr}.
In particular if the spectral sequence of \eqref{primfiltr} degenerates,
then so does that of \eqref{filtr}, too.
\end{rema}

In the next two subsections we shall apply these observations to
$\ol K = \Kaz_{prim}(k)$ and $K = \Kaz_{Morin}(k).$

\subsection{The spectral sequence for the Kazarian space $
\Kaz_{prim}(k)$
of prim maps}

\begin{lem}\
  \label{spectr-sequ-prim-kazar}
  \begin{enumerate}[a)]
  \item
    If the codimension $k$ of the maps is even, $k = 2\ell$,\ and
    $A = \Q [p_1, \dots , p_{\ell}],$ then the spectral sequence
    for prim maps is the following:
    \[
    \aligned
    &\begin{array}{c|c|c|c|c|c}
      i = & 0 & 1 & 2 & 3 \\
      \hline
      \bar E_1^{i,*} & A \oplus \chi A & \bar U_1 \cup (A \oplus \chi
      A) & \bar U_{2} \cup
      (A \oplus
      \chi A) & \bar U_{3} \cup (A \oplus \chi A)\\ 
      \hline
      \bar E^{i,*}_2 = \bar E_\infty^{i,*} & A & 0 & 0 & 0
      \end{array}\\
                \endaligned
    \]
    where $\chi$ is the Euler class of the universal oriented $k$-bundle 
$\gamma^{SO}_k$, and
$\bar U_{p}$ is the Thom class of the bundle $ \ol \xi_i = i \cdot
    (\gamma^{SO}_k \oplus \varepsilon ^1)$.
  \item
    If $k$ is odd, $k = 2 \ell + 1$, then
    $ \bar E^{i,*}_1 = \bar U_i \cup H^*(BSO(k))$
    and all the differentials vanish. 
    Hence $\bar E^{i,*}_\infty = \bar E^{i,*}_1.$
  \end{enumerate}
\end{lem}

\begin{rema}
The space $\Kaz_{prim} (k)$, which is the Kazarian space for all
oriented prim maps is  equal to the space $B\, SO(k + 1)$.
(To a prim map in codimension $k$ it can be associated an immersion of
codimension~$(k + 1)$, unique up to regular homotopy, its lift.)
So we know that the spectral sequence converges to
$H^*(B\, SO(k + 1); \Q)$.
\end{rema}

\begin{proof}[Proof for $k=2\ell$.]
The inclusion $B\, SO(k) \subset T i \gamma_k$
induces an injection 

$H^*(T i\gamma_k) \to H^*(B\, SO(k))$ with
image equal to the ideal generated by~$\chi^i.$

Hence we can rewrite the $\bar E_1$-term of the spectral sequence as follows:
\[
\begin{array}{c|c|c|c|c|c}
i = & \quad \ \, 0 \, \ \quad & 1 & 2 & 3 & \phantom{4} \\
\hline
E^{i,*}_1 & A \oplus \chi A & \chi A \oplus \chi^2 A &
\chi^2 A \oplus \chi^3 A & \chi^3 A \oplus \chi^4 A & \phantom{\chi^3 A \oplus \chi^4 A}
\end{array}
\]
This spectral sequence converges to $H^*(B\, SO(k + 1))$ which can
be identified with the subring $A$
in $ E^{0,*}_1 = A \oplus \chi A$. By dimensional reason we
have that the differential
\[
d_1 : E^{0,k}_1 \to E^{1,k}_1 \text{ maps the Euler class $\chi$ to $c\chi$ for
} c \neq 0,\ c \in \Q.
\]
(Indeed, all the further differentials $d_r$, $r > 1$ map the
group $E^{0,k}_r$ into a trivial group. Hence the only chance for the Euler
class to be killed is if its first differential is non-trivial.)
Similarly, $d_1 : E^{i,k}_1 \to E^{i + 1, k}_1$ maps $\chi^i$ to
a nonzero multiple of~$\chi^i$.
\end{proof}

\begin{lem}
\label{d1}
$d_1$ maps $\chi^{i + 1} \cdot A \subset E^{i,
*}_1$ onto $\chi^{i + 1} \cdot A \subset E^{i + 1, *}_1$ isomorphically.
\end{lem}

\begin{proof}
It is enough to show that $d_1 (\chi \cdot p_I) = d_1(\chi) \cup
p_I$ for any monomial $p_I$ of the Pontrjagin classes.
Let us denote by $\ol \Kaz_r$  the Kazarian space of the prim
$\Sigma^{1_r}$ maps (i.e.\ prim maps having only $\Sigma^0,
\Sigma^{1,0}, \dots, \Sigma^{1_r}$ type singularities), then
\[\ol
\Kaz_r = \bigcup\limits_{i \leq r} D\bigl(i(\gamma_k \oplus
\varepsilon^1)\bigr).
\]
(For $i = 0$ we put $D \bigl(i(\gamma_k \oplus \mathcal
E^1)\bigr) = B\, SO(k)$.)
The differential $d_1 : E^{i, *}_1 \to E^{i + 1, *}_1$ is by
definition the boundary operator $\delta$ in the exact sequence
of the triple $(\ol \Kaz_{i + 1},\ol \Kaz_i,\ol \Kaz_{i - 1})$.
The Pontrjagin classes $p_1, \dots, p_\ell$ in $E^{p, *}_1$ can
be identified with the Pontrjagin classes of the virtual normal
bundle $\nu^k$ over $\ol \Kaz_\infty$ restricted to the base space
of $D\bigl(i(\gamma_k \oplus \varepsilon^1)\bigr) \subset \ol \Kaz_\infty$.
Hence the inclusion $j : \ol \Kaz_i \hookrightarrow \ol \Kaz_{i + 1}$
maps them into each other: $j^*(p_\alpha) = p_\alpha$ for
$\alpha = 1, \dots, \ell$.
Now
\[
d_1(\chi \cup p_I) = \delta(\chi \cup p_I) = \delta \bigl(\chi\cup
j^* (p_I)\bigr) = \delta(\chi) \cup p_I + \chi \cup \delta\circ
j^*(p_I) = \delta(\chi) \cup p_I
\]

\end{proof}

\begin{proof}[Proof of Lemma~\ref{spectr-sequ-prim-kazar} for odd $k.$]
  Trivial computation shows that the rank of the $\bar E_1^{*,*}$
  is the same as that of $H^*BSO(k+1).$ Hence all the
  differentials vanish.
\end{proof}

\subsection{The spectral sequence for the Kazarian space $\Kaz_{Morin}(k)$
of arbitrary Morin maps}

Recall that $\xi_r$ is the universal normal bundle in the source of the stratum of
$\Sigma^{1_r}$ points, and $\ol \xi_r$ is that for prim maps. Let $U_r$ and
$\ol U_r$ denote their Thom classes, $\text {dim}\ U_r = \text {dim}\ \ol U_r = r(k+1).$ Further recall that
$A = \Q[p_1, \dots, p_{\ell}]$ is the subring of $H^*(BSO;\Q)$ generated by the
first $\ell$ Pontrjagin classes, and $\chi \in H^{2\ell}(BSO(2\ell);\Q)$
is the Euler class. 

\begin{lem}\
\label{E1}

  \begin{enumerate}[a)]
  \item
    For $k = 2 \ell$
    the group
    $H^*(T\xi_r; \Q)$
is the following:
    \[
    H^*(T\xi_r; \Q) = \begin{cases}
      U_{r} \cup A, &\text{if $r$ is odd},\\
      U_{r} \cup \chi \cdot A,\ &\text{if $r$ is even},
    \end{cases}
    \qquad r \geq 1.
    \]
    where $A = \Q[p_1, \dots, p_{\ell}].$
  \item
    For $k = 2\ell + 1$ this group is:
    $H^*(T\xi_r; \Q) =\begin{cases} 0 \ \text{if} \ r \ \text{is\ odd,} \\
    U_r \cup A \ \text{if}\ r \ \text{is\ even}\end{cases}.$
  \end{enumerate}
\end{lem}

\begin{proof}
This Lemma has been essentially proved in [R-Sz]. Here we recall the main steps.
The maximal compact subgroup  $G^O_r$ of the automorphism group
$\text{\rm Aut}\,(\eta_r)$, where $\eta_r : \bigl(\BR^{r(k + 1)}, 0\bigr) \to
\bigl(\BR^{(r + 1)k + r},0\bigr)$ has isolated $\Sigma^{1_r}$
singular point at the origin is $G^O_r = O(1) \times O(k)$ and its
representation $\lambda$ in the source is
\[
\left\lceil \frac{r - 1}{2}\right\rceil \cdot 1 + \left \lfloor
\frac{r + 1}{2} \right\rfloor \cdot \varrho_1 + \left\lceil
\frac{r}{2} \right\rceil (\varrho_1 \otimes \varrho_k) +
\left\lfloor \frac{r}{2} \right\rfloor \cdot \varrho_k
\]
($\varrho_i$ is the standard representation of $O(i)$ on $\BR^i$).

The ``oriented'' subgroup of $G_r^O$ (i.e. the subgroup that occurs as
structure group of the normal bundles of the stratum $\Sigma^{1_r}$ of maps
from an  oriented manifold into an oriented one) is
 $$G_r = \bigl\{(\varepsilon, A) \in O(1) \times
O(k)\mid \varepsilon^{r} \cdot \det \, A > 0 \bigr\}$$ and the
group $\ol G_r$ corresponding to ``oriented'' prim maps is $\ol
G_r = \bigl\{ (\varepsilon, A) \mid \varepsilon = 1$ and
$\det\, A > 0\bigr\}$.
The representations of $\ol G_r$ and $G_r$ are obtained
by restricting $\lambda$ to $SO(1) \oplus SO(k) = \ol G_r$
and to $G_r$, respectively.
Let us denote these representations by $\ol\lambda^{SO}$ and
$\lambda^{SO}$. From now on we omit the subindex $r$ in the notation of the groups.

The universal bundle associated to the representation
$\ol\lambda^{SO}$ of $\ol G$ is
\begin{equation}
\label{repr}
\left\lceil\frac{r -
1}{2}\right\rceil \cdot 1 \oplus \left\lfloor \frac{r + 1}{2}
\right\rfloor \gamma^{SO}_1 \oplus \left\lceil \frac{r}{2}
\right\rceil \bigl(\gamma^{SO}_1 \boxtimes \gamma^{SO}_k\bigr)
\oplus \left\lfloor \frac{r}{2} \right\rfloor \gamma^{SO}_k.
\end{equation}
(Note that here both $1$ and $\gamma^{SO}_1$ denote the trivial
line bundle over a point.
They are denoted differently because the above mentioned
$\Z_2$-action, (see Remark~\ref{involution}, and also the next paragraph) changes the orientation of $\gamma^{SO}_1$, while
keeps it on~$1$.)

The short exact sequence $1 \to \ol G \to G \to \Z_2
\to 1$ induces double covers $B\ol G \to BG$ of the
universal base spaces and the  bundles
$\ol \xi_r$ and $ \xi_r:$
\[
\ol \xi_r = \BR^{(r + 1)k} \underset{\ol \lambda^{SO}}{\times} E
\ol G \to \xi_r = \BR^{(r + 1)k}
\underset{\lambda^{SO}}{\times} EG.
\]
The action of the corresponding involution on $H^*(T\ol \xi_r; \Q)$ is 
what  we want to understand, because $\ol E^{r,*}_1 = H^*(T
\ol \xi_r; \Q)$ is the $E_1$-term of the spectral sequence of the Kazarian
space of oriented \underline{prim} maps
$\ol \Kaz_r$, and its $\Z_2$-invariant part is $H^*(T\xi_r;\Q)$ and that is
the  $E_1$ term  for the 
Kazarian space of oriented \underline{arbitrary Morin} maps.

\begin{enumerate}[\indent a]

\item
Let first $k$ be even, $k = 2 \ell.$ Then
\[
H^*(T\ol\xi_r ; \Q) =\ol U_r \cup H^*(B\,SO(k); \Q) = \ol U_r \cup
\Q[p_1, \dots, p_{\ell - 1}, \chi], \ \text{dim} \chi = k.
\]
The orientation of the fibre will change under this involution
by $(-1)^{\left\lfloor \frac{r + 1}{2} \right\rfloor}\cdot
(-1)^{\left\lceil \frac{r}{2} \right\rceil} = (-1)^{r + 1}$.
(The third summand in \eqref{repr} does not change its orientation
since $\varepsilon \cdot \det A > 0$, so $\gamma^{SO}_1$ and
$\gamma^{SO}_k$ change or keep their orientations at the same time.)
Hence the Thom class $\ol U_{r}$ of $\ol \xi_r$ is mapped into
$(-1)^{r + 1} \ol U_{r}$, and $\chi$ is mapped into $-\chi$.
Thus we proved part a). 

\item
The proof of part b) is similar.

\end{enumerate}
\end{proof}

\begin{proof} of Lemma~\ref{HK-Morin}.
The previous Lemma gives us the $E_1$ term of the  spectral sequence
for the Kazarian space of arbitrary Morin maps.
\begin{enumerate}[\indent a)]
\item
  Let $k = 2\ell$ and $A = \Q[p_1, \dots, p_{\ell}].$
  Then identifying $U_{r}$ with $\chi^r$ we obtain that the $E_1$-term
  looks as follows:
  \[
  \begin{array}{c|c|c|c|c|c|c}
    r  & \quad \ \, 0 \, \ \qquad & \quad \ \, 1 \, \ \quad  &  \
    \quad  2\quad \ \, & \, \ \quad 3 \, \ \quad  & \quad \ \, 4 \, \ \quad
    &\quad \ \, 5  \ \quad \phantom{6} \\
    \hline
    E^{r,*}_1 & A \oplus \chi A & \chi A & \chi^3 A & \chi^{3} A &
    \chi^5 A & \chi^5 A
  \end{array}
  \]
  The differentials map $\chi^r \cdot A$ onto $\chi^{r+1} \cdot A$ in the next
  column. This follows from Lemma~\ref{d1}.
  Hence $E^{0,*}_2 = E^{0,*}_\infty = A$,
  $E^{r,*}_2 = E^{r,*}_\infty = 0$ for $r > 0$.
We conclude that the composition map
$$ BO(k) \to BSO(k+1) = \Kaz_{prim}(k) \to \Kaz_{Morin}(k)$$
induces an isomorphism of the rational cohomology rings.
(Here the map $B\, O(k) \to B\, SO(k+1)$  induces to the direct sum of
the universal unoriented $k$-bundle with its determinant line bundle.)

Hence $\Kaz_{Morin}(k)$ and $BO(k)$ are rationally homotopically equivalent.
\item
  For $k= 2\ell+1$
  \[
  E_1^{r,j} = E_{\infty}^{r,j} =
  \begin{cases} 0 \ \text{if}\ r \ \text{is odd} \\
    p_{\ell+1}^i \cup H^{j} (B\, SO(k) ;\Q) \ \text {if}\ r = 2i
  \end{cases}
  \]
  where $p_{\ell+1}$ is the $\ell+1$-th Pontrjagin class of the virtual normal
  bundle.
Therefore $H^*(\Kaz_{Morin}(k);\Q) \approx H^*(BO(k+1);\Q).$
This isomorphism of cohomology groups can be induced by
the following map.
Consider the map  $$B\, O(k+1) \to B\, SO(k + 1) $$ corresponding to the oriented $k+1$-bundle
over $BO(k+1)$ equal to the tensor product of the universal (unoriented)
$k+1$-bundle with its determinant line bundle. Then compose it with the
composition map
 $$B\, SO(k + 1)  \cong \Kaz_{prim}(k) \to \Kaz_{Morin}(k).$$
 Hence
the spaces
 $\Kaz_{Morin}(k)$ and $BO(k+1)$ are rationally homotopically equivalent.

\end{enumerate}
\end{proof}
\begin{rema}
The natural forgetting map 
$$\Cob_{Prim}(n,k)\otimes \Q \to \Cob_{Morin}(n,k)
\otimes \Q$$ induced by considering an oriented  prim map just as an oriented Morin map, is an
epimorphism. Indeed, this is equivalent to the claim that
$H^*(\Kaz_{Morin}(k); \Q) \to H^*(\Kaz_{prim}(k); \Q)$
is injective, and this follows from Lemma~\ref{lem:2b}.
\end{rema}

\section{Thom polynomials of Morin singularities for Morin maps}

Having computed the rational Kazarian spectral sequence for all
Morin maps we have determined all the (higher) Thom polynomials
of all $\Sigma^{1_i}$- singularities for Morin maps.
First we recall the definitions.

\begin{defi}
Let $f : M^n \to P^{n+k}$ be a singular map and let $[\eta]$ be a
singularity class, let
$\eta(f)$ denote the set of $\eta$-points of~$f$.
Suppose that $\eta(f)$ is a manifold and its closure
$\ol\eta(f)$ represents a rational homology class dual to a cohomology
class $t$ of~$M^n$.

It has been shown by Thom that for generic $f$ the class $t$ is
a characteristic class of the virtual normal bundle of~$f$.
This characteristic class (expressed as a polynomial of the
basic characteristic classes, in our case of the Pontrjagin
classes) is called the \emph{Thom polynomial} of~$\eta$ and denoted by
$t = Tp_{\eta}(f).$

{\it {The higher Thom polynomials:}} Let $\eta$ be a highest singularity of
the $\tau$-map $f.$ Recall that the normal bundle of
$\eta(f)$ in $M^n$ can be induced by a homotopically unique map
$\eta(f) \overset{\nu_\eta}{\longrightarrow} BG_\eta$ from the
universal bundle over~$BG_\eta$.
Let $x \in H^*(BG_\eta ; \Q)$ be a cohomology class.
Let $i : \eta(f) \subset M^n$ be the inclusion.
Now we can ask how to express the push-forward class $i_!(\nu_\eta^*(x))$
as a polynomial of the characteristic classes of $M^n$ and $f^*(TP^{n+k})$.
This expression will be called {\it {the higher Thom polynomial
corresponding to the normal $G_\eta$-characteristic class $x$ of
the singularity~$\eta$}}.
Hence the set of all higher Thom polynomials associated to $\eta$ can be
identified with the image of the map $H^*(\Kaz_\tau,\Kaz_{\tau'}) \to
H^*(\Kaz_\tau)$, where  
$\tau' = \tau \setminus \{\eta\}.$
Now for a highest singularity $\eta'$ in $\tau'$ one can define the higher Thom
polynomials as above in $H^*(\Kaz_{\tau'}).$ But these are those {\it for
$\tau'$-maps.} If we want to define them for $\tau$-maps, then we have to
consider their preimages in $H^*(\Kaz_\tau)$, which are well defined modulo the
image of $H^*(\Kaz_\tau,\Kaz_{\tau'}).$
Hence
the higher Thom polynomials of the lower singularities are defined modulo
those
of the higher singularities.
Namely 
Kazarian has shown that
the $E_\infty$-terms of the cohomological spectral sequence of the
 Kazarian space $\Kaz_\tau$ can be identified
with the higher Thom polynomials of the singularities.
\end{defi}

Below we show for each element of the $E_\infty$ term of the spectral
sequence of the Kazarian space for Morin 
maps $\Kaz_{Morin}$ a singularity $\eta$ and its normal
$G_\eta$-characteristic class representing the given element.  
We start with the simpler case of prim maps.

\begin{enumerate}[A)]
\item {\bf {Prim maps.}} Let $f$ be the composition $f: M^n \overset{g}{\looparrowright} \BR^{n + k + 1}
\overset{\pi}{\longrightarrow} \BR^{n + k}$, where $g$ is an immersion.
Let $\chi = \chi_{k+1}$ be the
Euler class in 

$H^{k+1}(B\,SO(k+1);\Q).$

\begin{enumerate}[1)]
\item $k$ odd $= 2\ell + 1.$ Then $H^*(B\,SO(k); \Q) =
 \Q[p_1, \dots,p_\ell]$, recall that we 
denoted this ring by
$A.$ The Kazarian space $\Kaz_{prim}(k)$ for prim maps of codimension $k$, has the same
rational homotopy type as
$ B\,SO(k +~1). $ The Kazarian map
$\kappa_f:  M^n \to \Kaz_{prim}(k)$
corresponding to the prim map f can be identified with the normal map of
$g$, i.e 
$\kappa_f = \nu_g : M^n \to  B\,SO(k + 1).$
The cohomology ring of  $\Kaz_{prim}(k)$
is the polynomial ring over $A$ with
variable $\chi$, i.e. $H^*(B\,SO(k+1),Q) = A[\chi]$, where $\chi \in
H^{k+1}(B\,SO(k+1);Q)$ is the Euler class.
The Kazarian spectral sequence degenerates, i.e.for any $r$ the member $E^{*,*}_r$ looks as
follows:
\[
\begin{array}{c|c|c|c|c|c}
i = & 0 & 1 & 2 & 3 & \phantom{4}\rule{0pt}{10pt} \\
\hline
E^{i,*}_r & A & \chi A & \chi^2 A &
\chi^3 A &
\phantom{UUUU}\rule{0pt}{10pt} \\
\hline
\eta & \Sigma^0 & \Sigma^{1,0} & \Sigma^{1,1} & \Sigma^{1,1,1} &
\rule{0pt}{10pt}
\end{array}
\]

The last row in this table shows the singularity stratum, and the middle raw
shows which part of the cohomology ring $H^*(\Kaz_{prim},Q)$ is represented by the Thom polynomial
and the higher Thom polynomials of this singularity stratum $\eta$ (modulo
those of the higher singularities).
The minimal dimensional parts of each column give the Thom polynomial of 
the stratum.
Hence $\eta = \Sigma^{1_i}(f)$ represents $\chi^i$ in $H^*(M^n).$
\footnote{ A priori we obtain only that the Thom polynomial of
$\Sigma^{1_i}$ is $c \cdot \chi^i$ , where $c$ is a non-zero rational number.
But the similar spectral sequence can be considered with any coefficient ring,
containing  $\frac12$. Then we obtain, that for integer coefficients
the Thom polynomial of $\Sigma^{1_i}$ is equal to $\chi^i$, modulo $2$-primary
torsion elements.}

Let $x = p_I = p_{i_1} \dots p_{i_S}$, $i_1 \leq \dots \leq i_s
\leq \ell$ be a normal characteristic class of $\Sigma^{1_i}(f)$
\footnote{Note that these Pontrjagin classes are not those of the normal
bundle of the stratum. For each singularity $\eta$ here $G_\eta = SO(k)$,
and $p_i$ is the Pontrjagin class in $H^*(BSO(k)).$ But the normal bundle of
the stratum $\Sigma^{1_i}$ is $i(\gamma_k^{SO} \oplus \varepsilon^1).$} 
 and let
$j_i: \Sigma^{1_i}(f) \subset M^n$ be the embedding.
Then $j_{i!}(x) = \chi^i \cdot p_I$.

\bigskip

\item $k$ even $=2\ell$.

Let us denote again by $A$ the ring $\Q[p_1, \dots,p_{\ell}]$ and by $\chi = \chi_k$
the Euler class $\chi \in H^k(B\,SO(k),Q).$

Then $H^*(B\, SO(k)) = A \oplus \chi \cdot A.$

We have seen in
Lemma \ref{spectr-sequ-prim-kazar} that
the $E^{*,*}_\infty$ member of the spectral sequence is
\[
\begin{array}{c|c|c|c|c|c}
i & 0 & 1 & 2 & 3 & \phantom{4}\rule{0pt}{10pt} \\
\hline
E^{i,*}_\infty & A  & 0 & 0 & 0 &
\phantom{UUUU}\rule{0pt}{10pt} \\
\hline
\eta & \Sigma^0 & \Sigma^{1,0} & \Sigma^{1,1} & \Sigma^{1,1,1} &
\rule{0pt}{10pt}
\end{array}
\]

The last row shows the singularities corresponding to the blocks (by glueing)
of the Kazarian space and hence they correspond to the  columns of the
Kazarian spectral sequence.

 and  
 the map 
$H^*(BSO(k+1);\Q) \to H^*(BSO(k);\Q) = A$
is injective.)

{\bf {Conclusion.}}
All the singularity strata and all their normal characteristic cycles 
are null-homologous for a prim map of even codimension (except those of  the
regular stratum $\Sigma^0.$)

\end{enumerate}

\item {\bf {Arbitrary Morin maps.}}
\begin{enumerate}[1)]
\item
Case $k$ even $=2\ell.$
The above Conclusion remains true, i.e.\ all (higher) Thom
polynomials of the strata $\Sigma^{1_i}$, for $i \ge 1$ vanish. This follows from Lemma~\ref{lem:2b}.
\item
Case $k$ odd $= 2\ell + 1.$
The $E^{*,*}_r$ members of the spectral sequence are for any $r$ the
following. (Recall that $U_i$  denotes the Thom class of the normal bundle $\xi_i$
of the stratum $\Sigma^{1_i}.$, it has dimension $i(k+1).$)

\[
\begin{array}{c|c|c|c|c|c|c}
i & 0 & 1 & 2 & 3 & 4 & \phantom{5}\rule{0pt}{10pt} \\
\hline
E^{i,*}_r & A & 0 & U_{2} \cup A & 0 & U_{4} \cup A &
\phantom{UUUU}\rule{0pt}{10pt} \\
\hline
\eta & \Sigma^0 & \Sigma^{1,0} & \Sigma^{1,1} & \Sigma^{1,1,1} &
\Sigma^{1_4} &
\rule{0pt}{10pt}
\end{array}
\]

 (This is so for the $r=1$ by Lemma~\ref{E1}, and all the differentials vanish by
dimensional reason, since $A$ has no nonzero element of odd degree.)

The higher Pontrjagin classes vanish: $p_j(\nu_f) = 0$ for $j > \ell+1$, see
 Lemma~\ref{Pontrj}  below.
The cohomology class in $H^*(M^n;\Q)$ Poincare\`e dual
 to the homology class $[\ol
\Sigma^{1_i}]$ represented by the closure of the
singularity stratum $\Sigma^{1_i}$ is
\[
D_M\bigl[\ol\Sigma^{1_i}(f) \bigr] = \begin{cases}
0, &\text{if $i$ is odd,}\\
p^{i/2}_{\ell +1}, &\text{if $i$ is even.}
\end{cases}
\]
(Note that $H^* \bigl(BG_{\Sigma^{1_{2i}}}\bigr) = H^*(B\, SO(2\ell +
1)\bigr) = \Q[p_1, \dots, p_{\ell}]$).

If $x \in H^*\bigl(BG_{\Sigma^{1_{2i}}}; \Q\bigr)$ is a monomial on the
Pontrjagin classes $x = p_I =
p_{i_1} \dots p_{i_s}$, $i_1 \leq \dots \leq i_s \leq \ell$,
$(j_{2i})_!(p_I) = p^{i}_{\ell+1} \cdot p_I \in H^*(M^n)$.
Hence any characteristic class of the normal bundle  $\nu_f$ of the map $f$
is represented by a higher Thom polynomial.
\end{enumerate}
\end{enumerate}
\begin{lem}\label{Pontrj}
If $f : M^n \to P^{n + k}$ is a Morin map, the manifolds  $M^n$ and $ P^{n+k}$ are
oriented, the codimension $k$ is odd, $k =
2\ell + 1$ and $\nu_f : M^n \to B\, SO$ is the virtual normal
bundle, then $p_j(\nu_f) = 0$ for $j > \ell + 1$.
\end{lem}

\begin{proof}
Since $f$ is a Morin map, the stable normal map $\nu_f: M^n \to BSO$ can be decomposed as $\nu_f
= \pi \circ \kappa _f$, where $\kappa _f : M^n \to \Kaz_{Morin}(k)$
is the  map into the Kazarian space.

Now it remained to show that $\pi^*(p_j) = 0$ for $j > \ell+1$.

This follows from the diagram
\[
\aligned
&\phantom{AAA}\\
H^*(B\, SO(2\ell+2)) &\\[-7mm]
&\begin{CD}
 = H^* (\Kaz_{prim}) @< p^* <<
\hspace*{-3mm}<\hspace*{3mm} H^*(\Kaz_{Morin}(k))\\
@AA \ol \pi^* A @AA \pi^*A\\
H^*(B\,SO) @= H^*(B\, SO)
\end{CD}\hspace*{30mm}
\endaligned
\] 
Here $p : \Kaz_{prim}(k) \to \Kaz_{Morin}(k)$ is the quotient map.
Hence $p^*$ is injective.
Now we see that $\ol \pi^*(p_j) = 0$ for $j > \ell+1$, hence
$\pi^*(p_j) = 0$.
\end{proof}

\section {$\Sigma ^{1_r}$-maps.}

Let us denote by $\Cob_{\Sigma^{1_r}}(n,k)$ the oriented cobordism group
of Morin maps having at most $\Sigma^{1_r}$ singularities.
Recall that $ A = \Q[p_1, \dots , p_{\ell}].$

\begin{thm}
The group
$\Cob_{\Sigma^{1_r}}(n,k)\otimes \Q$ is isomorphic to the degree $n$ part of the
following graded ring:
$
\begin{cases} A\ \text {if}\ k = 2\ell\  \text{and}\ r \ \text{is \ odd}\\
 A \oplus \chi^{r+1}\cdot A\ \text{if}\ k = 2\ell\  \text{and}\ r \ \text{is \
   even}\\
A/ p_{\ell}^{[\frac{r+1}2]} = 0\ \text {if}\ k = 2\ell -1
\end{cases}$
\end{thm}

\begin{proof}
We have to consider the above Kazarian spectral sequence restricted to the
first $r+1$ columns.

\end{proof}
\section {Elimination of singularities}

Arnold and his coauthors in \cite{Ar1}, at page 212  put the question whether vanishing of
the Thom polynomial of a singularity $\eta$ for a map $f$ is enough for having a map $g$
homotopic to $f$ and having no $\eta$ points.

The answer to this question is negative, see \cite{Sz6}.
Here we consider the analogous question of elimination by $\tau$-cobordism
and give a complete solution to it.

Let $\tau$ be as before a set of stable singularities. Let $\eta$ be a top
singularity in $\tau.$ Let $f: M^n \to P^{n+k}$ be a $\tau$-map.
Then clearly the restriction of $f$ to the $\eta$-stratum $\eta(f)$
is an immersion and its normal bundle is induced from the bundle
 $\wt \xi_\eta.$ 
We show that the cobordism class of this immersion 
$f \big| _\eta(f)$ with normal bundle induced from $\wt \xi_i$ is a complete
(necessary and sufficient) obstruction to the elimination of $\eta$-points by
a $\tau$-cobordism.

Translating this obstruction into cohomologies we obtain 
the following:
Suppose that the Gysin map $f_!$ annihilates not only the Thom polynomial of $\eta$ for
the map $f$, but also all its higher Thom polynomials of $\eta$ (i.e. the
images of the characteristic classes of the normal bundle of $\eta(f) \subset
M^n$ under the Gysin map induced by this inclusion. 
 Then  a non-zero multiple of the $\tau$-cobordism class of $f$ contains an
$\eta$-free map.

\begin{nota}

Given a vector bundle $\zeta$ over a given base space and a smooth manifold
$P^{n+k}$, we will denote by
$\text{Imm}^{\zeta}(P^{n+k})$
the cobordism group of immersions in $P^{n+k}$ with normal
bundle induced from $\zeta .$
\end{nota}

\begin{thm}
\label{elim}
Let $f : M^n \to P^{n+k}$ be a $\tau$-map, $\tau = \{\eta_0 < \eta_1 <
\dots < \eta_i\}$, $\tau' = \{\eta_0 < \dots < \eta_{i - 1}\}$.
Recall that $f\big|_{\eta_i(f)} : \eta_i(f) \to P^{n+k}$ is an
immersion with normal bundle induced from the universal bundle
$\wt\xi_i = E G_i \underset{\wt\lambda_i}{\times} \BR^{c_i + k}$.
Hence it represents an element
$\theta(\eta, f) \in \text{\rm Imm}^{\wt\xi_i}(P^{n+k})$ in the
cobordism group of $\wt\xi_i$-immersions.
We claim that this element $\theta(\eta, f)$ vanishes if and only if 
there exists a $\tau'$-map $g : M^n \to P^{n+k}$
$\tau$-cobordant to $f.$
\end{thm}

\begin{proof}
Consider the cofibration
\begin{equation}
\label{virt}
V_{i - 1} \subset V_i \to V_i / V_{i -
1} \to T\wt \xi_\eta
\end{equation}
 of virtual complexes, 
 and apply the functor $\Gamma = \Omega^{\infty}
S^{\infty}.$
(Let us ignore for a while the fact that $V_i$ and $V_{i-1}$ are
virtual complexes, and let us treat them at first if they were usual complexes.)

 By a theorem of \cite{B--E} the functor $\Gamma$ turns a
cofibration into a fibration. Hence
we obtain a fibration

\begin{equation}
\label{key}
X_{\tau'} \subset X_\tau \to \Gamma
T\wt\xi_i
\end{equation}

Applying the homotopy functor $[P^{n+k}, \ - \ ]$ to this fibration, we
obtain that the sequence
\[
\Cob_{\tau'}(P^{n+k}) {\longrightarrow} \Cob_\tau(P^{n+k})
\overset{\beta}{\longrightarrow} \text{\rm Imm}^{\wt\xi_i} (P^{n+k})
\]
is exact.
$\beta$ maps $[f]$ to $\theta(\eta, f)$.

Now let us address the problem that  $V_i$ and $V_{i-1}$ are virtual
complexes.
For any natural number $m$ there exist:

a finite dimensional approximation  $V_{i,fin}$ and a
natural number $N$ such that the suspension $S^NV_{i,fin}$ exists and it
represents the $N+m$-homotopy type of $S^NV_i.$

Similarly we can suppose that for the same $m$ and $N$ the same holds for 
$i-1$ substituted for $i,$
and  
$S^N V_{i,fin} = S^NV_{i-1,fin} \bigcup D^N\times T\wt\xi_{\eta,fin}$ for a finite
approximation $T\wt\xi_{\eta,fin}$ of the Thom space $T\wt\xi_\eta$
($m$-equivalent to it).
Now instead of the virtual cofibration above we consider the (usual)
cofibration
$$S^NV_{i-1,fin} \subset S^NV_{i,fin} \to S^NT\wt\xi_{\eta,fin}$$
and apply the functor $\Gamma$ to it.
We obtain a fibration

$$\Gamma S^NV_{i-1,fin} \to \Gamma S^NV_{i,fin} \to \Gamma
S^NT\wt\xi_{\eta,fin}$$   

(Here the first space is the fiber the second is the total space and the third
one is the base space.)

Let us recall that whenever we have such a sequence of maps
$F \to E \to B$ we can continue it to the left infinitely and obtain an infinite to the
left sequence of fibrations (called the resolvent of the map $E\to B$).

$$\dots \Omega^2 F \to \Omega^2E \to  \Omega^2B \to \Omega F \to \Omega E \to \Omega B
\to F \to E \to B$$ 

In particular we can consider the fibration of $N$-th loop spaces

$$\Omega ^N\Gamma S^NV_{i-1,fin} \to \Omega ^N \Gamma S^NV_{i,fin} \to 
\Omega ^N\Gamma
S^NT\wt\xi_{\eta,fin}$$

Noticing the $\Omega ^N\Gamma S^N = \Gamma$, we obtain a usual fibration
$m$-equivalent to the fibration obtained from the cofibration of virtual
complexes.

\end{proof}

\begin{defi}
The fibre bundle \eqref{key}
will be called the
{\it {key fibration}} for the pair $(\tau, \tau')$.
\end{defi}

\begin{rema}
The existence of such a fibre bundle was conjectured by Endre Szab\'o.
This observation and the Kazarian conjecture are the main points of the paper.
\end{rema}

\begin{cor}
$[f] \otimes 1_\Q \in \Cob_\tau(P^{n+k})\otimes \Q$ contains a $\tau'$-map~$g$
(i.e.\ $\eta(g) = \emptyset$) iff for any $x \in H^*(BG_{\eta}; \Q)$
the class $(f\mid _{\eta(f)})_!(j^*(x)) = 0$ in $H^*(P^{n+k}; \Q)$, where $j : \eta(f)
\to BG_{\eta}$ is the map inducing the normal bundle of $\eta(f)$ in~$M^n$.
\end{cor}

\begin{proof}
Let $\theta_f \in \{P^{n+k}, T \wt \xi_{\eta}\}$ be the stable homotopy
class corresponding to cobordism class $\theta(\eta, f)$ of
 the immersion $f\big|_{\eta(f)} :
\eta(f) \looparrowright P^{n+k}$.
It is standard to show (see \cite{Sz6}) that $\theta^*_f(U \cup
x) = (f\mid \eta(f))_!(j^*(x))$, where $U \in H^*(T \wt\xi_\eta; \Q)$ is the Thom class.

The map
\[
\aligned
\text{\rm Imm}^{\wt\xi_\eta}(P^{n+k}) \otimes \Q = \{P^{n+k}, T\wt\xi_\eta\}
\otimes \Q &\longrightarrow \Hom\bigl(H^*(T\wt\xi_\eta; \Q), H^*(P^{n+k}; \Q)\bigr)\\
\alpha \quad &\longmapsto \quad \alpha^*
\endaligned
\]
is an isomorphism (see for example \cite{Sz6}, page 321).
Hence the class $\bigl[f\big|_{\eta(f)}\bigr] \otimes 1_\Q \in \text{\rm
Imm}^{\wt\xi_\eta}(P^{n+k}) \otimes \Q$ is equal to $\theta_f \otimes
1_\Q \in \{P^{n+k}, T\wt \xi_\eta\}$ and that corresponds to the map
$\theta_f^*.$
\end{proof}
\begin{rema}
This Corollary shows that while the answer was negative to the question of Arnold about the
elimination of the highest singularity by homotopy - vanishing of the
corresponding Thom polynomial is not enough for homotopical elimination - the reformulated
question about elimination by $\tau$-cobordism of the top singularity $\eta$ of a
$\tau$-map has an answer quite close to the one formulated in \cite{Ar1}.
Namely after taking the tensor product with the rational numbers $\Q$
of the cobordism group of $\tau$-maps into a fixed manifold $P^{n+k}$ we
obtain, that
an element $[f: M^n \to P^{n+k}]$ of this group contains an $\eta$-free map
iff all the characteristic numbers of $\tau$-maps ``coming from $\eta$'' vanish 
for $f,$ i.e. all the $(x,z)$-characteristic numbers 
$\bigl< \kappa^*_f(x) \cup f^*(z), [M^n]\bigr>$ vanish if $x \in \text{im}
(H^*(\Kaz_\tau, \Kaz_\tau',\Q) \to H^*(\Kaz_\tau;\Q))$, see Definition~\ref{char}.
\end{rema} 

\part{General theorems on the key fibration
 $X_\tau
\stackrel{X_{\tau'}}{\longrightarrow} \Gamma T\wt\xi_\eta.$}

Here we shall study this bundle from rational point of view and
give a general structure theorem for it.  There are
quite general conditions implying that this bundle is
a direct product rationally.
In this case 
for any oriented manifold~$P^{n+k}$
\[
\Cob_\tau(P^{n+k}) \otimes \Q \approx (\Cob_{\tau'}(P^{n+k}) \otimes \Q) \oplus
(\text{\rm Imm}^{\wt\xi_\eta}(P^{n+k}) \otimes \Q).
\]

\begin{defi}
We say that a fibration $\pi: E \overset{F}{\longrightarrow} B$
is an $H$-fibration, if all the spaces $E, F, B$ are $H$-spaces
of finite types
{\footnote {i.e. their homology groups are finitely generated in each 
dimension.}}
and the maps $F \subset E \to B$ are $H$-homomorphisms.
\end{defi}

\begin{defi}
We say that an $H$-map is $F$-trivial ($E$-trivial,
$B$-trivial respectively) if the space $F$ (respectively $E$ or $B$) is
contractible.
\end{defi}

\begin{lem}
[$\Q$-classification of $H$-bundles]
Any $H$-fibration can be represented rationally as a product of
three $H$-fibrations, which are $F$-, $E$- and $B$-trivial, respectively.
\end{lem}

The proof of this Lemma is a standard application of spectral sequences
therefore we omit it (see \cite{Sz10}).

\begin{thm}
\label{euler}
Let $\eta$ be a top singularity class in $\tau$, $\tau' = \tau
\setminus \{\eta\}$.
Suppose that the universal bundles $\xi_\eta$ and $\wt\xi_\eta$
are oriented and at least one of them has non-zero rational
Euler class: $e(\xi_\eta) \neq 0$, or $e\bigl(\wt\xi_\eta\bigr)
\neq 0$ in $H^*(BG_\eta; \Q)$.
Then the key fibre bundle $X_\tau
\overset{X_{\tau'}}{\longrightarrow} \Gamma T \wt\xi_\eta$
splits rationally, i.e.\ $X_\tau {\cong}_\Q
X_{\tau'}\times \Gamma T\wt\xi_\eta$.
\end{thm}

\begin{proof}
Let us consider the diagram (the homology groups are meant with
rational coefficients):
\[
\begin{array}{@{\hspace*{2pt}}c@{\hspace*{2pt}}c@{\hspace*{2pt}}c@{\hspace*{2
pt}}c@{\hspace*{2pt}}c@{\hspace*{2pt}}c@{\hspace*{2pt}}c@{\hspace*{2
pt}}}
H_*(\Kaz_{\tau'}) & \longrightarrow & H_*(\Kaz_\tau) & \longrightarrow
& H_*(T\xi_\eta) & \overset{\partial_K}{\longrightarrow} &
H_*(S\Kaz_{\tau'})  \longrightarrow \\[1.5mm]
\big\downarrow\approx &  & \big\downarrow\approx &  &
\big\downarrow\approx &  & \big\downarrow\approx   \\[1.5mm]
H_{* + k}(V_{\tau'}) & \longrightarrow & H_{* + k} (V_\tau) &
\longrightarrow & H_{* + k}(T\wt\xi_\eta) &
\overset{\partial_V}{\longrightarrow} & H_{* + k}(SV_{\tau'}) \\[1.5mm]
\big\downarrow Hu & & \big\downarrow & &
\big\downarrow && \big\downarrow \\[1.5mm]
\pi^s_{* + k}(V_{\tau'})\otimes \Q & \longrightarrow &
\pi^s_{* + k}(V_{\tau})\otimes \Q & \longrightarrow &
\pi^s_{* + k}(T\wt\xi_\eta )\otimes \Q & \longrightarrow &
\pi^s_{* + k}(SV_{\tau})\otimes \Q \\[1.5mm]
\big\downarrow\approx &  & \big\downarrow\approx &  &
\big\downarrow\approx &  & \big\downarrow\approx   \\[1.5mm]
\pi_{* + k}(X_{\tau'}) \otimes \Q & \longrightarrow &
\pi_{* + k}(X_{\tau}) \otimes \Q & \longrightarrow &
\pi_{* + k}(\Gamma T\wt\xi_\eta) \otimes \Q &
\overset{\partial}{\longrightarrow} &
\pi_{* + k - 1} (X_{\tau'}) \otimes \Q
\end{array}
\]
The two upper rows are essentially the homology exact sequences
of the pairs $(\Kaz_\tau,\Kaz_{\tau'})$ and $(V_\tau, V_{\tau'})$, respectively.
We are writing $H_*(S\Kaz_{\tau'})$ instead of $H_{* -
1}(\Kaz_{\tau'})$, because we consider this row as the dual of the
cohomology sequence, which is obtained from the Puppe  sequence
(see Chapter 9 in [Hu])
\[
\Kaz_{\tau'} \subset \Kaz_{\tau} \to T\xi_\eta \hookrightarrow S\Kaz_{\tau'}
\]
applying the functor $[-, K(\Q, *)]$.
Then $\partial_K$ is the dual of the map induced by the
inclusion $j : T \xi_\eta \hookrightarrow S\Kaz_{\tau'}$.

But the map $j^* : H^*(S\Kaz_{\tau'}) \to H^*(T\wt\xi_\eta)$ is
trivial if $e(\xi_\eta) \neq 0$.
Indeed, the multiplication in $H^*(S\Kaz_{\tau'})$ is trivial, while in
$H^*(T\xi_\eta)$ any element has the form $U_\eta \cup x$ for
some $x \in H^*(BG_\eta)$ and $(U_\eta \cup x_1) \cup (U_\eta
\cup x_2) = U_\eta \cup e(\xi_\eta)\cup x_1\cup x_2\neq 0.$ ($H^*(BG_{\eta};\Q)$
has no zero divisor since it is a subring of the ring of polynomials
$H^*(BT;\Q,)$ where $T$ is the maximal torus in $G_\eta.$)

If $e(\wt\xi_\eta) \neq 0$, then one can use similarly the
second row.
(The vertical arrows from the first row are the Thom
isomorphisms corresponding to the virtual bundle~$\nu^k$, those from
the second row are the stable Hurewicz homomorphisms (being
isomorphism modulo torsion).)
In any of the two cases ($e(\xi_\eta) \neq 0$ or $e(\wt\xi_\eta) \neq 0$)
the map $\partial : \pi_{* + k}(\Gamma T \wt\xi_\eta) \otimes
\Q \to \pi_{* + k - 1} (X_{\tau'}) \otimes \Q$ vanishes (see the bottom row).
Then $\pi_*(X_\tau) \otimes \Q = \pi_*(X_{\tau'})\otimes \Q
\oplus \pi_*(\Gamma T \wt\xi_\eta) \otimes \Q$.
The spaces $X_\tau$, $X_{\tau'}$, $\Gamma T\wt\xi_\eta$ -- being
$H$-spaces of finite types -- are rationally homotopically
equivalent to products of some rational Eilenberg--MacLane
spaces $K(\Q,i)$.
It follows that $X_\tau {\cong}_\Q X_{\tau'} \times
\Gamma T \wt\xi_\eta$.
\end{proof}

Now we give a general rational decomposition theorem  of the key bundle
 $X_\tau
\stackrel{X_{\tau'}}{\longrightarrow} \Gamma T\wt\xi_\eta.$
Let $\tau$, $\tau'$, $\eta$ be as above.
Let us denote the rational homotopy types of the spaces
$X_\tau$, $X_{\tau'}$, $\Gamma T\wt\xi_\eta$ by $X$, $X'$ and
$\Gamma$, respectively.
Recall that all these spaces are rational $H$-spaces, hence they
decompose into products of rational Eilenberg--MacLane spaces
$K(\Q, i)$ in a unique way.

\begin{thm}
  \label{H-space-decomposition}
  \
  \begin{enumerate}[a)]
  \item
    There is a rational $H$-space $B$
(i.e.. $B$ is a product of spaces $K(Q,i)$ for some $i$-s), 
such that $X \cong_\Q
    \frac{X' \times \Omega B}{\Omega \Gamma} \times B$.
    (Note that dividing by the space $\Omega\Gamma$ means that its
    Eilenberg--MacLane factors $K(\Q,j)$ all occur as factors of the
    numerator and we cancel these factors.)
  \item
    The rational space $B$ is defined in a unique way by the
    following formula:
    \[
    \pi_m(B) = E_{r, m - r - k}^\infty \otimes \Q,
    \]
    where $E_{r,*}^\infty$ is the last column (corresponding to the top
    singularity~$\eta =\eta_r$) in the homological spectral sequence of the
    Kazarian space~$\Kaz_\tau$.
  \end{enumerate}
\end{thm}

\begin{proof}
By the $\Q$-classification theorem of $H$-bundles we have:
There exist $H$-bundles: $*\overset{\Omega A}{\longrightarrow}
A$; $B \overset{*}{\longrightarrow} B$;
$C\overset{C}{\longrightarrow} *$ whose product is the bundle
$X{\longrightarrow} \Gamma$ with fibre $X'.$
In particular, $X = B \times C$; $X' = \Omega A \times C$
and $\Gamma = A \times B$
and so  $\Omega \Gamma = \Omega A \times \Omega B$.)
Part~a) is trivial now.

 To prove b) note that
$\pi_*(X) = \pi_*(B) \oplus \pi_*(C)$ and
$\pi_*(\Gamma) = \pi_*(B) \oplus \pi_*(A)$ and $p_* : \pi_*(X) \to \pi_*(\Gamma)$ is
isomorphism on $\pi_*(B)$, and this is the maximal such
subgroup both in $\pi_*(X)$ and in $\pi_*(\Gamma)$.
There is a commutative diagram
\[
\begin{CD}
\pi_m(X') @>>> \pi_m(X) @>p_*>> \pi_m(\Gamma)\\
@VVV @VVV @VVV\\
H_{m - k}(\Kaz_{\tau'}; \Q) @>>> H_{m - k}(\Kaz_\tau; \Q) @>\alpha >>
H_{m - k}(T\xi_\eta; \Q)
\end{CD}
\]
where the vertical arrows are the compositions of a stable Hurewicz
homomorphism and a Thom isomorphism. Hence they are isomorphisms.
For example the middle vertical arrow is the following:
\[
\gathered
\textstyle\pi_*(X)
= \textstyle \pi^s_*(V_\tau) \otimes \Q
\overset{h^{V_\tau}}{\longrightarrow} H_*(V_\tau; \Q) =
H_*(T\nu^k; \Q) \begin{array}{c}
\text{\small Thom}\\[-1mm]
\hbox to10mm{\rightarrowfill}\\[-1mm]
\approx
\end{array} H_{* - k}(\Kaz_\tau; \Q).\big)\\
\,H_{m - k}(\Kaz_{\tau}; \Q) = \bigoplus^r_{i = 0} E^\infty_{i,m - k - i} \\
E^\infty_{r,m - k - r} = \text{\rm im}\ \alpha :H_{m - k}(\Kaz_\tau; \Q) \to H_{m
- k} (\underbrace{T\xi_\eta}_{\Kaz_\tau / \Kaz_{\tau'}} ; \Q)\bigr .
\endaligned
\]
Now we have

$\pi_m(B) = \text{\rm {im}} (p_*) =
\text {im}\ \alpha =  E^\infty_{r,m - k - r}. $
\end{proof}
\begin{exam}
Under the conditions of Theorem~\ref{euler} the equality $B = \Gamma$ holds and so
$X = X' \times \Gamma$.
\end{exam}

\section{A Postnikov-like tower for the classifying space $X_\tau$ of $\tau$-maps}

\begin{thm}
  \label{postnikov-like-tower}
  \
  \begin{enumerate}[a)]
  \item
    The space $X_\tau$ classifying cobordisms of $\tau$-maps for

    $\tau = \{\eta_0 < \eta_1 < \dots < \eta_r\}$ can be obtained by
    a sequence of fibrations:
    \begin{equation}
      \label{post}
    \Gamma_r \overset{\Gamma_{r - 1}}{\hbox to9mm{\leftarrowfill}}
    A_2 \overset{\Gamma_{r - 2}}{\hbox to9mm{\leftarrowfill}} A_3
    \dots \overset{\Gamma_1}{\longleftarrow} X_\tau
    \end{equation}
    where $\Gamma_i = \Gamma T\wt\xi_{\eta_i}$.
  \item
    For any manifold $P^{n+k}$ there is a spectral sequence
    with $E_1$-term
    $$E^{i,j}_1 = \text{\rm Imm}^{\zeta_i}(P^{n+k} \times~\BR^j),$$
    where $\zeta_i = \wt\xi_{\eta_i}$ is the universal normal
    bundle of the $\eta_i$-stratum in the target, and converging to
    $\Cob_\tau(P^{n+k} \times \BR^*)$
    \[
    E^{i,j}_1 = \text{\rm Imm}^{\zeta_i} (P^{n+k} \times \BR^j)
    \Longrightarrow \Cob_\tau(P^{n+k} \times \BR^{i + j}).
    \]
Recall that given a vector bundle $\zeta$, and a manifold $Q$
$ \text{\rm Imm}^{\zeta_i}(Q)$ denotes the cobordism group of immersions of
closed oriented manifolds into $Q$ with normal bundle induced from $\zeta.$
  \end{enumerate}
\end{thm}

\begin{proof} \

\begin{enumerate}[a)]
  \item
Let $\tau_j$ be $\tau_j = \{\eta_0 < \eta_1 < \dots < \eta_j\}$
and denote by $V_j$ the virtual complex $V_{\tau_j}$, and by
$\Gamma_j$ the space $\Gamma T\wt\xi_j = \Gamma(V_j / V_{j - 1})$.
Now put $A_s \overset{\text{\rm def}}{=\!\!=} \Gamma(V_r / V_{r
- s})$ and recall that $\Gamma$ turns a cofibration of spaces
(or virtual spaces) into fibration.
 \item
In the same way as the Postnikov tower of a space $X$ induces a spectral sequence 
starting with $[S^jP, K(\pi,i)] = H^{i-j}(P^{n+k},\pi)$ and 
converging to
$[P^{n+k},X]$ for any $P^{n+k}$, (see \cite{M--T}), a Postnikov-like tower induces the spectral sequence in the Theorem.
\end{enumerate}
\end{proof}

Recall that $c_j$ denotes the codimension of the singularity
$\eta_j$ in the source.
Let us denote by $s_1, s_2, \dots, s_{\ell - 1}$ those indices,
where the parity of the number $c_j$ changes, and $s_\ell = r$,
i.e. \[
c_1 \equiv c_2 \equiv \dots \equiv c_{s_1} \not\equiv c_{s_1 +
1} \equiv c_{s_1 + 2} \equiv \dots \equiv c_{s_2} \not\equiv
c_{s_2 + 1} \equiv \dots  \mod 2.
\]
We shall say that the indices $i$ and $i'$, $1 \leq i < i' \leq
r$ belong to the same block if
$
 {\exists}\ t : s_t < i < i' \leq s_{t + 1}.
$

Let us consider the (homological) Kazarian spectral sequence
$E^1_{i,j}$ (with rational coefficients) of the space $\Kaz_\tau$
(denoted also by~$\Kaz_r$) induced by the filtration $\Kaz_0 \subset
\Kaz_1 \subset \dots \subset \Kaz_r$, where $\Kaz_\beta =
\bigcup\limits_{\alpha \leq \beta} D\xi_\alpha$ (i.e.\ $\Kaz_\beta$
is the Kazarian space for the set of singularity classes
$\{\eta_0 < \dots < \eta_{\beta} \}$ for $\beta \leq r$.

\begin{prop}\
  \begin{enumerate}[a)]
  \item
    A differential $d_t : E^t_{i,j} \to E^t_{i',j'}$
    vanishes in this spectral sequence, if the indices $i$, $i'$ belong to
    the same block.
  \item
    The quotient spaces $\Kaz_{s_{t + 1}}\bigm/\Kaz_{s_t}$ have
    the same rational homology groups as the wedge product
    $T\xi_{s_t + 1} \vee T\xi_{s_t + 2} \vee \dots \vee T\xi_{s_{t + 1}}$.
  \item
    The space $X_\tau$ admits a simplified Postnikov-like
    tower, containing $\ell$ fibrations ($\ell =$~the number of the
    blocks) with fibres
    \[
    \wh\Gamma_t = \prod_{s_t < \alpha \leq s_{t + 1}} \Gamma T \wt
    \xi_\alpha, \qquad t = 1,2,\dots, \ell.
    \]
    \[
    * \longleftarrow \wh \Gamma_\ell \overset{\wh\Gamma_{\ell -1}}{\hbox to
      9mm{\leftarrowfill}} \wh A_2 \overset{\wh\Gamma_{\ell-2}}{\hbox to
      9mm{\leftarrowfill}} \wh A_3 \dots \overset{\wh\Gamma_1}{\hbox to
      9mm{\leftarrowfill}} X_\tau.
    \]
  \end{enumerate}
\end{prop}

\begin{proof}
For any Lie group $G$ the cohomology ring $H^*(BG; \Q)$ has only
even dimen\-si\-o\-nal elements (since this is the invariant part under the
Weyl group action
of $H^*(BT; \Q)$, where $T$ is the maximal torus in $G$, and
$H^*(BT; \Q)$ is the tensor product of $\dim T$ copies of
the ring $H^*(CP^\infty; \Q)$).
It follows that $E^1_{i,j} = H_{i + j}(T\xi_i; \Q)$ can be
non-trivial only if $i + j \equiv c_i \mod 2$.
Hence $E^r_{i,j}$ can be non-trivial only if $i + j \equiv c_i
\mod 2$.
Now $d^r_{i,j} : E^r_{i,j} \to E^r_{i - r, j + r - 1}$ can be
non-trivial only if both groups are non-trivial.
This cannot be, if they are in the same block.

We have proven a) and from this it follows part~b), too.

\noindent
Part c):
Applying the functor $\Gamma$ to the filtration
$$
V_{s_1} \subset V_{s_2} \subset \dots V_{s_l}
$$
we obtain the sequence of fibrations in c).
The fact that the fibres decompose into products follows from the fact that
$V_{s_t}/V_{s_{t+1}}$ and the wedge product
$ T\wt \xi_{s_t+1} \vee \dots \vee~T\wt \xi_{s_{t+1}}$ 
are stably rationally homotopically equivalent.
\end{proof}

\begin{exam}
\label{exa:1b}
Consider Morin maps of odd codimension~$k$, and having at most
$\Sigma^{1_r}$ singularities, i.e. 
$\tau = \{\Sigma^0 < \Sigma^{1,0} < \Sigma^{1,1,0} < \dots
< \Sigma^{1_r,0}\}.$ 

Then $X_\tau {\cong}_\Q \prod_{i \leq r} \Gamma T
\wt\xi_i.$ 
For $i$ odd the space $T \xi_i$ is rationally trivial
(see Lemma~\ref{E1}), hence so is the space $\Gamma T\wt \xi_i.$
For $i = 2j$ we have $\Gamma T\wt \xi_i \cong_\Q \Gamma S^iT(i+1)
\gamma^{SO}_k.$

Hence for any $P^{n + k}$ manifold we have 
\[
\Cob_\tau(P^{n+k}) \otimes \Q \approx \bigoplus_{2i \leq r} \text{\rm
Imm}^{\wt\xi_{2i}} (P^{n+k}) \otimes \Q
\]
where $\wt \xi_{2i}$ is the universal normal bundle in the target of the 
stratum $\Sigma^{1_{2i}}.$
\end{exam}

\begin{exam}

Take for $\tau$ the previous set of singularities for any $r \geq 2$ 
and add
to it as a highest singularity the class of $\eta = III_{2,2}$ (see~\cite{Ma})
which is the simplest stable singularity of type $\Sigma^{2,0}.$
The universal normal bundle $\wt \xi_\eta$ can be found in~\cite{R--Sz},
Theorem 7)

\begin{enumerate}
\item
When the codimension  $k$ of the maps is odd, then  by the previous theorem we have again that  the classifying space 
$X_\tau$ splits into the product
\begin{equation}
\label{spl}
X_\tau = X_{\tau'} \times \Gamma T\wt \xi_\eta,
\end{equation}

where 
$\tau' = \{\leq \Sigma^{1_r}\}$, and 
$X_{\tau'}$ splits into product according to the previous example.
\item

Equality~\ref{spl} holds for $k$ even too if $r = \infty$, because in this
case 
the Kazarian space of all (oriented) Morin maps is rationally homotopy equivalent to $BO(k+1)$ and this has
only even dimensional rational cohomology groups, and the codimension of 
the $III_{(2,2)}$-stratum is even. (The cohomological exact sequence of the
pair $(\Kaz_\tau,\Kaz_{\tau'})$ with rational coefficients splits by dimensional reason.)

\end{enumerate}

\end{exam}

\begin{exam}[Quasi-holomorphic $\tau$-maps]
\label{exa:2b}
(This example is a bit out of the line of the paper.)
Let this time $\tau$ denote a set of
{it holomorphic} singularity classes, i.e. classes of holomorphic germs 
$(C^*, 0) \to C^{*+k},0).$
We call a (non-complex) smooth map $f$ of a smooth manifold into another
one {\it {quasi-holomorphic}} if it imitates the behaviour
of holomorphic maps in the following sense:

The source and target manifolds have smooth stratifications such
that $f$ maps strata into strata, further there exists a
diffeomorphism of the tubular neighbourhood of each stratum onto
a complex vector bundle such that the map $f$ gives a fiberwise
map from the tubular neighbourhoods of the strata in the source
 into those in  the target, 
and the arising mappings
from fibre to fibre are maps from a complex linear space into
another, and all these maps from fibre into fibre
 belong to the classes in $\tau.$
(This time singularity classes are formed using $\mathcal A$-equivalence  
by left-right holomorphic diffeomorphisms.)
The cobordism group of quasi holomorphic $\tau$-maps of
$n$-manifolds into $P^{n + k}$ can be defined and it will be
denoted by $\Cob^C_\tau(P^{n+k})$.
\end{exam}

A classifying space for the cobordisms of quasi holomorphic maps $X^C_\tau$ can be defined analogously
to the real case.

\begin{clai}
$X^C_\tau {\cong}_\Q \prod \Gamma T\wt\xi^C_i$, where
$\wt \xi^C_i$ are the universal complex vector bundles of the
strata in the target.
\end{clai}

\begin{proof}
Since each $\xi_i$ is a complex bundle, we have that each $c_i$
is even.
\end{proof}

\part{On the torsion of the groups $\Cob_\tau(P^{n+k})$}
\label{sec:torsion-part}

In the previous sections we gave a fairly complete computation
of the ranks of the cobordism groups of $\tau$-maps.
Here we present an approach  to the computation of the torsion part.
We shall see that the stable homotopy groups of spheres and the homology rings
of the groups $G_\eta$, \ $\eta\in \tau$ give an upper bound on the torsions of
$\Cob_\tau (n,k).$ The Thom polynomials $Tp_\eta$ evaluated on closed
manifolds will give an estimation of these torsion groups from below.

The stable Hurewicz homomorphism induces isomorphism not only
rationally but also between the $p$-components for sufficiently
high primes.
Arlettaz \cite{Arl} showed the following:
\begin{thm}[Arlettaz]
  \label{Arlettaz}
  For any $(b - 1)$ connected spectrum $X$, the stable Hurewicz
  homomorphism $h_m : \pi_m(X) \to H_m(X; Z)$ has the following properties.
  Let $\varrho_j$ be the exponent of $\pi_j S = \pi_{j + N}(S^N)$,
  $N \gg 1$.
  Then
  \begin{enumerate}[a)]
  \item
    $(\varrho_1 \dots \varrho_{m - b}) (\ker h_m) =
    0, \qquad \quad \forall m \geq b + 1$.
  \item
    $(\varrho_1 \dots \varrho_{m - b - 1}) (\coker h_m) =
    0, \quad \forall m \geq b + 2$.
  \end{enumerate}
\end{thm}
According to Serre \cite{Se} the numbers $\varrho_j$ are not
divisible by a prime $p$ if $p > \frac12 j + 1$.
Let us apply this theorem to (any spectrum defined by) the
virtual space~$T\nu^k$, see Remark~\ref{spect}.
Note that such a spectrum is $(k - 1)$-connected.
Put $b = k$ and $m = n + k$.
Then
\[
h_{n + k} : \pi^s_{n + k} (T\nu^k) \longrightarrow H_{n + k}
(T\nu^k) \approx H_n (\Kaz_\tau)
\]
has the property:
\[
\aligned
\varrho_1 \dots \varrho_n (\ker h_{n + k}) &= 0 \quad \forall n
\geq 1\\
\varrho_1 \dots \varrho_{n - 1} (\coker h_{n + k}) &= 0 \quad
\forall n \geq 2.
\endaligned
\]
Hence for $p > \frac{n}{2} + 1$ the map
$\varphi : \Cob_\tau(n, k) \to H_n(\Kaz_\tau)$ induces an
isomorphism of the $p$-components $(n \geq 2)$.
Analogously for any $(n + k)$-manifold $P^{n+k}$ the map
\[
\varphi : \Cob_\tau(P^{n+k}) \longrightarrow \bigoplus\limits_j H^j
\bigl(P^{n+k}, H_{j - k}(\Kaz_\tau)\bigr)
\]
described in Proposition~\ref{pro:1c} (and after) 
induces isomorphism of $p$-components for $p > \frac{n}{2} + 1$,
$n \geq 2$.
Indeed, Arlettaz' theorem can be reformulated in our case by
saying that the map $\Gamma T\nu^k \to SP\, T\nu^k$ induces a
$p$-homotopy equivalence (i.e. an isomorphism of the $p$-components of the
homotopy groups) until dimension $n + k$ if $p >
\frac{n}{2} + 1$.
Now the equalities
\[
\aligned
{}[P^{n+k}, \Gamma T\nu^k] &= \Cob_\tau(P^{n+k}) \quad \text{ and}\\
[P^{n+k}, SP\, T\nu^k] &= \bigoplus\limits_j H^j\bigl(P^{n+k}, H_{j - k}
(\Kaz_\tau)\bigr)
\endaligned
\] imply the statement.
\hfill $\square$

\begin{prop}
{\it There exists a double spectral sequence
\[
E^1_{p,q,i} = H_p(BG_i ; \pi^s_q) \ \text{ converging to }\
\Cob_\tau (n, k).
\]}
\end{prop}
\begin{expl}
Here --- as before --- $\tau = \{\eta_0 < \eta_1 < \dots < \eta_i
< \dots < \eta_r\}$, $G_i$ is the group $G_{\eta_i}$ (=~maximal
compact subgroup of the automorphism group of the root
of~$[\eta_i]$).
\begin{enumerate}[1)]
\item
  The groups $E^1_{p,q,i}$ for any fixed $i$ form the initial
  term of a spectral sequence, which converges to some groups
  $F_{m,i}$ $(m = p + q)$, and the latter form the initial term of
  spectral sequence that converges to~$\Cob_\tau(n,k)$.
  (Here $n = m+i$ and \ $k$ is fixed.)
\item
  Actually we can get a double spectral sequence from
  $H_p(BG_i; \pi^s_q)$ to $\Cob_\tau(*)$ also in a second way,
  namely we could start by fixing the index~$q$ first.
\end{enumerate}
\end{expl}

\begin{proof}
Let us consider the diagram
\[
\begin{array}{l@{\hspace*{1pt}}l}
H_p(BG_i; \pi^s_q)\
\begin{array}{c}
\approx\\[-1.5mm]
\hbox to10mm{\rightarrowfill}\\[-1.5mm]
\text{\small Thom}
\end{array}  H_{p + c_i} (T\xi_i; \pi^s_q)
\begin{array}{c}
{\text{\small Kazarian}}\\[-1.5mm]
=\!=\!=\!=\!=\!=\!=\!\!>\\
\text{\small in $H_*(\ ; \pi^s_q)$}
\end{array}
  &H_\mu(\Kaz_\tau; \pi^s_q)\quad
(\mu = p + c_i + i)\\[5mm]
\quad \approx \Big\downarrow \text{\rm Thom}
&\quad \approx \Big\downarrow \text{\rm Thom}\\[5mm]
H_{p + c_i + k} (T\wt\xi_i; \pi^s_q)  &H_{\mu + k} (V_\tau;
\pi^s_q) \\[5mm]
\qquad \Big\Downarrow \text{Atiyah--Hirz.}
&\qquad \Big\Downarrow \text{Atiyah--Hirz.} \\[5mm]
\hspace*{5mm}\pi^s_m(T\wt \xi_i) \hspace*{25mm} =\!=\!=\!=\!=\!=\!=\!\!> \quad
\pi^s_*(V_\tau) =  &\pi_*(X_\tau) = \Cob_\tau(*)
\end{array}
\]
$(m = p + c_i + k + q)$.

The two paths from $H_p(BG_i; \pi^s_q)$ to $\Cob_\tau(*)$ give
the two double spectral sequences. The double arrows indicate the
corresponding spectral sequences.
\end{proof}

Although these spectral sequences give a way in principle to
compute the group $\Cob_\tau(n,k)$ completely, the practical
computation seems to be rather difficult.
There are two very special cases when the groups $\Cob_\tau(n,k)$ are
completely computed.

\begin{thm}[\cite{E--Sz--T}]
\label{fold-cobordisms-with-torsion}
Let $\tau$ be $\{\Sigma^0, \Sigma^{1,0}\}$, i.e.\ $\tau$-maps
are fold maps. Then
\[
\Cob^{SO}_\tau(4m - 1, 2m - 1) \approx \Omega_{4m - 1} \oplus \Z_{3t(m)},
\]
where $t(m) = \min \bigl\{ j \mid \alpha_3 (2m + j) \leq 3j
\bigr\}$ and $\alpha_3(x)$ denotes the sum of digits of the
integer $x$ in triadic system.\hfill $\square$
\end{thm}

For the case of fold maps of $2K+2$ manifolds in the Euclidean space
$\BR^{3K+2}$ Terpai gave a complete computation in
  [\cite{TT}].

\section{The (left-right) bordism groups of $\tau$-maps}

By (left-right) bordism of $\tau$-maps we mean the version of
$\tau$-cobordism when not only the source but also the target
manifold of the $\tau$-map can be changed by an arbitrary
oriented cobordism.
We shall denote these groups by $Bord_\tau(n)$.
More precisely the definition is the following.

\begin{defi}\label{LR}
  Two $\tau$-maps $f_i: M_i^n \to P_i^{n+k}$ for $i = 0,\ 1 $
  are (left-right) bordant if
  \begin{enumerate}[\indent a)]
  \item
    there is a compact oriented $V^{n+1}$ manifold with boundary equal to
    
    $\partial V^{n+1} = -M_0^n \cup M_1^n,$
  \item
    there is a compact oriented $W^{n+k+1}$ manifold with boundary equal to
    $\partial W^{n+k+1} = -P_0^{n+k} \cup P_1^{n+k}$, and
  \item
    there is a $\tau$-map $H: V^{n+1}\to W^{n+k+1}$, mapping $M^n_i$ to
    $P^{n+k}_i$ by the map $f_i.$
\end{enumerate}
The set of equivalence classes is denoted by $Bord_\tau(n).$
The disjoint union defines a group operation on this set.
\end{defi}

Here we compute the ranks of these groups.

\begin{thm}
\label{left-right-bordisms}
$Bord_\tau(n) \approx \Omega_{n + k}(X_\tau)$.
\end{thm}

\begin{proof} The proof is completely analogous to the proof of the Pontrjagin
- Thom construction of $\tau$-maps.
\end{proof}

After this proposition the computation of the groups
$Bord_\tau(n)\otimes \Q$ is trivial - as soon as we suppose the rational homology
groups of the Kazarian space to be known. The answer can be formulated in the
simplest way using generating functions.
\footnote{I thank M. Kazarian for teaching me the method of generating
 functions giving a simple way of formulating results.}

\begin{nota}\
  \begin{enumerate}[1)]
  \item
    For any space $A$ we denote by $\mathcal P_A(t)$ the (reduced) rational Poincare\`e series of the space
    $A$, i.e.
    \[\mathcal P_A(t) = \sum_i \text{rk} \wt H_i(A,\Q)\cdot t^i.\]
  \item
    Let us put
    $$ \tau (t) = \sum \text{rk}( Bord_{\tau}(n))\cdot t^{n+k} $$
  \end{enumerate}
\end{nota}

\begin{lem}
Let us denote by $\mathcal {SP}_A(t)$
the rational Poincare\`e series of the infinite symmetric power $SP(A)$ of the
space A. Then $\mathcal{SP}_A(t)$ is the following:
$$
\mathcal {SP} _A(t) = 
\underset {i\ \text {even}} \prod
\Bigl(\frac{1}{1-t^i}\Bigr)^{b_i(A)}\cdot  \underset{i\ \text{odd}} \prod(1+t^i)
^{b_i(A)}
$$

where $b_i(A)$ is the $i$-th rational Betti number of the space $A,$ and $i$
runs over all natural numbers.
\end{lem} 

\begin{proof}
It is well-known, that $SP(X)$ is weakly homotopy equivalent to
$$\underset i \oplus K(H_i(X),i)$$ (see \cite{H}, page 472).
The cohomology ring of the space $K(\Q,i)$ with rational coefficients is
freely generated by an $i$-dimensional multiplicative generator.
By the K\"unneth formula the Poincare\`e series of a product is the product of
those of the factors. All these facts imply the statement easily.
\end{proof}   

Let us put
$$F_\tau(t) = \mathcal {SP}_{S^k(\Kaz_\tau ^+)} =
\underset {i\ \text {even}} \prod
\Bigl(\frac{1}{1-t^i}\Bigr)^{b_{i-k}(\Kaz_\tau)}\cdot  
\underset{i\ \text{odd}} \prod(1+t^i)^{b_{i-k}(\Kaz_\tau)}
$$

\begin{thm}
  $$
  \tau(t) = F_\tau(t) \cdot \mathcal P_{BSO}(t)
  $$
\end{thm}

\begin{proof}

Recall that $X_\tau {\cong}_\Q SP S^k (\Kaz_\tau^+).$

Now using the facts that
$$\Omega_i \otimes \Q \approx H_i(BSO;\Q)
\ \ \text{and} \ \
\Omega_i(A) \approx \bigoplus\limits_{a} H_a(A) \otimes \Omega_{i-a}\otimes \Q$$
we obtain the statement of the Theorem.
\end{proof}

\subsection{Geometric interpretation of the bordism groups of the Kazarian 
spaces.}

\begin{defi}

Let us replace in the previous Definition~\ref{LR} the $\tau$-maps by 
$\ell$-framed $\tau$-maps, and the target manifolds $P_i^{n+k}$ by some
closed oriented $n+k+\ell$-dimensional manifolds $Q_i^{n+k+\ell}.$
The obtained group will be denoted by $Bord_{\tau \oplus \ell}(n).$  
\end{defi}
\begin{prop}
For $\ell > n + k$ the group does not depend on $\ell$ and it is isomorphic to
the $n$-th bordism group of the space $\Kaz_\tau:$
$$Bord_{\tau \oplus\ell}(n) \approx \Omega_n(\Kaz_\tau).$$
The unoriented analogue of this statement holds as well.
\end{prop}
\begin{rema}
This Proposition can be considered as a geometric interpretation of the
homology groups of the Kazarian space because of the formulas connecting
homology groups and bordism groups, see \cite{C--F}.
\end{rema}
\begin{proof}
$$\Omega_n(\Kaz_\tau) \approx \Omega_{n+k}(T\nu^k) \approx
\Omega_{n+k+\ell}(S^\ell T\nu^k) \approx \Omega_{n+k+\ell}(S^\ell V_\tau)
\overset {(*)} \approx $$ $$ \overset {(*)} \approx
\Omega_{n+k+\ell}(\Gamma (S^\ell V_\tau) \approx
\Omega_{n+k+\ell}(B^{\ell}X_\tau)$$
Here $\Gamma$ is the functor $\Omega^\infty S^\infty$.
The isomorphism $(*)$ holds because for any $i$-connected space $A$
the pair $(\Gamma (A), A)$ is $2i$-connected (see \cite{B--E}).
$B^\ell$ means the $\ell$-times ``deloopization'' of the space $X_\tau$, i.e. the space 
$Z^\ell_\tau$, see Remark~\ref{deloop}.
\end{proof}

\section{Final remarks, open problems}

\subsection{Cobordism groups of singular maps as extraordinary cohomology theories.}

We have seen that all the classifying spaces $X_\tau$ and also
all the spaces $\Gamma T\wt\xi_\eta$ are infinite loop spaces.
Therefore they all determine extraordinary cohomology theories. Namely:
For $q \geq 0$ we put
\[
\aligned
h^{-q}_\tau (P^{n+k}) &\overset{\text{\rm def}}{=\!=} \Cob_\tau(P^{n+k}
\times \BR^q)\\
\text{and } \  h^q_\tau(P^{n+k}) &\overset{\text{\rm def}}{=\!=}
\Cob_{\tau \oplus q} (P^{n+k}).
\endaligned
\]
(Recall that $\Cob_{\tau \oplus q}(\ \ )$ denotes
the cobordism group of $\tau$-maps with $q$-framing.
Recall also that we have extended the
functors $\Cob_\tau(\ \ )$ and
$\Cob_{\tau\oplus q}(\ \ )$ from manifolds to
simplicial complexes.)

If $\alpha : P \to P'$ is any continuous map of simplicial complexes, then the induced
homomorphism $\alpha^* = h^q_\tau(\alpha)$ is defined by taking
transverse preimages.
The maps $\theta_{\tau', \tau}: X_{\tau'} \to X_\tau$ or
$\theta_{{\tau},{\eta}} : X_\tau \to \Gamma T\wt\xi_\eta$ when
$\tau'$ is a subset of $\tau$, and $\eta$ is a top singularity
in $\tau$ define {\it cohomology operations}.
For example we have seen that $[f] \in \Cob_\tau(P)$
contains an $\eta$-free map iff the cohomology operation
$\theta_{{\tau},{\eta}}$ vanishes on~$[f]$.
Extending the list of singularities included in $\tau$ we obtain better and better
approximations of the usual cobordism theory
\[
M\, SO^{-k}(P) = \lim\limits_{N \to \infty} \bigl[S^{N + k} P,
M\,SO(N) \bigr].
\]

We have computed all these cohomology theories modulo $p$-torsion for ``small''
primes $p$  --
by expressing them through the homologies of the
corresponding Kazarian space $\Kaz_\tau$, or equivalently through
the $SO$-equivariant cohomologies of the space of some polynomial
maps $\mathcal P_\tau = \lim\limits_{N \to \infty} \mathcal P_\tau(\BR^N,\BR^{N+k}).$ 
(The $SO$-action is one-sided, see the Appendix.)

\subsection{Open problems}

\begin{enumerate}[\bf 1)]
\item
  How to define and compute multiplicative structures on these
  cohomology theories?
\item
  Compute the torsion of $\Cob_\tau(P^{n+k})$, (see section~\ref{sec:torsion-part})
\item
  Extend all the results of this
  paper to unoriented cobordism groups.
\item
  Find analogues of exact sequences of Rohlin, Wall, Atiyah
  relating the oriented and unoriented versions of these groups (see theorems
  4.2 and 4.3 in \cite{At}).
\end{enumerate}

\section{Appendix: The stabilisation of the Kazarian space}

First we recall the Borel construction:

If $V$ is a $G$-space (i.e.\  a $G$-action $G \times V \to V$ is
given on it), then the Borel construction of $V$ is the space
$BV = EG \underset{G}{\times} V$.
If $\Xi \subset V$ is a $G$-invariant subspace, then we have
$B\Xi \subset BV$.

If $V = \bigcup\limits_\eta \Sigma_\eta$ is the decomposition of
$V$ into $G$-orbits (parameterised by $\eta$), then $BV =
\bigcup\limits_\eta B \Sigma_\eta$.

\begin{lem}[Kazarian, \cite{K1}]
If $V$ is contractible, and $G_x$ is the stabiliser of a point
$x \in V$, and $\Sigma = G(x)$ is the orbit of $x$, then
$B\Sigma = BG_x$, where $BG_x$ is the classifying space of the
group~$G_x$.
\end{lem}

\begin{proof}
$B \Sigma = EG \underset{G}{\times} \Sigma = EG
\underset{G}{\times} G / G_x = (EG) / G_x = BG_x$.
\end{proof}

\begin{nota}
By $\Kaz^{\prime}_{\tau}(n)$ we denoted the
``unstable'' Kazarian space which was a subset in $ J^K(\gamma^{SO}_n,
\gamma^{SO}_{n + k})$.
Namely: the fibre of the bundle $$J^K(\gamma^{SO}_n, \gamma^{SO}_{n +
k}) \longrightarrow B\, SO_n \times B\, SO_{n + k}$$ is
$J^K_0(\BR^n, \BR^{n + k})$ of the space of all polynomial maps $\BR^n
\to \BR^{n + k}$ of degree $\leq K$ and mapping the origin of
$\BR^n$ into that of $\BR^{n + k}$.
Let $V_\tau(n) \subset J^K_0(\BR^n, \BR^{n + k})$ consist
of the maps having a singularity at the origin from the list~$\tau$.

The union of these subsets $V_\tau(n)$ in each fibre
forms the space $ \Kaz'_\tau(n)$.
\end{nota}

Applying the previous Lemma we would get the decomposition
$\Kaz'_\tau(n) = \bigcup BG_x$, where $G_x$ is the stabiliser group of 
$x\in V_\tau(n)$ under the action of $G = SO_n \times SO_{n+k}.$
There is a small problem here. Namely that the stabiliser groups are not direct
products.  Namely if $x$ is representing a singularity class $[\eta]$, then the
stabiliser $G_x$ is an index $2$ subgroup of $G_\eta^O \times O(n-c(\eta))$,
(recall that $c(\eta)$ is the codimension of the $\eta$-stratum in the source). 
The same group $G_x$ can be described also as a $\Z_2$ extension of the group $G_\eta \times
SO(n-c(\eta)).$

In order to obtain the stabiliser subgroups as direct
products $G_\eta \times O(n-c(\eta))$, we apply the following trick.
We define a modification of the unstable Kazarian space $\Kaz'_\tau(n).$
Let us define the subgroup $H_n$ of index $2$ in $O_n \times O_{n+k}$ as
follows:
$$H_n = \{(A,B)| A\in O_n,\ B \in O_{n+k},\ \text{det}(A) \cdot \text{det}(B) >0\}.$$
Put
$$\Kaz^{co}_\tau(n) = V_\tau(n) \times _{H_n} EH_n.$$
(This is a version of the unstable Kazarian space $\Kaz'_\tau(n)$ for the Co-oriented
maps, i.e. for the maps having orientations on their virtual normal bundles, but
not necessarily both on the source and on the target.)

\begin{lem}
\label{lem:2}
$\Kaz^{co}_\tau(n) = \bigcup\limits_{\eta \in
\tau} B G_\eta \times BO(n - c(\eta))$.
\end{lem}

\begin{proof}
Let $\eta$ be a singularity class from $\tau$ and
let us denote by $\eta^0$ the root of $\eta$.
Recall that $\eta^0$ is a germ $(\BR^{c(\eta)}, 0) \to (\BR^{c(\eta)
+ k}, 0)$, and $\eta$ is locally equivalent to $\eta^0 \times
\text{\rm id}_{\BR^{n - c(\eta)}}$.
Then the automorphism group of $\eta$ in $H_n$ is $G_\eta \times O(n - c(\eta))$.
Hence if $x \in V_\tau(n)$ represents $\eta$ then the stabiliser is
$G^{(n)}_x = G_\eta
\times O(n - c(\eta))$. Now apply the previous Lemma.
\end{proof}

\begin{defi}
The (stable) Kazarian space is defined as:
\[
\Kaz_\tau \overset{\text{\rm def}}{=\!=} \lim\limits_{N \to \infty}
E\, SO(N + k) \underset{SO(N + k)}{\times} V_\tau (N)
\subset \lim\limits_{N \to \infty} J^K (\varepsilon^N, \gamma_{N
+ k})
\]
\end{defi}

Note that here the action of the group $SO(N+k)$ is one sided, we do not act on
the source space $\BR^N$  of the polynomial maps.

\begin{prop}
The (stable) Kazarian space
$\Kaz_\tau$
decomposes into the disjoint union of the base spaces $BG_\eta$
for $\eta \in \tau$.
(Here $G_\eta$ is the maximal compact subgroup of the
automorphism group of the roof of~$\eta$.)\hfill $\square$
\end{prop}

\begin{proof}

The double cover 
$B\,SO_n \times B\,SO_{n+k} \to BH_n$ induces a double cover
$\Kaz'_\tau(n) \to 
\Kaz^{co}_\tau(n).$
Let us denote by $\Kaz_\tau(n)$ the lift of the bundle 
$\Kaz'_\tau(n) \to B\,SO_n \times B\,SO_{n+k}$ to the 
space $ESO_n\times B\,SO(n+k) \cong BSO_{n+k}.$
Clearly the map $\Kaz_\tau(n) \to \Kaz'_\tau(n)$ is an $SO_n$-bundle.
The composition of this bundle map with the double cover 
$\Kaz'_\tau(n) \to 
\Kaz^{co}_\tau(n)$ is an $O_n$-bundle.
The part of this $O_n$ bundle over the stratum $BG_\eta \times BO(n-c)$
is the space of cosets $O_n/(O_{n-c} \times G_\eta)$, i.e the quotient of the
Stiefel manifold $V_c(R^n)$ by a free $G_\eta$ action.
Since the Stiefel manifold is $n-c-1$-connected the limit for $n \to \infty$
will be the classifying space  $BG_\eta.$ 
Hence the corresponding stratum in the stable Kazarian space $\Kaz_\tau$
will be $BG_\eta.$
\end{proof}

\begin {rema}
It follows from the construction 
of the stable Kazarian space $\Kaz_\tau$
that any (proper) $\tau$-map defines a
homotopically unique map of its source to the space $\Kaz_\tau.$
\end{rema}

\Addresses

\end{document}